\documentclass[11pt,reqno]{amsart} 
\usepackage{amsfonts, amsmath, amssymb, amsthm, color}
\usepackage[margin=1in]{geometry}
\usepackage[dvipsnames]{xcolor}
\usepackage[colorlinks=true]{hyperref}

\newtheorem{thm}{Theorem}[section]
\newtheorem{cor}[thm]{Corollary}
\newtheorem{prop}[thm]{Proposition}
\newtheorem{lemma}[thm]{Lemma}
\newtheorem{assum}[thm]{Assumption}

\theoremstyle{definition}
\newtheorem{definition}[thm]{Definition}
\newtheorem{remark}[thm]{Remark}

\newcommand{\N}{\mathbb{N}}
\newcommand{\R}{\mathbb{R}}
\renewcommand{\S}{\mathbb{S}}
\newcommand{\mbd}{\mathbf{d}}
\newcommand{\mbx}{\mathbf{x}}

\newcommand{\mca}{\mathcal{A}}

\newcommand{\mcr}{\mathcal{R}}
\newcommand{\mci}{\mathcal{I}}
\newcommand{\mcu}{\mathcal{U}}
\newcommand{\mcv}{\mathcal{V}}

\newcommand{\by}{\bar{y}}

\newcommand{\hps}{\hat{\psi}}
\newcommand{\tps}{\tilde{\psi}}
\newcommand{\whg}{\widehat{G}}
\newcommand{\whh}{\widehat{H}}

\newcommand{\td}{\tilde{d}}
\newcommand{\tr}{\tilde{r}}
\newcommand{\tx}{\tilde{x}}
\newcommand{\tga}{\tilde{\gamma}}
\newcommand{\tla}{\tilde{\lambda}}
\newcommand{\trh}{\tilde{\rho}}
\newcommand{\tta}{\tilde{\tau}}
\newcommand{\wtg}{\widetilde{G}}
\newcommand{\wth}{\widetilde{H}}
\newcommand{\wts}{\widetilde{S}}
\newcommand{\wtmbd}{\widetilde{\mbd}}

\newcommand{\pa}{\partial}
\newcommand{\ep}{\epsilon}
\newcommand{\ga}{\gamma}
\newcommand{\e}{\varepsilon}
\newcommand{\vrh}{\varrho}

\newcommand{\vsi}{\varsigma}

\newcommand{\bsde}{\boldsymbol{\delta}}
\newcommand{\bsmu}{\boldsymbol{\mu}}
\newcommand{\bsxi}{\boldsymbol{\xi}}

\newcommand{\dexi}{\bsde,\bsxi}
\newcommand{\tdx}{\widetilde{\mbd},\mbx}
\newcommand{\mux}{\bsmu,\mbx}

\newcommand{\ovh}{\overline{H}}
\newcommand{\ovom}{\overline{\Omega}}

\newcommand{\dist}{\textnormal{dist}}

\renewcommand{\(}{\left(}
\renewcommand{\)}{\right)}

\allowdisplaybreaks
\numberwithin{equation}{section}

\makeatletter
\@namedef{subjclassname@2020}{%
\textup{2020} Mathematics Subject Classification}
\makeatother

\begin{document}
\title[Positive solutions of the Lane-Emden system]{Asymptotic analysis on positive solutions of the \\ Lane-Emden system with nearly critical exponents}

\author{Seunghyeok Kim}
\address[Seunghyeok Kim]{Department of Mathematics and Research Institute for Natural Sciences, College of Natural Sciences, Hanyang University,
222 Wangsimni-ro Seongdong-gu, Seoul 04763, Republic of Korea}
\email{shkim0401@hanyang.ac.kr shkim0401@gmail.com}

\author{Sang-Hyuck Moon}
\address[Sang-Hyuck Moon]{National Center for Theoretical Sciences, National Taiwan University, Taipei 106, Taiwan}
\email{shmoon@ncts.ntu.edu.tw mshyeog@gmail.com}

\begin{abstract}
We concern a family $\{(u_{\e},v_{\e})\}_{\e > 0}$ of solutions of the Lane-Emden system on a smooth bounded convex domain $\Omega$ in $\R^N$
\[\begin{cases}
-\Delta u_{\e} = v_{\e}^p &\text{in } \Omega,\\
-\Delta v_{\e} = u_{\e}^{q_{\e}} &\text{in } \Omega,\\
u_{\e},\, v_{\e} > 0 &\text{in } \Omega,\\
u_{\e} = v_{\e} =0 &\text{on } \pa\Omega
\end{cases}\]
for $N \ge 4$, $\max\{1,\frac{3}{N-2}\} < p < q_{\e}$ and small
\[\e := \frac{N}{p+1} + \frac{N}{q_{\e}+1} - (N-2) > 0.\]
This system appears as the extremal equation of the Sobolev embedding $W^{2,(p+1)/p}(\Omega) \hookrightarrow L^{q_{\e}+1}(\Omega)$,
and is also closely related to the Calder\'on-Zygmund estimate. Under the a natural energy condition
\[\sup_{\e > 0} \(\|u_{\e}\|_{W^{2,{p+1 \over p}}(\Omega)} + \|v_{\e}\|_{W^{2,{q_{\e}+1 \over q_{\e}}}(\Omega)}\) < \infty,\]
we prove that the multiple bubbling phenomena may arise for the family $\{(u_{\e},v_{\e})\}_{\e > 0}$, and establish a detailed qualitative and quantitative description.
If $p < \frac{N}{N-2}$, the nonlinear structure of the system makes the interaction between bubbles so strong,
so the determination process of the blow-up rates and locations is completely different from that of the classical Lane-Emden equation.
If $p \ge \frac{N}{N-2}$, the blow-up scenario is relatively close to (but not the same as) that of the classical Lane-Emden equation, and only one-bubble solutions can exist.
Even in the latter case, the standard approach does not work well, which forces us to devise a new method.
Using our analysis, we also deduce a general existence theorem valid on any smooth bounded domains.
\end{abstract}

\date{\today}
\subjclass[2020]{Primary: 35J47, Secondary: 35B09, 35B33, 35B40, 35B44, 35B45}
\keywords{Lane-Emden system, critical hyperbola, positive solutions, asymptotic analysis, multi-bubbles, pointwise estimates}
\maketitle

\section{Introduction}
In this paper, we concern the Lane-Emden system, one of the simplest Hamiltonian-type elliptic systems. It is given as
\begin{equation}\label{eq:LE}
\begin{cases}
-\Delta u = v^p &\text{in } \Omega,\\
-\Delta v = u^q &\text{in } \Omega,\\
u,\, v > 0 &\text{in } \Omega,\\
u = v =0 &\text{on } \pa\Omega
\end{cases}
\end{equation}
where $\Omega$ is a smooth bounded domain in $\R^N$ for $N \ge 3$, $p \in (\frac{2}{N-2}, \frac{N+2}{N-2}]$, and $q \ge p$.
As the classical (scalar) Lane-Emden equation
\begin{equation}\label{eq:LEscalar}
\begin{cases}
-\Delta u = u^p &\text{in } \Omega,\\
u,\, v > 0 &\text{in } \Omega,\\
u = v = 0 &\text{on } \pa\Omega
\end{cases}
\end{equation}
for $p \in [1, \frac{N+2}{N-2}]$ is the extremal equation of the Sobolev embedding $W^{1,2}_0(\Omega) \hookrightarrow L^{p+1}(\Omega)$,
system \eqref{eq:LE} is that of the embedding $W^{2,(p+1)/p}(\Omega) \hookrightarrow L^{q+1}(\Omega)$ for a pair $(p,q)$ in the subcritical region
\begin{equation}\label{eq:cr-1}
\frac{1}{p+1} + \frac{1}{q+1} > \frac{N-2}{N},
\end{equation}
or on the (Sobolev) critical hyperbola\footnote{The critical hyperbola, introduced by Mitidieri \cite{Mi} and van der Vorst \cite{V}
(and soon treated by several researchers including Cl\'ement et al. \cite{CFM}),
has an essential role in determining the solution structure of \eqref{eq:LE}.}
\begin{equation}\label{eq:cr-hy}
\frac{1}{p+1} + \frac{1}{q+1} = \frac{N-2}{N}.
\end{equation}
System \eqref{eq:LE} is also closely related to the Calder\'on-Zygmund estimate in the sense that its solutions satisfy a circular relation
\begin{align*}
v \in L^{p+1}(\Omega) &\Rightarrow -\Delta u = v^p \in L^{p+1 \over p}(\Omega) \Rightarrow u \in W^{2,{p+1 \over p}}(\Omega) \subset L^{q+1}(\Omega) \\
&\Rightarrow -\Delta v = u^q \in L^{q+1 \over q}(\Omega) \Rightarrow v \in W^{2,{q+1 \over q}}(\Omega) \subset L^{p+1}(\Omega)
\end{align*}
provided \eqref{eq:cr-1} or \eqref{eq:cr-hy}.

\medskip
As the natural counterpart of \eqref{eq:LEscalar}, system \eqref{eq:LE} has received considerable attention for decades, and now abundant results on the existence theorem are available in the literature.
However, as far as the authors know, only a few works have examined the qualitative properties of the system, such as symmetry, concentration, or multiplicity of solutions.
The reader is advised to check the survey paper \cite{BMT}, which furnishes vast treatments on \eqref{eq:LE} known before the year 2014.

Our objective is to contribute to the latter direction by examining asymptotic behavior of solutions of \eqref{eq:LE}
provided that the pair $(p,q)$ satisfies $p \le q$, and \eqref{eq:cr-1}, and is close to the critical hyperbola.

\subsection{History and motivations of the problem}
A solution of \eqref{eq:LE} is realized as a positive critical point of the energy functional
\begin{equation}\label{eq:energy}
I_{p,q}(u,v) = \int_{\Omega} \nabla u \cdot \nabla v \, dx - \frac{1}{p+1} \int_{\Omega} |v|^{p+1} dx - \frac{1}{q+1} \int_{\Omega} |u|^{q+1} dx
\end{equation}
for $(u, v) \in W_0^{1,s}(\Omega) \times W_0^{1,\frac{s}{s-1}}(\Omega)$ for a suitable $s > 1$.
Exploiting the strong indefiniteness of $E$, Hulshof and van der Vorst \cite{HV} and de Figueiredo and Felmer \cite{dFF}
proved independently that if $p, q > 1$ and \eqref{eq:cr-1} holds, then \eqref{eq:LE} has a solution.
Moreover, by studying the associated minimization problem to the scalar equation, which is an equivalent formulation of \eqref{eq:LE},
\begin{equation}\label{eq:LEs}
\begin{cases}
-\Delta \((-\Delta u)^{1 \over p}\) = u^q &\text{in } \Omega,\\
u,\, -\Delta u > 0 &\text{in } \Omega,\\
u = \Delta u = 0 &\text{on } \pa \Omega,
\end{cases}
\end{equation}
Bonheure et al. \cite{BMR} extended the aforementioned existence results to cover any $p, q > 0$ satisfying $pq \ne 1$ and \eqref{eq:cr-1}.
On the other hand, a Pohozaev-type argument (refer to Mitidieri \cite{Mi}) yields that if $\Omega$ is star-shaped and
\[\frac{1}{p+1} + \frac{1}{q+1} \le \frac{N-2}{N},\]
then \eqref{eq:LE} has no solution.

\medskip
As shown in \cite{BMR}, \eqref{eq:LE} has a least energy solution (or a ground state), namely, a solution that attains the least value of $E$ among all nontrivial solutions, provided $pq \ne 1$ and \eqref{eq:cr-1}.
In fact, $u$ is a least energy solution of \eqref{eq:LEs} if and only if $(u, ((-\Delta) u)^{1/p})$ is a least energy solution of \eqref{eq:LE}.
In \cite{G}, Guerra investigated the asymptotic profile of least energy solutions of \eqref{eq:LE}, assuming that $\Omega$ is convex and $p \ge 1$.
Later, Choi and the first author of this paper \cite{CK} extended his results by covering $p < 1$ and arbitrary smooth bounded domains.
Their works are natural generalizations of Han \cite{H} and Rey \cite{R} that studied least energy solutions of the scalar equation \eqref{eq:LEscalar}, and summarized as follows:
Fix a number $p \in (\frac{2}{N-2}, \frac{N+2}{N-2}]$, and define $q_{\e}$ by the relation
\begin{equation}\label{eq:cr-2}
\frac{1}{p+1} + \frac{1}{q_{\e}+1} = \frac{N-2+\e}{N}
\end{equation}
where $\e > 0$ is small. Clearly, $(p,q_{\e})$ is in the subcritical region \eqref{eq:cr-1} and approaches the critical hyperbola \eqref{eq:cr-hy} as $\e \to 0$.
Let $q_0$ be the limit of $q_{\e}$ as $\e \to 0$ so that $(p,q_0)$ satisfies \eqref{eq:cr-hy} and $p \le q_0$.
If $\{(u_{\e},v_{\e})\}_{\e \in (0,\e_0)}$ is a family of least energy solutions of \eqref{eq:LE} with $q = q_{\e}$, then it blows up at an interior point $\xi_0$ of $\Omega$ as $\e \to 0$.
In other words, there exists a family $\{x_{\e}\}_{\e \in (0,\e_0)} \subset \Omega$ such that
\begin{equation}\label{eq:blowup}
u_{\e}(x_{\e}) = \|u_{\e}\|_{L^{\infty}(\Omega)} \to \infty, \quad x_{\e} \to \xi_0 \in \Omega,
\quad \text{and} \quad
u_{\e},\, v_{\e} \to 0 \quad \text{in } C^1_{\textnormal{loc}}(\Omega \setminus \{\xi_0\})
\end{equation}
as $\e \to 0$, after passing to a subsequence. Also, if we define
\begin{equation}\label{eq:albe}
\alpha_{\e} = \frac{2(p+1)}{pq_{\e}-1}, \quad \beta_{\e} = \frac{2(q_{\e}+1)}{pq_{\e}-1}, \quad \text{and} \quad \lambda_{\e} = u_{\e}^{1 \over \alpha_{\e}}({x_{\e}}),
\end{equation}
then, along a subsequence,
\[\(\lambda_{\e}^{-\alpha_{\e}} u_{\e}(\lambda_{\e}^{-1} \cdot + x_{\e}), \lambda_{\e}^{-\beta_{\e}} v_{\e}(\lambda_{\e}^{-1} \cdot + x_{\e})\) \to (U_{1,0},V_{1,0}) \quad \text{in } C^2_{\textnormal{loc}}(\R^N) \quad \text{as } \e \to 0.\]
Here, $(U_{1,0},V_{1,0})$ is the unique least energy solution of
\begin{equation}\label{eq:LERn}
\begin{cases}
-\Delta U_{1,0} = V_{1,0}^p,\ -\Delta V_{1,0} = U_{1,0}^{q_0} \quad \text{in } \R^N,\\
U_{1,0}, V_{1,0} > 0 \hspace{96pt} \text{in } \R^N,\\
(U_{1,0},V_{1,0}) \in \dot{W}^{2,\frac{p+1}{p}}(\R^N) \times \dot{W}^{2,\frac{q_0+1}{q_0}}(\R^N)
\end{cases}
\end{equation}
such that
\begin{equation}\label{eq:LERn2}
U_{1,0}(0) = 1 = \max_{x \in \R^N} U_{1,0}(x).\footnotemark
\end{equation}
\footnotetext{The unique existence of $(U_{1,0},V_{1,0})$ was proved by Hulshof and Van der Vorst \cite{HV2}. Throughout the paper, we call $(U_{1,0},V_{1,0})$ the standard bubble.}Let
$G$ and $\tau$ be the Green's function and the Robin function of the Dirichlet Laplacian $-\Delta$ in $\Omega$, respectively.
For $p \in (\frac{2}{N-2}, \frac{N}{N-2})$ and $\xi \in \Omega$, let $\wtg(\cdot,\xi)$ be the unique solution of
\begin{equation}\label{eq:wtg}
\begin{cases}
-\Delta_x \wtg(x,\xi) = G^p(x,\xi) &\text{for } x \in \Omega,\\
\wtg(x,\xi) = 0 &\text{for } x \in \pa\Omega,
\end{cases}
\end{equation}
$\wth$ be the $C^1$-regular part of $\wtg$, and $\tta(x) = \wth(x,x)$ for $x \in \Omega$; see Subsection \ref{subsec:G}.
Then $\xi_0$ must be a critical point of the Robin function $\tau$ for $p \in [\frac{N}{N-2}, \frac{N+2}{N-2}]$,
and that of $\tta$ for $p \in (\frac{2}{N-2}, \frac{N}{N-2})$.
Besides, the asymptotic behavior of $\|u_{\e}\|_{L^{\infty}(\Omega)}$, and a more exact asymptotic profile of $\{(u_{\e},v_{\e})\}_{\e \in (0,\e_0)}$ in $\Omega \setminus \{\xi_0\}$ than one in \eqref{eq:blowup}
can be given in terms of the quantities $N$, $(p,q_0)$, $\e$, $(U_{1,0},V_{1,0})$, $G(\cdot,\xi_0)$, $\tau(\xi_0)$, $\wtg(\cdot,\xi_0)$, and $\tta(\xi_0)$; refer to \cite[Theorem 1.1]{G}, \cite[Theorem B]{CK}, and \cite[Lemma 2.2]{KP}.
Note that if $p > \frac{N+2}{N-2}$, one can exchange the roles of $p$ and $q$.

\medskip
Recently, the first author and Pistoia \cite{KP} built multi-bubble solutions of \eqref{eq:LE}
by employing the non-degeneracy result on the standard bubble $(U_{1,0},V_{1,0})$ due to Frank and themselves \cite{FKP}; refer to Lemma \ref{lemma:FKP} below.
Among others, they proved the following result: Assume that $N \ge 4$, $p \in (1, \frac{N-1}{N-2})$, and $(p,q) = (p,q_{\e})$ satisfies \eqref{eq:cr-2}.
Suppose also that $\Omega$ is a dumbbell-shaped domain with $l-1$ necks for some $l \in \N$.\footnote{The precise definition of the dumbbell-shaped domain can be found in \cite[Section 1]{KP}.}
Given any $k \in \{1, \ldots, l\}$, there exists a small number $\e_0 > 0$ such that
for each for all $\e \in (0,\e_0)$, \eqref{eq:LE} has $N \choose k$ solutions that blow up at $k$ points as $\e \to 0$.

In the proof, they found a more exact asymptotic profile of the solutions, including their blow-up rate and location.
Also, while they focused the case $p \in (1, \frac{N-1}{N-2})$, their method can also cover $p \in [\frac{N-1}{N-2}, \frac{N+2}{N-2})$ as our analysis will indicate. See Theorem \ref{thm:mainby} below.

\medskip
The result of \cite{KP} implies the existence of higher-energy solutions (or excited states) of \eqref{eq:LE} in some instances.
Hence it is natural to ask how they behave as $(p,q)$ approaches the critical hyperbola.
The primary purpose of this paper is to give a complete description of the asymptotic behavior of an arbitrary family of solutions provided that they satisfy a natural energy condition, $\Omega$ is convex, and $p > \max\{1,\frac{3}{N-2}\}$.

More precisely, we will see that the multiple bubbling phenomena may arise, while the cluster or tower phenomena cannot.
Furthermore, performing a fine analysis, we will show that the asymptotic profile of the solutions depends on the value of $p$.
If $p$ is less than the Serrin exponent $\frac{N}{N-2}$, the interaction between bubbles turns out to be very strong.
We have to reflect this fact in determining the blow-up rates and locations, so the analysis becomes difficult.
If $p \ge \frac{N}{N-2}$, the blow-up scenario is relatively close to that of the classical Lane-Emden equation.
Nonetheless, we cannot follow the standard approach, because it only provides rough estimates that are unusable for $p$ close to $\frac{N}{N-2}$.
We will devise a new method to get over this technical issue.

Using our analysis, we will also obtain a general existence theorem that holds on arbitrary smooth bounded domains.

\subsection{Main theorems}
We now state the main theorems. The following is our general assumption for Theorems \ref{thm:main1}--\ref{thm:main3}.
\begin{assum}\label{assum:main}
We assume that $N \ge 4$, $\Omega$ is a smooth bounded convex domain in $\R^N$, $p \in (\max\{1,\frac{3}{N-2}\}, \frac{N+2}{N-2})$, $\e_0 > 0$ is a small number,
and $q_{\e}$ is the number determined by \eqref{eq:cr-2} for each $\e \in (0,\e_0)$.
Also, let $\{(u_{\e},v_{\e})\}_{\e \in (0,\e_0)}$ be a family of solutions of \eqref{eq:LE} with $q = q_{\e}$ satisfying a natural energy condition
\begin{equation}\label{eq:main0}
\limsup_{\e \to 0} \(\|u_{\e}\|_{W^{2,{p+1 \over p}}(\Omega)} + \|v_{\e}\|_{W^{2,{q_{\e}+1 \over q_{\e}}}(\Omega)}\) \le C
\end{equation}
for some constant $C > 0$.
\end{assum}

The first theorem describes the multi-bubble phenomenon of \eqref{eq:LE}.
\begin{thm}\label{thm:main1}
Suppose that Assumption \ref{assum:main} holds. Then, after passing to a subsequence, either
\begin{equation}\label{eq:uvconv}
(u_{\e},v_{\e}) \to (u_0,v_0) \quad \text{in } (C^2(\ovom))^2 \quad \text{as } \e \to 0
\end{equation}
where $(u_0,v_0)$ is a solution of \eqref{eq:LE} with $q = q_0$, or the followings hold:

\medskip \noindent \textnormal{(1)} $(u_{\e},v_{\e}) \rightharpoonup (0,0)$ weakly and not strongly in $W^{2,(p+1)/p}(\Omega) \times W^{2,(q_0+1)/q_0}(\Omega)$ as $\e \to 0$.

\medskip \noindent \textnormal{(2)} The family $\{(u_{\e},v_{\e})\}_{\e \in (0,\e_0)}$ blows up at $k \in \N$ distinct points $\xi_1, \ldots, \xi_k \in \Omega$ (apart from the boundary) as $\e \to 0$.
In other words, there exist $k$ families $\{x_{i\e}\}_{\e \in (0,\e_0)} \subset \Omega$ for $i = 1, \ldots, k$ such that
\begin{equation}\label{eq:uvconv0}
u_{\e}(x_{i\e}) \to \infty, \quad x_{i\e} \to \xi_i \in \Omega,
\quad \text{and} \quad
u_{\e},\, v_{\e} \to 0 \quad \text{in } C^1_{\textnormal{loc}}(\overline{\Omega} \setminus \{\xi_1, \ldots, \xi_k\})
\end{equation}
as $\e \to 0$.

\medskip \noindent \textnormal{(3)} Let $\alpha_{\e}$ and $\beta_{\e}$ be the numbers in \eqref{eq:albe}.
For a fixed index $i = 1, \ldots, k$, we set $\lambda_{i\e} = u_{\e}^{1/\alpha_{\e}}({x_{i\e}})$. Then
\begin{equation}\label{eq:uvconv1}
\(\lambda_{i\e}^{-\alpha_{\e}} u_{\e}(\lambda_{i\e}^{-1} \cdot + x_{i\e}), \lambda_{i\e}^{-\beta_{\e}} v_{\e}(\lambda_{i\e}^{-1} \cdot + x_{i\e})\) \to (U_{1,0},V_{1,0}) \quad \text{in } C^2_{\textnormal{loc}}(\R^N) \quad \text{as } \e \to 0
\end{equation}
where $(U_{1,0},V_{1,0})$ is the standard bubble. In addition, there exists a constant $C > 0$ such that
\begin{equation}\label{eq:uvconv11}
\limsup_{\e \to 0} \frac{\lambda_{i\e}}{\lambda_{j\e}} \le C \quad \text{for any } 1 \le i \ne j \le k.
\end{equation}
\end{thm}

The second theorem depicts the exact asymptotic profile of solutions to \eqref{eq:LE} as $(p,q)$ approaches the critical hyperbola provided $p < \frac{N}{N-2}$.
As can be seen, the determination process of the blow-up rates and locations are much more cumbersome than that of the scalar equation \eqref{eq:LEscalar} or the case $p \ge \frac{N}{N-2}$.
\begin{thm}\label{thm:main2}
Suppose that Assumption \ref{assum:main} holds, $p < \frac{N}{N-2}$, and Theorem \ref{thm:main1} (1) is valid.
By setting $d_{i\e} = \frac{\lambda_{1\e}}{\lambda_{i\e}}$ and passing to a subsequence, we may assume that
\[\delta_i = \lim_{\e \to 0} d_{i\e} \in (C^{-1},C) \quad \text{for some } C > 1.\]
Writing $(\dexi) = (\delta_1, \ldots, \delta_k,\xi_1, \ldots, \xi_k) \in (0,\infty)^k \times \Omega^k$, let $\wtg_{\dexi}: \Omega \to \R$ be the solution of
\begin{equation}\label{eq:wtgdx}
\begin{cases}
\displaystyle -\Delta \wtg_{\dexi} = \(\sum_{i=1}^k \delta_i^{\frac{N}{q_0+1}} G(\cdot,\xi_i)\)^p &\text{in } \Omega,\\
\wtg_{\dexi} = 0 &\text{on } \pa\Omega.
\end{cases}
\end{equation}
For $i = 1, \ldots, k$, let also $\wth_{\dexi,i}: \Omega \to \R$ be the local $C^1$-regularization of $\wtg_{\dexi}$ around $\xi_i$ whose precise definition is given in \eqref{eq:wthdxi}. Then the followings hold:

\medskip \noindent \textnormal{(1)} The $i$-th blow-up point $\xi_i \in \Omega$ depicted in Theorem \ref{thm:main1} (2) is a critical point of $\wth_{\dexi,i}$.

\medskip \noindent \textnormal{(2)} We have
\[\lim_{\e \to 0}\e u_{\e}^{p+1}(x_{i\e}) = \frac{N}{q_0+1} S^{\frac{p(q_0+1)}{1-pq_0}} \|U_{1,0}\|_{L^{q_0}(\R^N)}^{(p+1)q_0} \delta_i^{-\frac{Np}{q_0+1}} \wth_{\dexi,i}(\xi_i)\]
for all $i = 1, \ldots, k$, where
\begin{equation}\label{eq:S}
S := \inf_{U \in W^{2,\frac{p+1}{p}}(\R^N) \setminus \{0\}} \frac{\int_{\R^N} |\Delta U|^{\frac{p+1}{p}}}{\|U\|_{L^{q_0+1}(\R^N)}^{\frac{p+1}{p}}}
= \frac{\int_{\R^N} |\Delta U_{1,0}|^{\frac{p+1}{p}}}{\|U_{1,0}\|_{L^{q_0+1}}^{\frac{p+1}{p}}(\R^N)} = \|U_{1,0}\|_{L^{q_0+1}(\R^N)}^{\frac{pq_0-1}{p}}.
\end{equation}
In particular, $\delta_i^{N/(q_0+1)} \wth_{\dexi,i}(\xi_i) = \delta_j^{N/(q_0+1)} \wth_{\dexi,j}(\xi_j)$ for all $i, j = 1, \ldots, k$.

\medskip \noindent \textnormal{(3)} Without loss of generality, let us assume that $\lambda_{1\e}^{\alpha_{\e}} = \|u_{\e}\|_{L^\infty(\Omega)}$. Then
\[\lim_{\e \to 0} \|u_{\e}\|_{L^{\infty}(\Omega)} v_{\e}(x) = \|U_{1,0}\|_{L^{q_0}(\R^N)}^{q_0} \sum_{i=1}^k \delta_i^{\frac{N}{q_0+1}} G(x,\xi_i)\]
and
\[\lim_{\e \to 0} \|u_{\e}\|_{L^{\infty}(\Omega)}^p u_{\e}(x) = \|U_{1,0}\|_{L^{q_0}(\R^N)}^{pq_0} \wtg_{\dexi}(x)\]
in $C^1_{\textnormal{loc}}(\Omega \setminus \{\xi_1, \ldots, \xi_k\})$-sense.
\end{thm}

The third theorem shows that if $p \ge \frac{N}{N-2}$, then \eqref{eq:LE} behaves similarly to \eqref{eq:LEscalar}, and only one-bubble solutions exist for $\Omega$ convex.
\begin{thm}\label{thm:main3}
Suppose that Assumption \ref{assum:main} holds, $p \ge \frac{N}{N-2}$, and Theorem \ref{thm:main1} (1) is valid.
Then $\{(u_{\e},v_{\e})\}_{\e \in (0,\e_0)}$ blows up at only one point $\xi_0 \in \Omega$ and
\[\frac{\int_{\Omega} |\Delta u_{\e}|^{\frac{p+1}{p}}}{\|u_{\e}\|_{L^{q_{\e}+1}(\Omega)}^{p+1 \over p}} \to S \quad \text{as } \e \to 0\]
and $S > 0$ is the constant in \eqref{eq:S}. Also, the followings hold:

\medskip \noindent \textnormal{(1)} The blow-up point $\xi_0$ is a critical point of the function $\tta$.

\medskip \noindent \textnormal{(2)} We have
\[\begin{cases}
\displaystyle \lim_{\e \to 0^+} \frac{\e \|u_{\e}\|_{L^{\infty}(\Omega)}^{\frac{N}{N-2}+1}} {\log \|u_{\e}\|_{L^{\infty}(\Omega)}}
= (p+1) \left|\S^{N-1}\right| b_{N,p}^{\frac{N}{N-2}} S^{\frac{p(q_0+1)}{1-pq_0}} \|U_{1,0}\|_{L^{q_0}(\R^N)}^{q_0} \tau(\xi_0)
&\text{if } p = \frac{N}{N-2},\\
\displaystyle \lim_{\e \to 0^+} \e \|u_{\e}\|_{L^{\infty}(\Omega)}^{\frac{N}{(N-2)p-2}+1}
= (N-2) S^{\frac{p(q_0+1)}{1-pq_0}} \|U_{1,0}\|_{L^{q_0}(\R^N)}^{q_0} \|V_{1,0}\|_{L^p(\R^N)}^p \tau(\xi_0)
&\text{if } p \in (\frac{N}{N-2}, \frac{N+2}{N-2})
\end{cases}\]
where $b_{N,p} := \lim_{|x| \to \infty} |y|^{N-2} V_{1,0}(y) \in (0, \infty)$ and $|\S^{N-1}|$ is the surface measure of the unit sphere $\S^{N-1}$ in
$\R^N$.\footnote{The well-definedness of the number $b_{N,p}$ was proved by Hulshof and Van der Vorst \cite{HV2}; see Lemma \ref{lemma:3} below.}

\medskip \noindent \textnormal{(3)} We have
\[\lim_{\e \to 0} \|u_{\e}\|_{L^{\infty}(\Omega)} v_{\e}(x) = \|U_{1,0}\|_{L^{q_0}(\R^N)}^{q_0} G(x,\xi_0)\]
and
\[\begin{cases}
\displaystyle \lim_{\e \to 0} \frac{\|u_{\e}\|_{L^{\infty}(\Omega)}^{N \over N-2}}
{\log \|u_{\e}\|_{L^{\infty}(\Omega)}} u_{\e}(x) = \frac{p+1}{N-2} \left|\S^{N-1}\right| b_{N,p}^{N \over N-2} G(x,\xi_0)
&\text{if } p = \frac{N}{N-2},\\
\displaystyle \lim_{\e \to 0} \|u_{\e}\|_{L^{\infty}(\Omega)}^{N \over (N-2)p-2} u_{\e}(x) = \|V_{1,0}\|_{L^p(\R^N)}^p G(x,\xi_0)
&\text{if } p \in (\frac{N}{N-2}, \frac{N+2}{N-2})
\end{cases}\]
in $C^1_{\textnormal{loc}}(\Omega \setminus \{\xi_0\})$-sense.
\end{thm}

Finally, we deduce the following existence theorem extending the result in \cite[Section 8]{KP}, which does not require the convexity assumption on $\Omega$.
Given any parameter $(\mu,\xi) \in (0,\infty) \times \R^N$, we set a bubble
\begin{equation}\label{eq:bubble}
(U_{\mu,\xi}(x),V_{\mu,\xi}(x)) = \(\mu^{-\frac{N}{q_0+1}} U_{1,0}(\mu^{-1}(x-\xi)), \mu^{-\frac{N}{p+1}} V_{1,0}(\mu^{-1}(x-\xi))\) \quad \text{for } x \in \R^N,
\end{equation}
which solves \eqref{eq:LERn}, and its projection $(PU_{\mu,\xi},PV_{\mu,\xi})$ on $W_0^{1,s}(\Omega) \times W_0^{1,\frac{s}{s-1}}(\Omega)$ for a suitable $s > 1$; refer to \eqref{eq:PUPV} for its precise definition.
\begin{thm}\label{thm:mainby}
Assume that $N \ge 4$, $\Omega$ is a smooth bounded domain in $\R^N$, $k \in \N$, $p \in (\max\{1,\frac{3}{N-2}\}, \frac{N+2}{N-2})$, and $q_{\e}$ is the number determined by \eqref{eq:cr-2}.
Given any $(\tdx) := (\td_1, \ldots, \td_k, x_1, \ldots, x_k) \in (0,\infty)^k \times \Omega^k$ and $i = 1, \ldots, k$,
let $\wth_{\tdx,i}: \Omega \to \R$ be the local $C^1$-regularization of $\wtg_{\tdx}$ around $x_i$ in \eqref{eq:wthdxi}. Then we set
\begin{equation}\label{eq:Upsilon}
\begin{aligned}
&\ \Upsilon_{p,k}(\tdx) \\
&= \begin{cases}
{\displaystyle \sum_{i=1}^k \td_i^{\frac{N}{q_0+1}} \wth_{\tdx,i}(x_i) - C_0 \log\(\td_1\cdots \td_k\)} \hspace{80pt} \text{if } p \in (\max\{1,\frac{3}{N-2}\}, \frac{N}{N-2}),\\
\displaystyle \sum_{i=1}^k \td_i^{N-2} \tau(x_i) - \frac{1}{2} \sum_{i,j=1, i \ne j}^k \(\td_i^{\frac{N}{p+1}} \td_j^{\frac{N}{q_0+1}} + \td_i^{\frac{N}{q_0+1}} \td_j^{\frac{N}{p+1}}\) G(x_i,x_j) - C_0 \log\(\td_1 \cdots \td_k\) \\
\hspace{305pt} \text{if } p \in [\frac{N}{N-2}, \frac{N+2}{N-2})
\end{cases}
\end{aligned}
\end{equation}
where the value of $C_0 > 0$ is given in Section \ref{sec:existence}.
Suppose that $(\widetilde{\bsde},\bsxi) := (\tilde{\delta}_1, \ldots, \tilde{\delta}_k, \xi_1,\ldots,\xi_k) \in (0,\infty)^k \times \Omega^k$
is an isolated critical point of $\Upsilon_{p,k}$ and $\Lambda$ is an open neighborhood of $(\widetilde{\bsde},\bsxi)$ such that
\begin{equation}\label{eq:deg}
\nabla \Upsilon_{p,k}(\tdx) \ne 0 \quad \text{for all } (\tdx) \in \pa \Lambda \quad \text{and} \quad \deg(\nabla \Upsilon_{p,k},\Lambda,0) \ne 0
\end{equation}
where $\deg$ designates the Brouwer degree. Then there exist a small number $\e_0 > 0$ and a family $\{(u_{\e},v_{\e})\}_{\e \in (0,\e_0)}$ of solutions of \eqref{eq:LE} with $q = q_{\e}$ such that
\[(u_{\e},v_{\e}) \simeq \sum_{i=1}^k (PU_{\mu_{i\e},x_{i\e}}, PV_{\mu_{i\e},x_{i\e}}) \quad \text{in } \Omega\]
at the main order. Here, $\td_{i\e} \to \tilde{\delta}_i$ and $x_{i\e} \to \xi_i$ as $\e \to 0$, and
\begin{equation}\label{eq:muie}
\mu_{i\e} := \begin{cases}
\e^{1 \over (N-2)p-2} \td_{i\e} &\text{if } p \in (\max\{1,\frac{3}{N-2}\}, \frac{N-1}{N-2}),\\
\displaystyle \left[\frac{-(N-2)\e}{W_{-1}(-(N-2)\e)}\right]^{1 \over N-2} \td_{i\e} &\text{if } p = \frac{N}{N-2},\\
\e^{1 \over N-2} \td_{i\e} &\text{if } p \in (\frac{N}{N-2}, \frac{N+2}{N-2})
\end{cases}
\end{equation}
for each $i = 1, \ldots, k$, where $W_{-1}$ stands for the Lambert $W$ function.\footnote{If $x \in [-\frac{1}{e},0)$, then $y = W_{-1}(x)$ is the unique solution of $ye^y = x$ such that $y \in (-\infty,-1]$.
It is known that $W_{-1}$ is a strictly deceasing function and $W_{-1}(x) \simeq \log(-x)$ 
for $x$ near $0$. In particular, $\left[\frac{-(N-2)\e}{W_{-1}(-(N-2)\e)}\right]^{1 \over N-2} \simeq \left[\frac{(N-2)\e}{|\log \e|}\right]^{1 \over N-2}$ for $\e > 0$ small.}
Particularly, the family $\{(u_{\e},v_{\e})\}_{\e \in (0,\e_0)}$ satisfies the uniform energy bound \eqref{eq:main0}.
\end{thm}

\begin{remark}
A couple of remarks regarding the above theorems are in order.

\medskip \noindent (1) In Theorems \ref{thm:main1} and \ref{thm:main2}, we use the convexity of the domain $\Omega$ only to show that the blow-up points are away from the boundary $\pa \Omega$.

Equation \eqref{eq:LEscalar} is well-behaved under the Kelvin transform, and in particular,
the moving plane method yields that the blow-up points for \eqref{eq:LEscalar} are away from $\pa \Omega$ for any smooth bounded domain.
Such an argument no longer works for system \eqref{eq:LE} in general, because it is rather ill-behaved under the Kelvin transform (see, e.g. \cite[Proposition 4.4]{CS}).

In \cite{CK}, Choi and the first author developed an alternative method to exclude the boundary blow-up behavior for the family of least energy solutions of \eqref{eq:LEs}, which mainly uses local Pohozaev identities.
Unfortunately, their idea is not directly applicable to our general situation.
Hence one has a challenging question: Is it possible to lift the convexity assumption on $\Omega$ in Theorems \ref{thm:main1} and \ref{thm:main2}?
If it is true, one will be also able to derive a general blow-up result for $p \in [\frac{N}{N-2}, \frac{N+2}{N-2})$ owing to our proof.

\medskip \noindent (2) If $p = q_0 = \frac{N+2}{N-2}$, then system \eqref{eq:LE} with $q = q_0$ and $u = v$ is reduced to the critical Lane-Emden equation.
By virtue of the invariance of \eqref{eq:LEscalar} under the Kelvin transform, an energy condition corresponding to \eqref{eq:main0} automatically holds; refer to Li \cite[Theorem 0.2]{Li}.
The analogous results to Theorems \ref{thm:main1} and \ref{thm:main2} were established by Barhi et al. \cite{BLR},
and one to Theorem \ref{thm:main3} was deduced by Grossi and Takahashi \cite[Theorem 2.5]{GT}.

If $p = 1$, then \eqref{eq:LE} is reduced to the biharmonic Lane-Emden equation.
For this case, Geng \cite{Ge} proved results corresponding to Theorems \ref{thm:main1} and \ref{thm:main2}, assuming that $\Omega$ is convex and $N \ge 5$.

\medskip \noindent (3) The condition $p \ge 1$ is often used in Section \ref{sec:dist}:
For instance, we need it when we apply the moving plane method to exclude the boundary blow-up.
We also require it in the proof of Lemmas \ref{lemma:Harnack} and \ref{lemma:7}, and Propositions \ref{prop:isos} and \ref{prop:dist}.

Moreover, the conditions $N \ge 4$ and $p > \frac{3}{N-2}$ are used in Section \ref{sec:loc}; see e.g. the proof of Lemmas \ref{lemma:loc1} and \ref{lemma:loc2}.
We wonder if these can be lifted in Theorems \ref{thm:main1}--\ref{thm:main3}.

\medskip \noindent (4) In Theorems \ref{thm:main2} and \ref{thm:main3}, we corrected some coefficients which were not properly given in \cite{G,CK}.

\medskip \noindent (5) For $p \in (1,\frac{N-1}{N-2})$, the function $\Upsilon_{p,k}$ in \eqref{eq:Upsilon} coincides with the function $F_{1\e}$ in \cite[(5.2)]{KP} up to a constant multiple.
\end{remark}

\subsection{Structure of the paper}
In Section \ref{sec:pre}, we provide several definitions and auxiliary lemmas needed throughout the paper.
In Section \ref{sec:dist}, we carry out a blow-up analysis of $(u_{\e},v_{\e})$ by modifying the argument of Schoen \cite{Sc}
(see also, e.g. \cite{Ge, Li, LZ}) to be suitable for examining the coupled system \eqref{eq:LE}. As a consequence, we show Theorem \ref{thm:main1}.
In Section \ref{sec:point}, we deduce various information on the shape of $(u_{\e},v_{\e})$.
From it, we find necessary conditions for the blow-up rates and locations in Section \ref{sec:loc}, establishing Theorems \ref{thm:main2} and \ref{thm:main3}.
In Section \ref{sec:existence}, we prove the existence result stated in Theorem \ref{thm:mainby}.

\subsection{Notations}
We list some notational conventions which will be used throughout the paper.

\medskip \noindent - Given $N \in \N$, we write $B(\xi,r) = \{x \in \R^N: |x-\xi| < r\}$ for each $\xi \in \R^N$ and $r > 0$.

\medskip \noindent - $\R^N_+ := \R^{N-1} \times (0,\infty)$ is the upper-half space in $\R^N$.

\medskip \noindent - The surface measure is denoted as $dS$.

\medskip \noindent - Let $D$ be a domain with boundary $\pa D$. The outward unit normal vector on $\pa D$ is written as $\nu$.

\medskip \noindent - Let (A) be a condition. We set $\mathbf{1}_{\text{(A)}} = 1$ if (A) holds and $0$ otherwise.

\medskip \noindent - For a function $f: \R^N \to \R$ and $l = 1, \ldots, N$, we write $\pa_l f(x) = \pa_{x_l}f(x)$.

\medskip \noindent - $C > 0$ is a generic constant independent of small numbers $\e > 0$ that may vary from line to line.

\section{Preliminary results}\label{sec:pre}
\subsection{Green's function $G$ and its related functions}\label{subsec:G}
Let $G$ be the Green's function of the Laplacian $-\Delta$ in $\Omega$ with the Dirichlet boundary condition, and $H$ the regular part of $G$, that is,
\begin{equation}\label{eq:H}
\begin{cases}
-\Delta_x H(x,\xi) = 0 &\text{for } x \in \Omega,\\
\displaystyle H(x,\xi) = \frac{\ga_N}{|x-\xi|^{N-2}}
&\displaystyle \text{for } x \in \pa\Omega \quad \text{where } \ga_N := \frac{1}{(N-2)\left|\S^{N-1}\right|}
\end{cases}
\end{equation}
for each $\xi \in \Omega$. It is well-known that
\begin{equation}\label{eq:Green}
0 < G(x,\xi) = G(\xi,x) = \frac{\ga_N}{|x-\xi|^{N-2}}-H(x,\xi) < \frac{\ga_N}{|x-\xi|^{N-2}} \quad \text{for } (x,\xi) \in \Omega \times \Omega, \ x \ne \xi.
\end{equation}
The Robin function $\tau$ is given by $\tau(x) = H(x,x)$ for $x \in \Omega$.

If $p \in (\frac{2}{N-2}, \frac{N}{N-2})$, then there is a number $\ep > 0$ depending only on $N$ and $p$ such that $G^p(\cdot,\xi) \in L^{1+\ep}(\Omega)$ for each $\xi \in \Omega$.
Accordingly, \eqref{eq:wtg} has the unique solution $\wtg(\cdot,\xi) \in W^{2,1+\ep}(\Omega)$ that can be represented by an integral.
Define the regular part $\wth: \Omega \times \Omega \to \R$ of $\wtg$ by
\[\wth(x,\xi) = \begin{cases}
\displaystyle \frac{\tga_{N,p,1}}{|x-\xi|^{(N-2)p-2}} - \wtg(x,\xi) &\text{if } p \in (\frac{2}{N-2}, \frac{N-1}{N-2}),\\
\displaystyle \frac{\tga_{N,p,1}}{|x-\xi|^{(N-2)p-2}} - \frac{\tga_{N,p,2} H(x,\xi)}{|x-\xi|^{(N-2)p-N}} - \wtg(x,\xi) &\text{if } p \in [\frac{N-1}{N-2}, \frac{N}{N-2})
\end{cases}\]
where
\begin{equation}\label{eq:tga}
\begin{cases}
\ga_N^p = \tga_{N,p,1}[(N-2)p-2][N-(N-2)p]&\text{for } p \in (\frac{2}{N-2}, \frac{N}{N-2}),\\
p\ga_N^{p-1} = \tga_{N,p,2}[(N-2)p-2(N-1)][N-(N-2)p] &\text{for } p \in [\frac{N-1}{N-2}, \frac{N}{N-2}).
\end{cases}
\end{equation}
Then, an application of elliptic regularity theory to the equation of $\wth$ yields that the map $x \in \Omega \mapsto \wth(x,\xi)$ belongs to $C^1_{\text{loc}}(\Omega)$ for each $\xi \in \Omega$; refer to \cite[Lemma 2.3]{CK} for the proof.
Set $\tta(x) = \wth(x,x)$ for $x \in \Omega$.

Given $\bsde = (\delta_1, \ldots, \delta_k) \in (0,\infty)^k$ and $\bsxi = (\xi_1, \ldots, \xi_k) \in \Omega^k$ such that $\xi_i \ne \xi_j$ for $1 \le i \ne j \le k$,
let $\wtg_{\dexi}:\Omega \to \R$ be the function defined by \eqref{eq:wtgdx} and $\wth_{\dexi}: \Omega \to \R$ the global $C^1$-regularization of $\wtg_{\dexi}$ given as
\begin{equation}\label{eq:wthdx}
\wth_{\dexi}(x) = \begin{cases}
\displaystyle \sum_{i=1}^k \frac{\tga_{N,p,1} \delta_i^{\frac{Np}{q_0+1}}}{|x-\xi_i|^{(N-2)p-2}} - \wtg_{\dexi}(x) &\text{if } p \in (\frac{2}{N-2}, \frac{N-1}{N-2}),\\
\displaystyle \sum_{i=1}^k \left[\frac{\tga_{N,p,1} \delta_i^{\frac{Np}{q_0+1}}}{|x-\xi_i|^{(N-2)p-2}}
- \frac{\tga_{N,p,2} A_{\dexi,i} \delta_i^{\frac{N(p-1)}{q_0+1}}}{|x-\xi_i|^{(N-2)p-N}}\right] - \wtg_{\dexi}(x) &\text{if } p \in [\frac{N-1}{N-2}, \frac{N}{N-2})
\end{cases}
\end{equation}
where
\begin{equation}\label{eq:Adxi}
A_{\dexi,i} := \delta_i^{\frac{N}{q_0+1}} \tau(\xi_i) - \sum_{j=1, j \ne i}^k \delta_j^{\frac{N}{q_0+1}} G(\xi_i,\xi_j).
\end{equation}
For each $i = 1, \ldots, k$, we also define the local $C^1$-regularization $\wth_{\dexi,i}: \Omega \to \R$ of $\wtg_{\dexi}$ around $\xi_i$ by
\begin{equation}\label{eq:wthdxi}
\wth_{\dexi,i}(x) = \begin{cases}
\displaystyle \frac{\tga_{N,p,1} \delta_i^{\frac{Np}{q_0+1}}}{|x-\xi_i|^{(N-2)p-2}} - \wtg_{\dexi}(x) &\text{if } p \in (\frac{2}{N-2}, \frac{N-1}{N-2}),\\
\displaystyle \frac{\tga_{N,p,1} \delta_i^{\frac{Np}{q_0+1}}}{|x-\xi_i|^{(N-2)p-2}}
- \frac{\tga_{N,p,2} A_{\dexi,i} \delta_i^{\frac{N(p-1)}{q_0+1}}}{|x-\xi_i|^{(N-2)p-N}} - \wtg_{\dexi}(x) &\text{if } p \in [\frac{N-1}{N-2}, \frac{N}{N-2})
\end{cases}
\end{equation}
Reasoning as in the proof of \cite[Lemma 2.11]{KP}, we observe that $\wth_{\dexi} \in C^1(\ovom)$ and $\wth_{\dexi,i}$ is continuously differentiable in a neighborhood of $\xi_i$.

\medskip
For future use, we further introduce the following functions: Let $\whh: \Omega \times \Omega \to \R$ be the function such that
\begin{equation}\label{eq:whh1}
-\Delta_x \whh(x,\xi) = 0 \quad \text{for } x \in \Omega
\end{equation}
and
\begin{equation}\label{eq:whh2}
\whh(x,\xi) = \begin{cases}
\displaystyle \frac{\ga_N \log|x-\xi|}{|x-\xi|^{N-2}} &\text{for } p = \frac{N}{N-2} \text{ and } x \in \pa\Omega,\\
\displaystyle \frac{\ga_N}{|x-\xi|^{(N-2)p-2}} &\text{for } p \in (\frac{2}{N-2}, \frac{N}{N-2}) \text{ and } x \in \pa\Omega
\end{cases}
\end{equation}
for each $\xi \in \Omega$. Furthermore, we set
\begin{equation}\label{eq:whg}
\whg(x,\xi) = \frac{\ga_N}{|x-\xi|^{(N-2)p-2}} - \whh(x,\xi) \quad \text{for } x, \xi \in \Omega,\, x \ne \xi.
\end{equation}
Given $p \in [\frac{N-1}{N-2}, \frac{N}{N-2})$, let $\ovh: \Omega \times \Omega \to \R$ be the function satisfying
\begin{equation}\label{eq:ovh}
\begin{cases}
-\Delta_x \ovh(x,\xi) = 0 &\text{for } x \in \Omega,\\
\ovh(x,\xi) = |x-\xi|^{N-(N-2)p} &\text{for } x \in \pa\Omega
\end{cases}
\end{equation}
for each $\xi \in \Omega$.

\subsection{Some results regarding the standard bubble $(U_{1,0},V_{1,0})$}
Recall that the standard bubble $(U_{1,0},V_{1,0})$ is the unique least energy solution of \eqref{eq:LERn} satisfying \eqref{eq:LERn2}.
The first result on this pair is due to Hulshof and Van der Vorst \cite{HV2}.
\begin{lemma}\label{lemma:3}
The pair $(U_{1,0},V_{1,0})$ is unique, radially symmetric, and decreasing in the radial variable.
Also, there exist numbers $a_{N,p},\, b_{N,p} > 0$ depending only on $N$ and $p$ such that
\begin{equation}\label{eq:HV}
\begin{cases}
\lim\limits_{|y| \to \infty} |y|^{N-2} U_{1,0}(y) = a_{N,p} &\text{if } p \in (\frac{N}{N-2}, \frac{N+2}{N-2}],\\
\lim\limits_{|y| \to \infty} \dfrac{|y|^{N-2}}{\log |y|} U_{1,0}(y) = a_{N,p} &\text{if } p = \frac{N}{N-2},\\
\lim\limits_{|y| \to \infty} |y|^{(N-2)p-2} U_{1,0}(y) = a_{N,p} &\text{if } p \in (\frac{2}{N-2}, \frac{N}{N-2})
\end{cases}
\quad \text{and} \quad
\lim_{|y| \to \infty} |y|^{N-2} V_{1,0}(y) = b_{N,p}.
\end{equation}
If $p \in (\frac{2}{N-2}, \frac{N}{N-2})$, we also have that
\begin{equation}\label{eq:HVab}
b_{N,p}^p = a_{N,p}[(N-2)p-2][N-(N-2)p].
\end{equation}
\end{lemma}

In the next lemma, we present a refinement of \eqref{eq:HV}, which generalizes \cite[Corollaries 2.6 and 2.7]{KP}.
Our proof is different from and simpler than one in \cite{KP}.
For the sake of convenience, we write $U_{1,0}(y) = U_{1,0}(|y|)$ and $V_{1,0}(y) = V_{1,0}(|y|)$.
\begin{lemma}\label{lemma:UVdecay}
Assume that $r = |y| \ge 1$ and $(U_{1,0}',V_{1,0}')$ denotes the derivative of $(U_{1,0},V_{1,0})$ with respect to $r$. There exists a constant $C > 0$ depending only on $N$ and $p$ such that
\begin{equation}\label{eq:V10est}
\left|V_{1,0}(r) - \frac{b_{N,p}}{r^{N-2}}\right| \le \frac{C}{r^N} \quad \text{and} \quad \left|V_{1,0}'(r) + \frac{(N-2)b_{N,p}}{r^{N-1}} \right| \le \frac{C}{r^{N+1}}.
\end{equation}
Besides, if $p \in (\frac{N}{N-2}, \frac{N+2}{N-2}]$, then
\begin{equation}\label{eq:U10est1}
\left|U_{1,0}(r) - \frac{a_{N,p}}{r^{N-2}}\right| \le \frac{C}{r^{N-2+\kappa_0}} \quad \text{and} \quad \left|U_{1,0}'(r) + \frac{(N-2)a_{N,p}}{r^{N-1}} \right| \le \frac{C}{r^{N-1+\kappa_0}}
\end{equation}
where $\kappa_0 := (N-2)p-N > 0$. If $p = \frac{N}{N-2}$, then
\begin{equation}\label{eq:U10est2}
\left|U_{1,0}(r) - \frac{a_{N,p} \log r}{r^{N-2}}\right| \le \frac{C}{r^{N-2}} \quad \text{and} \quad \left|U_{1,0}'(r) + \frac{(N-2)a_{N,p}\log r}{r^{N-1}}\right| \le \frac{C}{r^{N-1}}.
\end{equation}
If $p \in (\frac{2}{N-2}, \frac{N}{N-2})$, then
\begin{equation}\label{eq:U10est3}
\left|U_{1,0}(r) - \frac{a_{N,p}}{r^{(N-2)p-2}}\right| \le \frac{C}{r^{(N-2)p-2+\kappa_1}} \quad \text{and} \quad \left|U_{1,0}'(r) + \frac{((N-2)p-2)a_{N,p}}{r^{(N-2)p-1}}\right| \le \frac{C}{r^{(N-2)p-1+\kappa_1}}
\end{equation}
where $\kappa_1 \in (0, \min\{N-(N-2)p,((N-2)p-2)q_0-N)\})$.\footnote{We have that $((N-2)p-2)q_0 > N+2$ as shown in \cite[Lemma A.1]{KP}.}
\end{lemma}
\begin{proof}
In (3.22)--(3.24) of \cite{HV2}, it was found that
\begin{equation}\label{eq:HV2}
\lim_{r \to \infty} \frac{rU_{1,0}'(r)}{U_{1,0}(r)} = \begin{cases}
2-N &\text{if } p \in [\frac{N}{N-2}, \frac{N+2}{N-2}],\\
2-(N-2)p &\text{if } p \in (\frac{2}{N-2}, \frac{N}{N-2}),
\end{cases}
\quad \text{and} \quad
\lim_{r \to \infty} \frac{rV_{1,0}'(r)}{V_{1,0}(r)} = 2-N.
\end{equation}
Applying \eqref{eq:HV}, \eqref{eq:HV2}, \eqref{eq:LERn}, and the radial symmetry of $(U_{1,0},V_{1,0})$, we see
\begin{equation}\label{eq:V10esta}
\begin{aligned}
(2-N)b_{N,p} - r^{N-1}V_{1,0}'(r) &= \int_r^\infty (t^{N-1}V_{1,0}'(t))' dt = -\int_r^{\infty} U_{1,0}^{q_0}(t) t^{N-1}dt \\
&= \begin{cases}
\displaystyle O(r^{N-(N-2)q_0}) &\text{if } p \in (\frac{N}{N-2}, \frac{N+2}{N-2}],\\
\displaystyle O(r^{N-(N-2)q_0} \log r) &\text{if } p = \frac{N}{N-2},\\
\displaystyle O(r^{N-((N-2)p-2)q_0}) &\text{if } p \in (\frac{2}{N-2}, \frac{N}{N-2}),
\end{cases} \\
&= O(r^{-2})
\end{aligned}
\end{equation}
as $r \to \infty$, so the second inequality in \eqref{eq:V10est} holds.
By integrating it over $[r,\infty)$, we obtain the first inequality in \eqref{eq:V10est}.

To study the decay of $U_{1,0}(r)$ and $U_{1,0}'(r)$ as $r \to \infty$, we consider three different cases.

\medskip \noindent \textbf{Case 1: $p = (\frac{N}{N-2}, \frac{N+2}{N-2}]$.}
In this case, arguing as in \eqref{eq:V10esta}, we discover \eqref{eq:U10est1}.

\medskip \noindent \textbf{Case 2: $p = \frac{N}{N-2}$.} By \eqref{eq:V10est},
\begin{align*}
r^{N-1}U_{1,0}'(r) - U_{1,0}'(1) &= - \int_1^r V_{1,0}^p(t) t^{N-1}dt = -\int_1^r t^{N-1} \(b_{N,p}\, t^{2-N}+O(t^{-N})\)^p dt\\
&=-\int_1^r t^{N-1} \(b_{N,p}^p t^{-N}+O(t^{-N-2})\) dt = -b_{N,p}^p \log r + O(1),
\end{align*}
which implies that
\[U_{1,0}'(r) = -b_{N,p}^p r^{1-N} \log r + O(r^{1-N}).\]
Integrating it over $[r,\infty)$, we obtain
\begin{equation}\label{eq:U10est21}
U_{1,0}(r) = \frac{b_{N,p}^p}{N-2} r^{2-N} \log r + O(r^{2-N}).
\end{equation}
Comparing \eqref{eq:U10est21} with \eqref{eq:HV}, we get that $b_{N,p}^p = (N-2)a_{N,p}$. From this, we deduce \eqref{eq:U10est2}.

\medskip \noindent \textbf{Case 3: $p \in (\frac{2}{N-2}, \frac{N}{N-2})$.} From \eqref{eq:V10esta}, we observe
\begin{align*}
(r^{N-1}U_{1,0}'(r))' &= -V_{1,0}^p(r) r^{N-1} = -r^{N-1}(b_{N,p}r^{2-N}+O(r^{2-((N-2)p-2)q_0}))^p \\
&= -b_{N,p}^p r^{N-1-(N-2)p} + O(r^{N-1-(N-2)p+N-((N-2)p-2)q_0}).
\end{align*}
Let
\[U^*(r) = \frac{b_{N,p}^p}{[(N-2)p-2][N-(N-2)p]} r^{2-(N-2)p} = a_{N,p}\, r^{2-(N-2)p}\]
where the second equality comes from \eqref{eq:HVab}. Then
\[\(r^{N-1}\(U_{1,0}'-(U^*)'\)(r)\)' = O(r^{N-1-(N-2)p+N-((N-2)p-2)q_0}),\]
so integrating it over $[1,r]$ gives
\[U_{1,0}'(r) + a_{N,p}((N-2)p-2) r^{1-(N-2)p} = O(r^{1-N}) + O(r^{1-(N-2)p} \cdot r^{N-((N-2)p-2)q_0} \cdot \log r).\]
This equality leads to \eqref{eq:U10est3}.
\end{proof}

Using the previous lemma, one can also relate the numbers $a_{N,p}$ and $b_{N,p}$ when $p \in [\frac{N}{N-2}, \frac{N+2}{N-2}]$; cf. \eqref{eq:HVab}.
\begin{lemma}
Let
\begin{equation}\label{eq:A12}
A_1 = \int_{\R^N} U_{1,0}^{q_0} \quad \text{and} \quad A_2 = \int_{\R^N} V_{1,0}^p.
\end{equation}
If $p \in (\frac{N}{N-2}, \frac{N+2}{N-2}]$, then
\begin{equation}\label{eq:ab}
a_{N,p} A_1 = b_{N,p} A_2.
\end{equation}
If $p = \frac{N}{N-2}$, then
\begin{equation}\label{eq:ab1}
b_{N,p}^p = (N-2)a_{N,p} \quad \text{and} \quad a_{N,p}A_1 = a_{N,p}b_{N,p} \ga_N^{-1} = b_{N,p}^{p+1}\left|\S^{N-1}\right|.
\end{equation}
\end{lemma}
\begin{proof}
Suppose that $p \in (\frac{N}{N-2}, \frac{N+2}{N-2}]$. If $U_{1,0}^*$ denotes the Kelvin transform of $U_{1,0}$, then the representation formula gives
\[U_{1,0}^*(y) = \int_{\R^N} \frac{\ga_N}{|y-x|x|^{-2}|^{N-2}} V_{1,0}^p(x) \frac{dx}{|x|^{N-2}}, \quad \text{and so} \quad U_{1,0}^*(0) = \ga_N A_2\]
where $\ga_N > 0$ is the value in \eqref{eq:H}. We discover from \eqref{eq:HV} that
\begin{equation}\label{eq:agaA2}
a_{N,p} = \lim_{|y| \to \infty} |y|^{N-2} U_{1,0}(y) = \lim_{|y| \to 0} U_{1,0}^*(y) = U_{1,0}^*(0) = \ga_N A_2.
\end{equation}
Similarly, we have $b_{N,p} = \ga_N A_1$, from which we obtain \eqref{eq:ab}.

\medskip
Assume that $p = \frac{N}{N-2}$. The first equality in \eqref{eq:ab1} was shown in the proof of Lemma \ref{lemma:UVdecay}.
Also, since the relation $b_{N,p} = \ga_N A_1$ still holds, the second equality is true. The third equality is a consequence of $\ga_N^{-1} = (N-2)|\S^{N-1}|$ and the first one.
\end{proof}

Let
\[\(\Psi_{1,0}^0(x), \Phi_{1,0}^0(x)\) = \(x \cdot \nabla U_{1,0}(x) + \frac{NU_{1,0}(x)}{q_0+1}, x \cdot \nabla V_{1,0}(x) + \frac{NV_{1,0}(x)}{p+1}\)\]
and
\[\(\Psi_{1,0}^l(x), \Phi_{1,0}^l(x)\) = (\pa_l U_{1,0}(x), \pa_l V_{1,0}(x))\]
for $x \in \R^N$ and $l = 1, \ldots, N$. Also, given any parameter $(\mu,\xi) \in (0,\infty) \times \Omega$, we set
\begin{equation}\label{eq:bubble21}
\(\Psi_{\mu,\xi}^0(x), \Phi_{\mu,\xi}^0(x)\) = \(\mu^{-\frac{N}{q_0+1}} \Psi_{1,0}^l (\mu^{-1}(x-\xi)), \mu^{-\frac{N}{p+1}} \Phi_{1,0}^l (\mu^{-1}(x-\xi))\)
\end{equation}
and
\begin{equation}\label{eq:bubble22}
\(\Psi_{\mu,\xi}^l(x), \Phi_{\mu,\xi}^l(x)\) = \(\mu^{-\frac{N}{q_0+1}-1} \Psi_{1,0}^l (\mu^{-1}(x-\xi)), \mu^{-\frac{N}{p+1}-1} \Phi_{1,0}^l (\mu^{-1}(x-\xi))\)
\end{equation}
for $x \in \R^N$ and $l = 1, \ldots, N$. In \cite{FKP}, Frank et al. proved a non-degeneracy result for the bubbles.
\begin{lemma}\label{lemma:FKP}
Recall the pair $(U_{\mu,\xi}, V_{\mu,\xi})$ defined in \eqref{eq:bubble}. The space of solutions of the linearized system
\[\begin{cases}
-\Delta \Psi = pV_{\mu,\xi}^{p-1} \Phi \quad \ \ \text{in } \R^N,\\
-\Delta \Phi = q_0U_{\mu,\xi}^{q_0-1} \Psi \quad \text{in } \R^N,\\
(\Psi,\Phi) \in \dot{W}^{2,\frac{p+1}{p}}(\R^N) \times \dot{W}^{2,\frac{q_0+1}{q_0}}(\R^N) \text{ or } \lim\limits_{|x| \to \infty} (\Psi(x),\Phi(x)) = (0,0)
\end{cases}\]
is spanned by
\[\left\{\(\Psi_{\mu,\xi}^0, \Phi_{\mu,\xi}^0\), \(\Psi_{\mu,\xi}^1,\Phi_{\mu,\xi}^1\), \ldots, \(\Psi_{\mu,\xi}^N,\Phi_{\mu,\xi}^N\)\right\}.\]
\end{lemma}

\subsection{Algebraic property of $(p,q_{\e})$}
We shall need the following elementary lemma on $(p,q_{\e})$.
\begin{lemma}
Let $(\alpha_{\e}, \beta_{\e})$ be the pair in \eqref{eq:albe}. If $\e > 0$ is small enough, then
\begin{equation}\label{eq:pqrel}
\max\{\alpha_{\e},\beta_{\e}\} < \min\{N-2, (N-2)p-2\}
\end{equation}
and
\begin{equation}\label{eq:pqrel2}
[(N-2)p-2]\beta_{\e} > (N-2)\alpha_{\e}.
\end{equation}
\end{lemma}
\begin{proof}
By \eqref{eq:cr-hy}, we have that $(N-2)(pq_0-1) = 2(p+q_0+2)$. Thus, if we let $\alpha_0$ and $\beta_0$ be the numbers given by \eqref{eq:albe} with $\e = 0$, then
\[\alpha_0 = \frac{2(p+1)}{\frac{2}{N-2}(p+1)+\frac{2}{N-2}(q_0+1)} < N-2 \quad \text{and} \quad \beta_0 < N-2.\]
Besides, because $\beta_0 p -\alpha_0 = \alpha_0 q_0 - \beta_0 = 2$, we see
\[\alpha_0 q_0 - pq_0(N-2) + 2q_0 = \beta_0 - pq_0(N-2) + 2(q_0+1) = \beta_0-(N-2)-2(p+1) < 0,\]
which implies that $\alpha_0 < (N-2)p-2$. Similarly, $\beta_0 < (N-2)p-2$.
Therefore \eqref{eq:pqrel} is true for every small $\e > 0$.

On the other hand, we have
\[[(N-2)p-2](q_0+1) = N(p+1) > (N-2)(p+1),\]
so \eqref{eq:pqrel2} is valid for $\e > 0$ small.
\end{proof}

\subsection{Local Pohozaev-type identity}
The following Pohozaev-type identity will be useful in the sequel.
\begin{lemma}
Let $\{(u_{\e},v_{\e})\}_{\e \in (0,\e_0)}$ be a family of solutions of \eqref{eq:LE} with $q = q_{\e}$.
If $B(x_0,\rho) \subset \Omega$ for some $x_0 \in \Omega$ and $\rho > 0$, then
\begin{multline}\label{eq:poho}
\e \int_{B(x_0,\rho)} \nabla u_{\e} \cdot \nabla v_{\e} = \rho \int_{\pa B(x_0,\rho)} \(2 \frac{\pa u_{\e}}{\pa \nu} \frac{\pa v_{\e}}{\pa \nu}-\nabla u_{\e} \cdot \nabla v_{\e}\) dS \\
+ \rho \int_{\pa B(x_0,\rho)} \(\frac{v_{\e}^{p+1}}{p+1} + \frac{u_{\e}^{q_{\e}+1}}{q_{\e}+1}\) dS
+ N\int_{\pa B(x_0,\rho)} \(\frac{v_{\e}}{p+1}\frac{\pa u_{\e}}{\pa \nu} + \frac{u_{\e}}{q_{\e}+1}\frac{\pa v_{\e}}{\pa \nu}\) dS.
\end{multline}
\end{lemma}
\begin{proof}
By multiplying the first equation of \eqref{eq:LE} by $(x-x_0) \cdot \nabla v_{\e}$ or $v_{\e}$,
and integrating the both sides over $B(x_0,\rho)$, we obtain two identities.
Also, by multiplying the second equation of \eqref{eq:LE} by $(x-x_0) \cdot \nabla u_{\e}$ or $u_{\e}$,
and integrating the both sides over $B(x_0,\rho)$, we obtain two more identities.
Combining them and applying \eqref{eq:cr-2}, we deduce \eqref{eq:poho}.
\end{proof}

\section{Estimates on distances between blow-up points}\label{sec:dist}
Throughout this section, we assume that $N \ge 3$, $p \in (\frac{2}{N-2}, \frac{N+2}{N-2})$, and $p \ge 1$.

Let $\{(u_{\e},v_{\e})\}_{\e \in (0,\e_0)}$ be a family of solutions of \eqref{eq:LE} with $q = q_{\e}$.
The convexity of the domain $\Omega$, the moving plane method, and elliptic estimates
guarantee the existence of a small number $r_0 > 0$ and a large number $M_0 > 0$ such that
\[0 < (u_{\e} + v_{\e})(x) \le M_0 \quad \text{for all } x \in \Omega \text{ such that } \dist(x,\pa\Omega) \le 2r_0;\]
refer to Page 188 in Guerra \cite{G} for the proof. In particular, all blow-up points of $\{(u_{\e},v_{\e})\}_{\e \in (0,\e_0)}$ are away from the boundary $\pa\Omega$.

\subsection{Various types of blow-up points}
We present the notion of a blow-up point, an isolated blow-up point, and an isolated simple blow-up point of solutions to system \eqref{eq:LE} by altering the original definition in \cite{Sc}.
\begin{definition}\label{def:blowup}
Let $\{\e_n\}_{n \in \N}$ be a sequence of positive small numbers such that $\e_n \to 0$ as $n \to \infty$.
Let also $\{(u_n,v_n) = (u_{\e_n},v_{\e_n})\}_{n \in \N}$ be a sequence of solutions of \eqref{eq:LE} with $q = q_n := q_{\e_n}$.

\medskip \noindent \textnormal{(1)} A point $\xi \in \Omega$ is called a blow-up point of $\{(u_n,v_n)\}_{n \in \N}$
if there exists a sequence $\{x_n\}_{n \in \N} \subset \Omega$ such that $x_n$ is a local maximum point of $u_n$,
\[u_n(x_n) \to \infty \quad \text{and} \quad x_n \to \xi \in \Omega \quad \text{as } n \to \infty.\]
For the sake of convenience, we will often say that $x_n \to \xi$ is a blow-up point of $\{(u_n,v_n)\}_{n \in \N}$.

\medskip \noindent \textnormal{(2)} A point $\xi \in \Omega$ is called an isolated blow-up point of $\{(u_n,v_n)\}_{n \in \N}$
if $\xi$ is a blow-up point and there exist numbers $\rho_0 \in (0,r_0]$ small and $C > 0$ (independent of $n \in \N$) such that
\begin{equation}\label{eq:iso}
u_n(x) \le C |x-x_n|^{-\alpha_n} \quad \text{and} \quad v_n(x) \le C |x-x_n|^{-\beta_n} \quad \text{for all } x \in B(x_n,\rho_0) \setminus \{x_n\}.
\end{equation}
Here, $\alpha_n := \alpha_{\e_n}$ and $\beta_n := \beta_{\e_n}$ are the numbers in \eqref{eq:albe}.

\medskip \noindent
\textnormal{(3)} Define the spherical average of $u_n$ and $v_n$ by
\[\bar{u}_n(r) = \frac{1}{|\pa B(x_n,r)|} \int_{\pa B(x_n,r)} u_n\, dS \quad \text{and} \quad \bar{v}_n(r) = \frac{1}{|\pa B(x_n,r)|} \int_{\pa B(x_n,r)} v_n\, dS\]
for $r \in (0,\rho_0)$.
We say that an isolated blow-up point $\xi \in \Omega$ of $\{(u_n,v_n)\}_{n \in \N}$ is simple
if $\xi$ is an isolated blow-up point and each map $r \mapsto r^{\alpha_n} \bar{u}_n$ and $r \mapsto r^{\beta_n} \bar{v}_n$
has only one critical point in the interval $(0,\rho_0)$ after the value of $\rho_0 > 0$ is suitably reduced.
\end{definition}

\subsection{Isolated blow-up points}
As a preparation, we deduce an annular Harnack inequality.
\begin{lemma}\label{lemma:Harnack}
Let $x_n \to \xi$ be an isolated blow-up point of $\{(u_n,v_n)\}_{n \in \N}$. For any $r \in (0,\frac{\rho_0}{3})$, we have
\begin{equation}\label{eq:Harnack1}
\max_{x \in B(x_n,2r) \setminus B(x_n,{r \over 2})} u_n(x) \le C \min_{x \in B(x_n,2r) \setminus B(x_n,{r \over 2})} u_n(x)
\end{equation}
and
\begin{equation}\label{eq:Harnack2}
\max_{x \in B(x_n,2r) \setminus B(x_n,{r \over 2})} v_n(x) \le C \min_{x \in B(x_n,2r) \setminus B(x_n,{r \over 2})} v_n(x)
\end{equation}
for $C > 0$ independent of $n \in \N$ and $r \in (0,\frac{\rho_0}{3})$.
\end{lemma}
\begin{proof}
Given $r \in (0,\frac{\rho_0}{3})$, we set
\[(U_{n;r}(y), V_{n;r}(y)) = \(r^{\alpha_n} u_n(ry+x_n), r^{\beta_n} v_n(ry+x_n)\) \quad \text{for } y \in B(0,3).\]
Since $x_n \to \xi$ is an isolated blow-up point, we see from \eqref{eq:iso} that
\[U_{n;r}(y) \le C|y|^{-\alpha_n} \le C \quad \text{and} \quad V_{n;r}(y) \le C |y|^{-\beta_n} \le C \quad \text{for all } y \in B\(0,\frac{5}{2}\) \setminus \overline{B\(0,\frac{1}{4}\)}.\]
Also, $(U_{n;r}, V_{n;r})$ is a solution of a weakly coupled cooperative elliptic system
\begin{equation}\label{eq:UVer}
\begin{pmatrix}
\Delta U_{n;r} \\ \Delta V_{n;r}
\end{pmatrix} +
\begin{pmatrix}
0 & V_{n;r}^{p-1} \\
U_{n;r}^{q_n-1} & 0
\end{pmatrix}
\begin{pmatrix}
U_{n;r} \\ V_{n;r}
\end{pmatrix} = \begin{pmatrix}
0 \\ 0
\end{pmatrix}
\quad \text{in } B\(0,\frac{5}{2}\) \setminus \overline{B\(0,\frac{1}{4}\)}.
\end{equation}
Under the assumption that $p \ge 1$, system \eqref{eq:UVer} satisfies all necessary conditions to apply Harnack's inequality in \cite[Theorem 1.1]{CZ}.
Consequently, \eqref{eq:Harnack1}--\eqref{eq:Harnack2} is valid.
\end{proof}

We next show that in the case of an isolated blow-up point, a renormalization of $(u_n,v_n)$ tends to the standard bubble $(U_{1,0},V_{1,0})$ as $n \to \infty$.
\begin{lemma}\label{lemma:7}
Let $x_n \to \xi$ be an isolated blow-up point of $\{(u_n,v_n)\}_{n \in \N}$ and $\lambda_n = u_n^{1/\alpha_n}({x_n})$.
Assume that $\{R_n\}_{n \in \N}$ and $\{\vsi_n\}_{n \in \N}$ are families of numbers such that $R_n \to \infty$ and $\vsi_n \to 0$ as $n \to \infty$.
After possibly passing to a subsequence, we have
\begin{equation}\label{eq:72}
\left\|\lambda_n^{-\alpha_n} u_n(\lambda_n^{-1}\cdot + x_n) - U_{1,0}\right\|_{C^2(\overline{B(0,2R_n)})}
+ \left\|\lambda_n^{-\beta_n} v_n(\lambda_n^{-1}\cdot + x_n) - V_{1,0}\right\|_{C^2(\overline{B(0,2R_n)})} \le \vsi_n
\end{equation}
and $R_n \lambda_n^{-1} \to 0$ as $n \to \infty$.
\end{lemma}
\begin{proof}
The proof is similar to that of \cite[Proposition 2.1]{Li} or \cite[Lemma 3.2]{LZ}. We sketch it to point out the necessary modifications.

\medskip
Let
\[(U_n(y),V_n(y)) = \(\lambda_n^{-\alpha_n} u_n(\lambda_n^{-1}y + x_n), \lambda_n^{-\beta_n} v_n(\lambda_n^{-1}y + x_n)\) \quad \text{for } y \in \lambda_n(\Omega-x_n).\]
By virtue of \eqref{eq:LE}, Definition \ref{def:blowup} (1), and \eqref{eq:iso}, we have
\begin{equation}\label{eq:71}
\begin{cases}
-\Delta U_n = V_n^p,\ -\Delta V_n = U_n^{q_n} \quad \text{in } B(0,\lambda_n\rho_0),\\
U_n(0) = 1, \ \nabla U_n(0) = 0,\\
U_n(y) \le C|y|^{-\alpha_n} \quad \text{and} \quad V_n(y) \le C |y|^{-\beta_n} \quad \text{for } y \in B(0,\lambda_n\rho_0).
\end{cases}
\end{equation}
On the other hand, we infer from the maximum principle and \eqref{eq:Harnack1} that $U_n(y) \le C$ for $|y| \le 1$.
Combining this with the last inequality of \eqref{eq:71}, we obtain that $U_n(y) \le C$ for $|y| \le \lambda_n\rho_0$.
Moreover, the standard elliptic regularity theory and \eqref{eq:pqrel} yield
\[\|V_n\|_{L^{\infty}(B(y,1))} \le C\(\|V_n\|_{L^1(B(y,2))} + \|U_n^{q_n}\|_{L^{\infty}(B(y,2))}\) \le C \quad \text{for any } |y| \le 1,\]
which implies that $V_n(y) \le C$ whenever $|y| \le \lambda_n\rho_0$.
A further application of elliptic regularity shows that there exists a pair $(\mcu_0,\mcv_0)$ of smooth positive functions in $\R^N$ such that
\[(U_n,V_n) \to (\mcu_0,\mcv_0) \quad \text{in } C^2_{\text{loc}}(\R^N) \times C^2_{\text{loc}}(\R^N) \quad \text{as } n \to \infty,\]
along a subsequence. It holds that
\[\begin{cases}
-\Delta \mcu_0 = \mcv_0^p,\ -\Delta \mcv_0 = \mcu_0^{\,q_0} \quad \text{in } \R^N,\\
\mcu_0(0) = 1, \ \nabla \mcu_0(0) = 0,\\
\displaystyle \mcu_0(y) \le \frac{C}{1+|y|^{\alpha_0}} \quad \text{and} \quad \mcv_0(y) \le \frac{C}{1+|y|^{\beta_0}} \quad \text{for } y \in \R^N.
\end{cases}\]
Using the fact that $\alpha_0q_0 > 2$ and $\beta_0p > 2$, we conclude
\[\mcu_0(x) = \int_{\R^N} \frac{\ga_N}{|x-y|^{N-2}} \, \mcv_0^p(y) dy \quad \text{and} \quad \mcv_0(x) = \int_{\R^N} \frac{\ga_N}{|x-y|^{N-2}} \, \mcu_0^{q_0}(y) dy \quad \text{for } x \in \R^N\]
where $\ga_N > 0$ is the constant in \eqref{eq:H}. Thanks to the classification theorem of Chen et al. \cite{CLO},
it follows that $(\mcu_0,\mcv_0) = (U_{1,0},V_{1,0})$ where $(U_{1,0},V_{1,0})$ is the standard bubble provided $p \ge 1$.
\end{proof}
\begin{remark}\label{rmk:7}
The above lemma and \eqref{eq:pqrel} imply that each map $r \mapsto r^{\alpha_n} \bar{u}_n$ and $r \mapsto r^{\beta_n} \bar{v}_n$ has a critical point in $(0,R_n\lambda_n^{-1})$.
In particular, if $\xi$ is isolated simple blow-up point of $\{(u_n,v_n)\}_{n \in \N}$, then we may assume that
\begin{equation}\label{eq:uvmono}
(r^{\alpha_n} \bar{u}_n(r))' < 0 \quad \text{and} \quad (r^{\beta_n} \bar{v}_n(r))' < 0 \quad \text{for any } R_n \lambda_n^{-1} \le r \le \rho_0 \text{ and } n \in \N
\end{equation}
by discarding a first few $(u_n,v_n)$'s.
\end{remark}

\subsection{Isolated simple blow-up points}
We concern a decay estimate of $(u_n,v_n)$ in a small ball (of a fixed radius) centered at an isolated simple blow-up point.
\begin{lemma}\label{lemma:8}
Assume that Lemma \ref{lemma:7} holds with $R_n \to \infty$ and $0 < \vsi_n < e^{-R_n}$.
If $x_n \to \xi$ is an isolated simple blow-up point of $\{(u_n,v_n)\}_{n \in \N}$, then for a sufficiently small $\eta > 0$,
there exist $C > 0$ and $\rho_1 \in (0,\frac{\rho_0}{3})$ independent of $n \in \N$ (but dependent on $\eta$) such that
\begin{equation}\label{eq:801}
v_n(x) \le C \lambda_n^{\beta_n-(N-2)+\eta} |x-x_n|^{-(N-2)+\eta}
\end{equation}
and
\begin{equation}\label{eq:802}
u_n(x) \le \begin{cases}
C\lambda_n^{\alpha_n-(N-2)+\eta} |x-x_n|^{-(N-2)+\eta} &\text{if } p \in (\frac{N}{N-2}, \frac{N+2}{N-2}),\\
C \lambda_n^{\alpha_n-(N-2)+\eta p} |x-x_n|^{-(N-2)+\eta p} \log (\lambda_n |x-x_n|) &\text{if } p = \frac{N}{N-2},\\
C\lambda_n^{\alpha_n+2-(N-2)p+\eta p} |x-x_n|^{2-(N-2)p+\eta p} &\text{if } p \in (\frac{2}{N-2}, \frac{N}{N-2})
\end{cases}
\end{equation}
for $r_n := R_n \lambda_n^{-1} \le |x-x_n| \le \rho_1$.
\end{lemma}
\begin{proof}
Throughout the proof, we write $r = |x-x_n|$.

\medskip
By \eqref{eq:Harnack1}, Definition \ref{def:blowup} (3), \eqref{eq:72}, and \eqref{eq:HV},
\begin{equation}\label{eq:uedecay}
\begin{aligned}
r^{\alpha_n} u_n(x) &\le C r^{\alpha_n} \bar{u}_n(r) \le C r_n^{\alpha_n} \bar{u}_n(r_n) \\
&\le \begin{cases}
C r_n^{\alpha_n} \lambda_n^{\alpha_n} R_n^{2-N} = C R_n^{\alpha_n-(N-2)} &\text{if } p \in (\frac{N}{N-2}, \frac{N+2}{N-2}),\\
C r_n^{\alpha_n} \lambda_n^{\alpha_n} R_n^{2-N} \log R_n = C R_n^{\alpha_n-(N-2)} \log R_n &\text{if } p = \frac{N}{N-2},\\
C r_n^{\alpha_n} \lambda_n^{\alpha_n} R_n^{2-(N-2)p} = C R_n^{\alpha_n+2-(N-2)p} &\text{if } p \in (\frac{2}{N-2}, \frac{N}{N-2})
\end{cases}
\end{aligned}
\end{equation}
for $r_n \le r \le \rho_1$. We have similar inequalities for $v_n$.

Define
\[L_n(\phi,\psi) = \(\Delta \phi + v_n^{p-1}\psi, \Delta \psi + u_n^{q_n-1}\phi\)\]
for a pair $(\phi,\psi)$ of functions on $\overline{B(x_n,\rho_1)} \setminus B(x_n,r_n)$.
We know that $u_n$ and $v_n$ are positive on their domains, and $L_n(u_n,v_n) = 0$.
Therefore a minor modification of the proof of \cite[Theorem 1.2]{dFM} shows that the maximum principle holds for $L_n$.
Let $b_n = \beta_n(p-1)$. Because of the condition $1 \le p < \frac{N+2}{N-2}$, we may assume that
\begin{equation}\label{eq:be}
0 \le b_n = \frac{2(p-1)(q_n+1)}{pq_n-1} < 2 \quad \text{and} \quad \alpha_n (q_n-1) = 4-b_n > 2 \quad \text{for all } n \in \N.
\end{equation}
Also, if $N = 4$, then $p > \frac{2}{N-2} = 1$ so $b_n$ is away from $0$. If $N = 3$, then $p > 2$ and so $N-4+b_n = b_n-1$ is positive and away from $0$.

Given a sufficiently small number $\eta \in (0,1)$, a direct calculation with \eqref{eq:pqrel} and \eqref{eq:be} shows
\begin{align*}
\Delta (r^{-\eta}) + v_n^{p-1} r^{-\eta-2+b_n} &\le -\eta (N-2-\eta) r^{-\eta-2} + o(1) r^{-\beta_n(p-1)} r^{-\eta-2+b_n} \\
&\le - \frac{\eta}{2} (N-2-\eta) r^{-\eta-2} < 0
\end{align*}
and
\begin{align*}
\Delta (r^{-\eta-2+b_n}) + u_n^{q_n-1} r^{-\eta} &\le -(\eta+2-b_n)(N-4-\eta+b_n) r^{-\eta-4+b_n} + o(1) r^{-\alpha_n (q_n-1)} r^{-\eta} \\
&\le - \frac12 (\eta+2-b_n)(N-4-\eta+b_n) r^{-\eta-4+b_n} < 0
\end{align*}
for $r_n \le r \le \rho_1$. Hence
\[L_n\(r^{-\eta}, r^{-\eta-2+b_n}\) < 0 \quad \text{for } \eta \in (0,1) \text{ small}.\]
The above computations also yield
\[L_n\(r^{-N+4-b_n+\eta}, r^{-(N-2)+\eta}\) < 0 \quad \text{for } \eta \in (0,1) \text{ small}.\]

At this moment, we divide the cases according to the value of $p$.

\medskip \noindent \textbf{Case 1: $p \in (\frac{N}{N-2}, \frac{N+2}{N-2})$.}
By the maximum principle, \eqref{eq:uedecay}, and an analogous estimate for $r^{\beta_n} v_n$,
\begin{equation}\label{eq:811}
u_n(x) \le M_n r^{-\eta} + Q_n r^{-N+4-b_n+\eta} \quad \text{and} \quad v_n(x) \le M_n r^{-\eta-2+b_n} + Q_n r^{-(N-2)+\eta}
\end{equation}
for $x$ with $r_n \le r=|x-x_n| \le \rho_1$. Here,
\begin{equation}\label{eq:812}
M_n := \(\rho_1^\eta + \rho_1^{2-b_n+\eta}\) \max\left\{\max_{x \in \pa B(x_n,\rho_1)} u_n(x),\, \max_{x \in \pa B(x_n,\rho_1)} v_n(x)\right\}
\end{equation}
and
\begin{equation}\label{eq:813}
\begin{aligned}
Q_n := C\lambda_n^{\beta_n-(N-2-\eta)} R_n^{-\eta}
&= C\max\left\{\lambda_n^{\alpha_n-(N-4+b_n-\eta)}R_n^{-2+b_n-\eta},\, \lambda_n^{\beta_n-(N-2-\eta)}R_n^{-\eta}\right\} \\
&\ge \max\left\{r_n^{N-4+b_n-\eta} \max_{x \in \pa B(x_n,r_n)} u_n(x),\,
r_n^{N-2-\eta} \max_{x \in \pa B(x_n,r_n)} v_n(x)\right\}
\end{aligned}
\end{equation}
for a suitable $C > 0$, where the equality in \eqref{eq:813} follows from \eqref{eq:be} and
\begin{equation}\label{eq:815}
\alpha_n-b_n = \alpha_n-p\beta_n+\beta_n = \beta_n-2.
\end{equation}
Furthermore, from \eqref{eq:Harnack1}--\eqref{eq:Harnack2}, \eqref{eq:uvmono}, and \eqref{eq:811}, we observe
\[\begin{cases}
\displaystyle \rho_1^{\alpha_n} \max_{x \in \pa B(x_n,\rho_1)} u_n(x) \le C \theta^{\alpha_n} \(M_n \theta^{-\eta} + Q_n \theta^{-N+4-b_n+\eta}\),\\
\displaystyle \rho_1^{\beta_n} \max_{x \in \pa B(x_n,\rho_1)} v_n(x) \le C \theta^{\beta_n} \(M_n \theta^{-\eta-2+b_n} + Q_n \theta^{-(N-2)+\eta}\)
\end{cases}
\quad \text{for all } r_n \le \theta \le \rho_1.\]
Choosing small $\theta$ and $\eta$ (independent of $m$) and applying $\beta_np > 2$, we get
\begin{equation}\label{eq:814}
M_n \le C Q_n = C \lambda_n^{\beta_n-(N-2)+\eta} R_n^{-\eta}.
\end{equation}
By plugging \eqref{eq:814} into \eqref{eq:811}, we derive
\[\begin{aligned}
v_n(x) &\le C \lambda_n^{\beta_n-(N-2)+\eta} R_n^{-\eta} \(r^{-\eta-2+b_n}+r^{-(N-2)+\eta}\) \\
&\le C\lambda_n^{\beta_n-(N-2)+\eta} R_n^{-\eta} r^{-(N-2)+\eta}
\end{aligned}
\quad \text{for } r_n \le r \le \rho_1,\]
which implies \eqref{eq:801}.

It remains to verify \eqref{eq:802}. By \eqref{eq:801}, the assumption $p > \frac{N}{N-2}$, and \eqref{eq:815},
\[-\Delta u_n = v_n^p \le C \lambda_n^{(\beta_n-(N-2)+\eta)p} r^{(-(N-2)+\eta)p} \le C \lambda_n^{\alpha_n-(N-2)+\eta} r^{-N+\eta}
\quad \text{for } r_n \le r \le \rho_1.\]
Thus
\[\Delta \(u_n - C \lambda_n^{\alpha_n-(N-2)+\eta} r^{-(N-2)+\eta}\) \ge 0 \quad \text{for } r_n \le r \le \rho_1.\]
Besides, \eqref{eq:uedecay} gives
\[u_n \le C \lambda_n^{\alpha_n-(N-2)} r_n^{2-N} \le C \lambda_n^{\alpha_n-(N-2)+\eta} r_n^{-(N-2)+\eta}.\]
It follows from the maximum principle that
\begin{equation}\label{eq:816}
u_n \le C \(\lambda_n^{\alpha_n-(N-2)+\eta} r^{-(N-2)+\eta} + \bar{u}_n(\rho_1) \rho_1^{\eta} r^{-\eta}\) \quad \text{for } r_n \le r \le \rho_1.
\end{equation}
Arguing as in the derivation of \eqref{eq:814} and applying \eqref{eq:816} at this time, we see
\begin{equation}\label{eq:817}
\bar{u}_n(\rho_1) \le C \lambda_n^{\alpha_n-(N-2)+\eta}.
\end{equation}
Inserting \eqref{eq:817} into \eqref{eq:816}, we deduce \eqref{eq:802}.

\medskip \noindent \textbf{Case 2: $p = \frac{N}{N-2}$.} A slight modification of the argument in Case 1 gives \eqref{eq:801}--\eqref{eq:802}.

\medskip \noindent \textbf{Case 3: $p \in (\frac{2}{N-2}, \frac{N}{N-2})$.}
Arguing as in Case 1, together with \eqref{eq:815} and $b_n \le (N-2)(p-1)$, we see that \eqref{eq:811} holds with $M_n$ in \eqref{eq:812}
and $Q_n = C\lambda_n^{\beta_n-(N-2-\eta)}R_n^{-\eta}$ for some large $C > 0$. From this fact, we further deduce that
\[u_n(x),\, v_n(x) \le C\lambda_n^{\beta_n-(N-2)+\eta} R_n^{-\eta} r^{-(N-2)+\eta} \quad \text{for } r_n \le r \le \rho_1,\]
which in particular implies \eqref{eq:801}.

Let us check \eqref{eq:802}. By \eqref{eq:815},
\[-\Delta u_n = v_n^p \le C \lambda_n^{(\beta_n -(N-2)+\eta)p} r^{-(N-2-\eta)p} = C\lambda_n^{\alpha_n+2 -((N-2)-\eta)p} r^{-(N-2-\eta)p}\]
for $r_n \le r \le \rho_1$. Also, since $p \in (\frac{2}{N-2}, \frac{N}{N-2})$, it holds that
\[-\Delta\(r^{2-(N-2-\eta)p}\) = ((N-2)p-2-\eta p)(N-(N-2)p+\eta p) r^{-(N-2-\eta)p} \ge C r^{-(N-2-\eta)p}.\]
Therefore 
\begin{equation}\label{eq:831}
u_n(x) \le C\(\lambda_n^{\alpha_n+2-(N-2)p+\eta p} r^{2-(N-2)p+\eta p} + \bar{u}_n(\rho_1) \rho_1^\eta r^{-\eta}\) \quad \text{for } r_n \le r \le \rho_1.
\end{equation}
Furthermore, as in \eqref{eq:817}, we can show that
\begin{equation}\label{eq:832}
\bar{u}_n(\rho_1) \le C \lambda_n^{\alpha_n+2-(N-2)p+\eta p}.
\end{equation}
Combining \eqref{eq:831} and \eqref{eq:832}, we conclude that \eqref{eq:802} is true.
\end{proof}

\begin{lemma}\label{lemma:2}
In the setting of Lemma \ref{lemma:8}, we reduce the value of $\rho_1$ if needed.
Then there exists a constant $C = C(\rho) > 0$ independent of $n \in \N$ such that
\begin{equation}\label{eq:21}
\limsup_{n \to \infty} \max_{x \in \pa B(x_n,\rho)} \lambda_n^{N-2-\beta_n} v_n(x) \le C(\rho)
\end{equation}
for each $\rho \in (0,\frac{\rho_1}{2})$.
\end{lemma}
\begin{proof}
We infer from the comparison principle and \eqref{eq:72} that
\begin{equation}\label{eq:22}
v_n(x) \ge C \lambda_n^{\beta_n -(N-2)} \(|x-x_n|^{2-N}-\rho_1^{2-N}\) \quad \text{in } B(x_n,\rho_1) \setminus B(x_n,r_n).
\end{equation}
By \eqref{eq:802}, \eqref{eq:22}, \eqref{eq:pqrel2}, and \eqref{eq:72},
\begin{equation}\label{eq:23}
u_n^{\beta_n} \le Cv_n^{\alpha_n} \quad \text{in } \overline{B\(x_n,\frac{\rho_1}{2}\)}
\end{equation}
for all $n \in \N$.

Given a fixed $\rho \in (0,\frac{\rho_1}{2})$, we write
\[(\hat{u}_n, \hat{v}_n) = \(\hat{\lambda}_n^{-\alpha_n} u_n, \hat{\lambda}_n^{-\beta_n} v_n\)
\quad \text{where } \hat{\lambda}_n := \max_{x \in \pa B(x_n,\rho)} v_n^{\frac{1}{\beta_n}}(x) \to 0 \quad \text{as } n \to \infty.\]
Then, by \eqref{eq:815} and $\alpha_nq_n-\beta_n=2$,
\[-\Delta \hat{u}_n = \hat{\lambda}_n^2 \hat{v}_n^p \quad \text{and} \quad -\Delta \hat{v}_n = \hat{\lambda}_n^2 \hat{u}_n^{q_n} \quad \text{in } B(x_n,\rho_1).\]
In light of the Harnack inequality \eqref{eq:Harnack1}--\eqref{eq:Harnack2} and \eqref{eq:23}, if $K$ is a compact subset of $\overline{B(x_n,\rho_1/2)} \setminus \{\xi\}$, there exists $C = C(K) > 1$ such that
\[C(K)^{-1} \le \hat{v}_n \le C(K) \quad \text{and} \quad \hat{u}_n \le C \hat{v}_n^{\alpha_n \over \beta_n} \le C(K) \quad \text{on } K.\]
Hence standard elliptic estimates yield
\[\hat{v}_n \to v_0 \quad \text{in } C^2_{\text{loc}}\(\overline{B\(x_n,\frac{\rho_1}{2}\)} \setminus \{\xi\}\)\]
passing to a subsequence, so that $v_0$ is harmonic and positive in $B(x_n,\rho_1/2) \setminus \{\xi\}$.
Also, by \eqref{eq:uvmono}, we have that $(r^{\beta_0} \bar{v}_0(r))' < 0$ for $0 < r \le \frac{\rho_1}{2}$, so $v_0$ is singular at $0$. It follows from B\^{o}cher's theorem that
\[-\int_{B(x_n,\rho)} \Delta \hat{v}_n = -\int_{\pa B(x_n,\rho)} \nabla \hat{v}_n \cdot \nu = -\int_{\pa B(\xi,\rho)} \nabla v_0 \cdot \nu +o(1) > c\]
for some number $c > 0$, and so
\[c < -\int_{B(x_n,\rho)} \Delta \hat{v}_n = \int_{B(x_n,\rho)} \hat{\lambda}_n^{-\beta_n} u_n^{q_n}
\le C\hat{\lambda}_n^{-\beta_n} \lambda_n^{\alpha_n q_n-N} = C\hat{\lambda}_n^{-\beta_n} \lambda_n^{\beta_n-(N-2)}\]
where we used \eqref{eq:72} and \eqref{eq:802} for the last inequality. As a result, \eqref{eq:21} holds.
\end{proof}

The following is the main result of this subsection.
\begin{prop}\label{prop:uvdecay}
In the setting of Lemma \ref{lemma:8}, we reduce the value of $\rho_1$ if needed.
Then there is a constant $C > 0$ independent of $n \in \N$ such that
\begin{equation}\label{eq:901}
v_n(x) \le C \lambda_n^{\beta_n-(N-2)} |x-x_n|^{2-N}
\end{equation}
and
\begin{equation}\label{eq:902}
u_n(x) \le \begin{cases}
C\lambda_n^{\alpha_n-(N-2)} |x-x_n|^{2-N} &\text{if } p \in (\frac{N}{N-2}, \frac{N+2}{N-2}),\\
C\lambda_n^{\alpha_n-(N-2)} |x-x_n|^{2-N} \log (\lambda_n |x-x_n|+2) &\text{if } p = \frac{N}{N-2},\\
C\lambda_n^{\alpha_n+2-(N-2)p} |x-x_n|^{2-(N-2)p} &\text{if } p \in (\frac{2}{N-2}, \frac{N}{N-2})
\end{cases}
\end{equation}
for $|x-x_n| \le \frac{\rho_1}{3}$.

Furthermore, there exist numbers $a,\, b > 0$ and $a' < 0$ such that
\begin{equation}\label{eq:903}
\begin{cases}
\lambda_n^{N-2-\alpha_n}u_n \to a|\cdot-\xi|^{2-N} + h_1 &\text{if } p \in (\frac{N}{N-2}, \frac{N+2}{N-2}),\\
\lambda_n^{N-2-\alpha_n}(\log\lambda_n)^{-1} u_n \to a|\cdot-\xi|^{2-N} + h_1 &\text{if } p = \frac{N}{N-2},\\
\lambda_n^{(N-2)p-2-\alpha_n}u_n \to a|\cdot-\xi|^{2-(N-2)p} + a'|\cdot-\xi|^{N-(N-2)p}h_2 + h_1 &\text{if } p \in [\frac{N-1}{N-2}, \frac{N}{N-2}),\\
\lambda_n^{(N-2)p-2-\alpha_n}u_n \to a|\cdot-\xi|^{2-(N-2)p} + h_1 &\text{if } p \in (\frac{2}{N-2}, \frac{N-1}{N-2})
\end{cases}
\end{equation}
and
\begin{equation}\label{eq:904}
\lambda_n^{N-2-\beta_n}v_n \to b|\cdot-\xi|^{2-N} + h_2
\end{equation}
in $C^2_{\textnormal{loc}}\(B(\xi,\frac{\rho_1}{3}) \setminus \{\xi\}\)$ as $n \to \infty$, up to a subsequence. Here,
\begin{itemize}
\item[-] if $p \in [\frac{N}{N-2}, \frac{N+2}{N-2})$, then $h_1 \in C^{\infty}(B(\xi,\frac{\rho_1}{3}))$ is a harmonic function.
   If $p \in (\frac{2}{N-2}, \frac{N}{N-2})$, then $h_1 \in C^{1,\sigma}(B(\xi,\frac{\rho_1}{3}))$ for some $\sigma \in (0,1)$;
\item[-] $h_2 \in C^{\infty}(B(\xi,\frac{\rho_1}{3}))$ is a harmonic function;
\item[-] we have
\begin{equation}\label{eq:909}
\begin{cases}
pb^{p-1} = a'[(N-2)p-2(N-1)][N-(N-2)p] &\text{for } p \in [\frac{N-1}{N-2}, \frac{N}{N-2}),\\
b^p = a[(N-2)p-2][N-(N-2)p]&\text{for } p \in (\frac{2}{N-2}, \frac{N}{N-2}).
\end{cases}
\end{equation}
\end{itemize}
\end{prop}
\begin{proof}
First, we prove \eqref{eq:901}. If it is not true, then there exists $\tx_n \in B(x_n,\frac{\rho_1}{3})$ such that
\[v_n(\tx_n) \lambda_n^{N-2-\beta_n}|\tx_n-x_n|^{N-2} \to \infty \quad \text{as } n \to \infty,\]
along a subsequence. By \eqref{eq:72}, we have that $r_n < |\tx_n-x_n| \le \frac{\rho_1}{3}$.
If we set $\tr_n = |\tx_n -x_n|$, then $0$ is an isolated simple blow-up point of
$\{(\tr_n^{\alpha_n} u_n(\tr_n\cdot + x_n), \tr_n^{\beta_n} v_n(\tr_n\cdot + x_n))\}$. By \eqref{eq:21},
\[\max_{x \in \pa B(x_n,1)} (\lambda_n \tr_n)^{N-2-\beta_n} \tr_n^{\beta_n} v_n(\tr_nx + x_n)
= \max_{x \in \pa B(x_n,\tr_n)} \lambda_n^{N-2-\beta_n} \tr_n^{N-2} v_n(x) < \infty,\]
which contradicts the assumption. Thus, the assertion follows.

\medskip
In the rest of the proof, we will verify the claims \eqref{eq:902}, \eqref{eq:903}, and \eqref{eq:904} by splitting the cases according to the value of $p$.

\medskip \noindent \textbf{Case 1: $p \in (\frac{N}{N-2}, \frac{N+2}{N-2})$.}
In this case, by arguing as in the proof of Lemma \ref{lemma:2}, we observe
\[\limsup_{n \to \infty} \max_{x \in \pa B(x_n,\rho)} \lambda_n^{N-2-\alpha_n} u_n(x) \le C(\rho)\]
for all $\rho \in (0,\frac{\rho_1}{2})$. Using this, we can imitate the proof of \eqref{eq:901} to deduce \eqref{eq:902}.

We have
\[\begin{cases}
-\Delta\(\lambda_n^{N-2-\alpha_n} u_n\) = \lambda_n^{N-(N-2)p} \(\lambda_n^{N-2-\beta_n}v_n\)^p,\\
-\Delta\(\lambda_n^{N-2-\beta_n} v_n\) = \lambda_n^{N-(N-2)q_n} \(\lambda_n^{N-2-\alpha_n}u_n\)^{q_n}.
\end{cases}\]
From the above equations, \eqref{eq:901}--\eqref{eq:902}, and elliptic regularity, we see that
\[\lambda_n^{N-2-\alpha_n}u_n \to w_1 \quad \text{and} \quad \lambda_n^{N-2-\beta_n}v_n \to w_2
\quad \text{in } C^2_{\text{loc}}\(B\(\xi,\frac{\rho_1}{3}\) \setminus \{\xi\}\)\]
after passing to a subsequence, where $w_1$ and $w_2$ are harmonic in $B(\xi,\frac{\rho_1}{3}) \setminus \{\xi\}$.
By the comparison principle, we obtain \eqref{eq:22} and
\[u_n(x) \ge C\lambda_n^{\alpha_n-(N-2)}\(|x-x_n|^{2-N}-\rho_1^{2-N}\) \quad \text{in } B(x_n,\rho_1) \setminus B(x_n,r_n),\]
so $w_1, w_2>0$. It follows from B\^{o}cher's theorem that
\[w_1 = a|\cdot-\xi|^{2-N}+h_1 \quad \text{and} \quad w_2 = b|\cdot-\xi|^{2-N} + h_2 \quad \text{in } B\(\xi,\frac{\rho_1}{3}\) \setminus \{\xi\}\]
where $h_1$ and $h_2$ are harmonic in $B(\xi,\frac{\rho_1}{3})$. Also, by \eqref{eq:uvmono}, $w_1$ and $w_2$ are singular at $\xi$, which reads $a, b > 0$.
As a consequence, \eqref{eq:903} and \eqref{eq:904} are true.

\medskip \noindent \textbf{Case 2: $p = \frac{N}{N-2}$.}
According to \eqref{eq:72},
\[u_n(x) \le C \lambda_n^{\alpha_n-(N-2)} |x-x_n|^{2-N}\log(\lambda_n |x-x_n|+2) \quad \text{for } x \in B(x_n,r_n).\]
Since $-\Delta (r^{2-N}\log(\lambda_n r)) = (N-2) r^{-N}$, the comparison argument, \eqref{eq:901}, and \eqref{eq:uvmono} show
\[u_n(x) \le C \lambda_n^{\alpha_n-(N-2)} |x-x_n|^{2-N}\log(\lambda_n |x-x_n|) \quad \text{for } x \in B\(x_n,\frac{\rho_1}{3}\) \setminus B(x_n,r_n),\]
and we deduce \eqref{eq:902}.

We have
\[\begin{cases}
-\Delta\left[\lambda_n^{N-2-\alpha_n} (\log \lambda_n)^{-1} u_n\right] = \lambda_n^{N-(N-2)p} (\log \lambda_n)^{-1} \(\lambda_n^{N-2-\beta_n}v_n\)^p,\\
-\Delta\(\lambda_n^{N-2-\beta_n} v_n\) = \lambda_n^{N-(N-2)q_n} (\log \lambda_n)^{q_n} \left[\lambda_n^{N-2-\alpha_n}(\log \lambda_n)^{-1} u_n\right]^{q_n}.
\end{cases}\]
From the above equations, \eqref{eq:901}--\eqref{eq:902}, and elliptic regularity, we see that
\[\lambda_n^{N-2-\alpha_n}(\log\lambda_n)^{-1} u_n \to w_1 \quad \text{and} \quad \lambda_n^{N-2-\beta_n}v_n \to w_2
\quad \text{in } C^2_{\text{loc}}\(B\(\xi,\frac{\rho_1}{3}\) \setminus \{\xi\}\)\]
after passing to a subsequence, where $w_1$ and $w_2$ are harmonic in $B(\xi,\frac{\rho_1}{3}) \setminus \{\xi\}$. By the comparison principle, we obtain \eqref{eq:22} and
\[u_n(x) \ge C\lambda_n^{\alpha_n-(N-2)} \left[|x-x_n|^{2-N}\log(\lambda_n |x-x_n|) - \rho_1^{2-N} \log(\lambda_n \rho_1)\right]
\quad \text{in } B(x_n,\rho_1) \setminus B(x_n,r_n),\]
so $w_1, w_2 > 0$. Reasoning as in Case 1, we conclude that \eqref{eq:903} and \eqref{eq:904} hold.

\medskip \noindent \textbf{Case 3: $p \in (\frac{2}{N-2}, \frac{N-1}{N-2})$.}
According to \eqref{eq:72},
\[u_n(x) \le C \lambda_n^{\alpha_n+2-(N-2)p} |x-x_n|^{2-(N-2)p} \quad \text{for } x \in B(x_n,r_n).\]
Employing this inequality, the comparison argument, and \eqref{eq:uvmono}, we deduce \eqref{eq:902}.

We have
\begin{equation}\label{eq:905}
\begin{cases}
-\Delta\(\lambda_n^{(N-2)p-2-\alpha_n} u_n\) = \(\lambda_n^{N-2-\beta_n}v_n\)^p,\\
-\Delta\(\lambda_n^{N-2-\beta_n} v_n\) = \lambda_n^{N-(N-2)pq_n+2q_n} \(\lambda_n^{(N-2)p-2-\alpha_n}u_n\)^{q_n}.
\end{cases}
\end{equation}
Note that
\[N-(N-2)pq_n+2q_n = -2p-2+\e_n(p+1)(q_n+1) < 0.\]
Therefore, arguing as in Case 1, we obtain \eqref{eq:904}.
Also, by \eqref{eq:905}, \eqref{eq:901}--\eqref{eq:902}, elliptic regularity, and the assumption $p < \frac{N}{N-2}$,
\begin{equation}\label{eq:906}
\lambda_n^{(N-2)p-2-\alpha_n} u_n \to w_1 \quad \text{in } C^2_{\text{loc}}\(B\(\xi,\frac{5\rho_1}{12}\) \setminus \{\xi\}\),
\end{equation}
up to a subsequence, and $w_1$ satisfies
\[-\Delta w_1 = \(\frac{b}{|\cdot-\xi|^{N-2}} + h_2\)^p \quad \text{in } B\(\xi,\frac{5\rho_1}{12}\).\]
From \eqref{eq:909}, we see that
\begin{equation}\label{eq:907}
-\Delta\(w_1 - \frac{a}{|\cdot-\xi|^{(N-2)p-2}}\) = \(\frac{b}{|\cdot-\xi|^{N-2}} + h_2\)^p - \frac{b^p}{|\cdot-\xi|^{(N-2)p}} \quad \text{in } B\(\xi,\frac{5\rho_1}{12}\).
\end{equation}
Since $p < \frac{N-1}{N-2}$, the right-hand side of \eqref{eq:907} belongs to $L^{N+\sigma'}(B(\xi,\frac{5\rho_1}{12}))$ for some $\sigma' \in (0,1)$.
By elliptic regularity, \eqref{eq:903} holds.

\medskip \noindent \textbf{Case 4: $p \in [\frac{N-1}{N-2}, \frac{N}{N-2})$.}
In this case, the proof goes along the same lines as that of Case 3 until the derivation of \eqref{eq:906}.
If we set $a$ and $a'$ as in \eqref{eq:909}, then
\begin{equation}\label{eq:908}
\begin{aligned}
&\ -\Delta\(w_1(x) - \frac{a}{|x-\xi|^{(N-2)p-2}} - a'|x-\xi|^{N-(N-2)p}h_2(x)\) \\
&= \left[\(\frac{b}{|x-\xi|^{N-2}} + h_2\)^p - \frac{b^p}{|x-\xi|^{(N-2)p}} - \frac{pb^{p-1}h_2(x)}{|x-\xi|^{(N-2)(p-1)}}\right] \\
&\ - \frac{2pb^{p-1} \nabla h_2(x)}{2(N-1)-(N-2)p} \cdot \frac{x-\xi}{|x-\xi|^{(N-2)(p-1)}}
\end{aligned}
\quad \text{for } x \in B\(\xi,\frac{5\rho_1}{12}\)
\end{equation}
where we used the fact that $h_2$ is harmonic in $B(\xi,\frac{5\rho_1}{12})$. Since $p < \frac{N}{N-2}$,
the right-hand side of \eqref{eq:908} belongs to $L^{N+\sigma'}(B(\xi,\frac{5\rho_1}{12}))$ for some $\sigma' \in (0,1)$.
By elliptic regularity, \eqref{eq:903} holds.
\end{proof}

\subsection{Exclusion of bubble accumulation}
We will prove that each blow-up point of $\{(u_n,v_n)\}_{n \in \N}$ is isolated simple and distant from the other ones. The following result is a key ingredient of the proof.
\begin{lemma}\label{lemma:10}
Let $x_n \to \xi$ be an isolated simple blow-up point of $\{(u_n,v_n)\}_{n \in \N}$, and $h_1$ and $h_2$ be the functions in Proposition \ref{prop:uvdecay}.
Assume that $h_1(x) = A+O(|x-\xi|)$ and $h_2(x) = B+O(|x-\xi|)$ as $x \to \xi$.

\medskip \noindent
\textnormal{(1)} If $p \in [\frac{N}{N-2}, \frac{N+2}{N-2})$, then
\begin{equation}\label{eq:100}
\frac{aB}{p+1} + \frac{bA}{q_0+1} \le 0.
\end{equation}

\medskip \noindent
\textnormal{(2)} If $p \in (\frac{2}{N-2}, \frac{N}{N-2})$, then
\begin{equation}\label{eq:101}
A \le 0.
\end{equation}
\end{lemma}
\begin{proof}
(1) One can prove it by following the argument in \cite[Proposition 3.2]{LZ}.

\medskip \noindent (2) Fix $\rho \in (0,\frac{\rho_1}{3})$. By Pohozaev identity \eqref{eq:poho} and \eqref{eq:LE}, we have
\begin{multline}\label{eq:poho2}
\rho \int_{\pa B(x_n,\rho)} \(2 \frac{\pa u_n}{\pa \nu} \frac{\pa v_n}{\pa \nu}-\nabla u_n \cdot \nabla v_n\) dS
+ \rho \int_{\pa B(x_n,\rho)} \(\frac{v_n^{p+1}}{p+1} + \frac{u_n^{q_n+1}}{q_n+1}\) dS \\
+ N\int_{\pa B(x_n,\rho)} \left[\(\frac{1}{p+1} - \frac{\e_n}{N}\) v_n \frac{\pa u_n}{\pa \nu} + \frac{u_n}{q_n+1}\frac{\pa v_n}{\pa \nu}\right] dS
= \e_n \int_{B(x_n,\rho)} v_n^{p+1} \ge 0.
\end{multline}
We will derive \eqref{eq:101} by plugging the estimates in Proposition \ref{prop:uvdecay} into \eqref{eq:poho2}.

First, applying \eqref{eq:903}--\eqref{eq:904}, we compute
\begin{multline}\label{eq:102}
\lim_{n \to \infty} \lambda_n^{(N-2)(p+1)-2-(\alpha_n+\beta_n)} \rho \int_{\pa B(x_n,\rho)} \(2 \frac{\pa u_n}{\pa \nu} \frac{\pa v_n}{\pa \nu}-\nabla u_n \cdot \nabla v_n\) dS \\
= ab((N-2)p-2)(N-2) \left|\S^{N-1}\right| \rho^{2-(N-2)p} + O\(\rho^{N-(N-2)p} + \rho\)
\end{multline}
as $\rho \to 0$. Second, from \eqref{eq:pqrel2} and \eqref{eq:815}, we obtain
\begin{multline}\label{eq:103}
\lim_{n \to \infty} \lambda_n^{(N-2)(p+1)-2-(\alpha_n+\beta_n)} \rho \int_{\pa B(x_n,\rho)} \(\frac{v_n^{p+1}}{p+1} + \frac{u_n^{q_n+1}}{q_n+1}\) dS \\
= \left[\frac{b^{p+1}}{p+1} \left|\S^{N-1}\right| \rho^{2-(N-2)p} + O\(\rho^{N-(N-2)p}\)\right] + 0.
\end{multline}
Lastly, we see
\begin{equation}\label{eq:104}
\begin{aligned}
&\ \lim_{n \to \infty} \lambda_n^{(N-2)(p+1)-2-(\alpha_n+\beta_n)} N\int_{\pa B(x_n,\rho)}
\left[\(\frac{1}{p+1} - \frac{\e_n}{N}\) v_n \frac{\pa u_n}{\pa \nu} + \frac{u_n}{q_n+1}\frac{\pa v_n}{\pa \nu}\right] dS \\
&= -\left|\S^{N-1}\right| \frac{abN}{p+1}((N-2)p-2) \rho^{2-(N-2)p} \\
&\ -\left|\S^{N-1}\right| \left[\frac{abN}{q_0+1}(N-2) \rho^{2-(N-2)p} + \frac{bN}{q_0+1}(N-2)A\right] + O\(\rho^{N-(N-2)p} + \rho\).
\end{aligned}
\end{equation}
By \eqref{eq:cr-hy} and \eqref{eq:909},
\begin{equation}\label{eq:105}
\begin{aligned}
ab((N-2)p-2)(N-2) + \frac{b^{p+1}}{p+1} - \left[\frac{abN}{p+1}((N-2)p-2) + \frac{abN}{q_0+1}(N-2)\right] = 0.
\end{aligned}
\end{equation}
Putting \eqref{eq:102}--\eqref{eq:105} into \eqref{eq:poho2}, we derive
\begin{equation}\label{eq:106}
O\(\rho^{N-(N-2)p} + \rho\) - \frac{bN}{q_0+1}(N-2)\left|\S^{N-1}\right| A \ge 0.
\end{equation}
Finally, invoking $p < \frac{N}{N-2}$ and taking $\rho \to 0$ in \eqref{eq:106}, we establish \eqref{eq:101}.
\end{proof}

\begin{remark}\label{rmk:10}
Using \eqref{eq:poho2} with $\rho = \frac{\rho_1}{3}$ and
\[\e_n \int_{B(x_n,\frac{\rho_1}{3})} v_n^{p+1} \ge C \e_n \lambda_n^{\beta_n(p+1)-N} \int_{B(0,\frac{\lambda_n \rho_1}{3})} V_{1,0}^{p+1} \ge C \e_n \lambda_n^{\beta_n(p+1)-N},\]
one can derive
\begin{equation}\label{eq:10}
\e_n \le \begin{cases}
C \lambda_n^{2-N} &\text{if } p \in (\frac{N}{N-2}, \frac{N+2}{N-2}),\\
C \lambda_n^{2-N} \log \lambda_n &\text{if } p = \frac{N}{N-2},\\
C \lambda_n^{2-(N-2)p} &\text{if } p \in (\frac{2}{N-2}, \frac{N}{N-2}).
\end{cases}
\end{equation}
In particular, $\lambda_n^{\e_n} \to 1$ as $n \to \infty$, and one can replace each $\alpha_n$ and $\beta_n$ in \eqref{eq:901}--\eqref{eq:904} with $\alpha_0$ and $\beta_0$, respectively.
\end{remark}

\begin{prop}\label{prop:isos}
We reduce the value of $\rho_1$ if needed. Then every isolated blow-up point of $\{(u_n,v_n)\}_{n \in \N}$ is in fact isolated simple.
\end{prop}
\begin{proof}
Let $x_n \to \xi$ be an isolated blow-up point of $\{(u_n,v_n)\}$.
We assert that each map $r \mapsto r^{\alpha_n} \bar{u}_n$ and $r \mapsto r^{\beta_n} \bar{v}_n$
has only one critical point in the interval $(0,\rho_2)$ for some $\rho_2 \in (0, \frac{\rho_1}{3})$.

Suppose not. In view of Remark \ref{rmk:7}, after passing to a subsequence of $(u_n,v_n)$,
either the map $r \mapsto r^{\alpha_n} \bar{u}_n$ has the second critical point $\zeta_{1n} \ge R_n\lambda_n^{-1}$
or the map $r \mapsto r^{\beta_n} \bar{v}_n$ has the second critical point $\zeta_{2n} \ge R_n\lambda_n^{-1}$,
and $\zeta_n := \min\{\zeta_{1n}, \zeta_{2n}\} \to 0$ as $n \to \infty$.\footnote{If $r \mapsto r^{\alpha_n} \bar{u}_n$ has no critical point on $[R_n\lambda_n^{-1}, \frac{\rho_1}{3})$, then we set $\zeta_{1n} = \infty$.
Similarly, if $r \mapsto r^{\beta_n} \bar{v}_n$ has no critical point on $[R_n\lambda_n^{-1}, \frac{\rho_1}{3})$, then we set $\zeta_{2n} = \infty$.} Let
\[(w_n,z_n)(y) = \(\zeta_n^{\alpha_n} u_n(\zeta_n y + x_n),\,
\zeta_n^{\beta_n} v_n(\zeta_ny + x_n)\) \quad \text{for } |y| \le \frac{\rho_1}{3\zeta_n}.\]
It satisfies
\[\begin{cases}
-\Delta w_n = z_n^p,\ -\Delta z_n = w_n^{q_n} &\text{in } B\(0,\frac{\rho_1}{3\zeta_n}\),\\
w_n, z_n > 0 &\text{in } B\(0,\frac{\rho_1}{3\zeta_n}\),\\
w_n(y) \le C|y|^{-\alpha_n},\ z_n(y) \le C |y|^{-\beta_n} &\text{for } y \in B\(0,\frac{\rho_1}{3\zeta_n}\).
\end{cases}\]
Furthermore, each map $r \mapsto r^{\alpha_n} \bar{w}_n$ and $r \mapsto r^{\beta_n} \bar{z}_n$ has precisely one critical point in the interval $(0,1)$, and
\begin{equation}\label{eq:isos1}
\text{either } \left. \frac{d}{dr} [r^{\alpha_n} w_n(r)]\right|_{r=1} = 0
\quad \text{or} \quad \left. \frac{d}{dr} [r^{\beta_n} z_n(r)]\right|_{r=1} = 0.
\end{equation}
In particular, $0 \in \R^N$ is an isolated simple blow-up point of $\{(w_n,z_n)\}_{n \in \N}$.

By Proposition \ref{prop:uvdecay}, \eqref{eq:903} and \eqref{eq:904} hold in $C^2_{\text{loc}}(\R^N \setminus \{0\})$
if $(u_n,v_n)$, $\lambda_n$, and $\xi$ are replaced with $(w_n,z_n)$, $\tla_n := w_n^{1/\alpha_n}(x_n)$, and $0$, respectively. In addition,
\[\begin{cases}
w_n(0) = \zeta_n^{\alpha_n} u_n(x_n) = (\zeta_n \lambda_n)^{\alpha_n} \ge R_n^{\alpha_n} \to \infty,\\
z_n(0) = \zeta_n^{\beta_n} v_n(x_n) \ge C \zeta_n^{\beta_n} u_n^{\beta_n/\alpha_n}(x_n)
= C(\zeta_n \lambda_n)^{\beta_n} \ge R_n^{\beta_n} \to \infty \quad \text{(by \eqref{eq:23})}
\end{cases}\]
as $n \to \infty$. Therefore the numbers $a$ and $b$ in \eqref{eq:903}--\eqref{eq:904} are positive and $\liminf_{|x| \to \infty} h_2(x)\ge 0$.
Since $h_2$ is harmonic in $\R^N$, it follows that $h_2(x) = B \ge 0$ for all $x \in \R^N$; see Lemma \ref{lemma:10}.

At this moment, we divide the cases according to the value of $p$.

\medskip \noindent \textbf{Case 1: $p \in [\frac{N}{N-2}, \frac{N+2}{N-2})$.}
Since $h_1$ is harmonic in $\R^N$, we have that $h_1(x) = A \ge 0$ for all $x \in \R^N$; see Lemma \ref{lemma:10}.

Suppose that $\zeta_n = \zeta_{1n}$ for infinitely many $n \in \N$. Then the first equality of \eqref{eq:isos1} holds for such $m$'s, and so
\[0 = \left.\frac{d}{dr} \left[r^{\alpha_0} \(ar^{2-N} + h_1\)\right]\right|_{r=1} = - a(N-2-\alpha_0) + A\alpha_0.\]
Owing to \eqref{eq:pqrel}, it holds that $A > 0$. From \eqref{eq:100}, it follows that
\[0 < \frac{aB}{p+1}+\frac{bA}{q_0+1} \le 0,\]
a contradiction. Consequently, $\zeta_n = \zeta_{2n}$ for all but finitely many $n \in \N$.
Using the second equality of \eqref{eq:isos1}, we see that this cannot happen either. The assertion must be true.

\medskip \noindent \textbf{Case 2: $p \in (\frac{2}{N-2}, \frac{N-1}{N-2})$.}
By the first equation of \eqref{eq:905} (with a suitable change of notations) and \eqref{eq:909},
\begin{equation}\label{eq:isos2}
-\Delta h_1(x) = \(b|x|^{2-N} + B\)^p - b^p|x|^{-(N-2)p} \ge 0 \quad \text{for } x \in \R^N
\end{equation}
and $\liminf_{|x| \to \infty} h_1(x) \ge 0$.
Thus the strong maximum principle yields that $h_1 \ge 0$ in $\R^N$, and either $A = h_1(0) > 0$ or $h_1 = 0$ in $\R^N$.
In light of \eqref{eq:101}, the first possibility cannot happen, and so $h_1 = 0$ in $\R^N$.
If the first equality of \eqref{eq:isos1} holds for infinitely many $n \in \N$, then
\[0 = \left.\frac{d}{dr} \left[r^{\alpha_0} \(ar^{2-(N-2)p} + 0\)\right]\right|_{r=1} = - a((N-2)p-2-\alpha_0) < 0 \quad \text{(by \eqref{eq:pqrel})},\]
a contradiction. If the second equality of \eqref{eq:isos1} holds for infinitely many $n \in \N$, then
\[0 = \left.\frac{d}{dr} \left[r^{\beta_0} \(br^{2-N} + B\)\right]\right|_{r=1} = - b(N-2-\beta_0) + \beta_0 B,\]
from which we obtain that $B > 0$. However, it is again absurd, because the left-hand side of \eqref{eq:isos2} is 0, while the right-hand side is nonzero.
As a result, the assertion holds.

\medskip \noindent \textbf{Case 3: $p \in [\frac{N-1}{N-2}, \frac{N}{N-2})$.}
We have
\begin{equation}\label{eq:isos3}
-\Delta h_1(x) = \(b|x|^{2-N}+B\)^p - b^p|x|^{-(N-2)p} - pb^{p-1}|x|^{-(N-2)(p-1)}B \ge 0 \quad \text{for } x \in \R^N
\end{equation}
where we used $p \ge 1$ to get the inequality. Moreover, since $a > 0$, $a' < 0$, and $h_2 = B \ge 0$ in $\R^N$,
\[\text{either } \left[B > 0 \text{ and } \liminf_{|x| \to \infty} h_1(x) = +\infty\right] \text{ or } B = 0.\]
The former situation cannot happen in view of \eqref{eq:isos3} and the maximum principle. The latter situation is also impossible because of the argument in Case 2.
The assertion now follows.
\end{proof}

\begin{lemma}
Let $\{(u_n,v_n)\}_{n \in \N}$ be a sequence of solutions of \eqref{eq:LE} with $q = q_n$.
Given any $\eta > 0$ small and $R > 0$ large, there exists a constant $C > 0$ depending on $\eta$ and $R$ such that if
\[\max_{x \in \Omega} \max\left\{u_n^{1/\alpha_n}(x),v_n^{1/\beta_n}(x)\right\} \ge C,\]
then the followings hold for some $k \in \N$ determined by the energy condition \eqref{eq:main0}, up to a subsequence:

\medskip \noindent \textnormal{(1)} There exists a set $\{x_{n,1}, \ldots, x_{n,k}\}$ of local maxima of $u_n$
such that $\{B(x_{n,i},R\lambda_{n,i}^{-1})\}_{i=1}^k$ is a disjoint collection of subsets in $\Omega$. Here, $\lambda_{n,i} := u_n^{1/\alpha_n}(x_{n,i})$.

\medskip \noindent \textnormal{(2)} For $i = 1, \ldots, k$, it holds that
\[\left\|\lambda_{n,i}^{-\alpha_n} u_n(\lambda_{n,i}^{-1}\cdot + x_{n,i}) - U_{1,0}\right\|_{C^2(\overline{B(0,R)})}
+ \left\|\lambda_{n,i}^{-\beta_n} v_n(\lambda_{n,i}^{-1}\cdot + x_{n,i}) - V_{1,0}\right\|_{C^2(\overline{B(0,R)})} \le \eta.\]

\medskip \noindent \textnormal{(3)} We have that
\begin{equation}\label{eq:41}
\(\min_{i = 1, \ldots, k} |x-x_{n,i}|^{\alpha_n}\) u_n(x) + \(\min_{i = 1, \ldots, k} |x-x_{n,i}|^{\beta_n}\) v_n(x) \le C \quad \text{for all } x \in \Omega.
\end{equation}
\medskip \noindent \textnormal{(4)} We have that
\begin{equation}\label{eq:42}
\frac{u_n^{1/\alpha_n}(x_{n,i})}{v_n^{1/\beta_n}(x_{n,i})} \to U_{1,0}^{1/\alpha_0}(0) \text{ or } V_{1,0}^{-1/\beta_0}(0) \quad \text{as } n \to \infty \text{ for each } i = 1, \ldots, k.
\end{equation}
\end{lemma}
\begin{proof}
By suitably modifying the argument in the proof of \cite[Lemma 5.1]{LZ} and applying the energy condition \eqref{eq:main0}, one can prove the next claim:
Given any $\eta > 0$ small and $R > 0$ large, there exists a constant $C > 0$ depending on $\eta$ and $R$ such that
if $\{K_n\}_{n \in \N}$ is a sequence of compact subsets of $\Omega$ and $\{(u_n,v_n)\}_{n \in \N}$ is a sequence of solutions of \eqref{eq:LE} with $q = q_n$ satisfying
\[\max_{x \in \Omega \setminus K_n} \max\left\{\dist(x,K_n) u_n^{1/\alpha_n}(x), \dist(x,K_n) v_n^{1/\beta_n}(x)\right\} \ge C,\]
then, along a subsequence, either
\[\left\|\lambda_{1n}^{-\alpha_n} u_n(\lambda_{1n}^{-1}\cdot + x_{1n}) - U_{1,0}\right\|_{C^2(\overline{B(0,R)})}
+ \left\|\lambda_{1n}^{-\beta_n} v_n(\lambda_{1n}^{-1}\cdot + x_{1n}) - V_{1,0}\right\|_{C^2(\overline{B(0,R)})} \le \eta\]
for a local maximum point $x_{1n}$ of $u_n$ in $\Omega \setminus K_n$, or
\[\left\|\lambda_{2n}^{-\alpha_n} u_n(\lambda_{2n}^{-1}\cdot + x_{2n}) - U_{1,0}\right\|_{C^2(\overline{B(0,R)})}
+ \left\|\lambda_{2n}^{-\beta_n} v_n(\lambda_{2n}^{-1}\cdot + x_{2n}) - V_{1,0}\right\|_{C^2(\overline{B(0,R)})} \le \eta\]
for a local maximum point $x_{2n}$ of $v_n$ in $\Omega \setminus K_n$.
Here, $\lambda_{1n} := u_n^{1/\alpha_n}(x_{1n})$, $\lambda_{2n} := v_n^{1/\beta_n}(x_{2n})$, and $\dist(x,\emptyset) := 1$
where $\emptyset$ denotes the empty set.\footnote{If a classification theorem \cite{CLO, CL} for the critical Lane-Emden system in $\R^N$ continues to hold under a local integrability condition,
we will not need \eqref{eq:main0} in Theorems \ref{thm:main1}--\ref{thm:main3}.
However, because system \eqref{eq:LERn} is not invariant under the Kelvin transform, deducing it is a very challenging problem.
For related results, see the references \cite{So, CL, CH} and others devoted to the Lane-Emden conjecture.}

Applying the above claim, we follow the proof of \cite[Proposition 5.1]{LZ}. Then we can obtain the desired result.
\end{proof}

\begin{prop}\label{prop:dist}
Let $x_{n,i} \to \xi_i$ be blow-up points of $\{(u_n,v_n)\}_{n \in \N}$ for $i = 1, \ldots, k$ where $k \in \N$.
Then there exists $\eta_0 > 0$ independent of $n \in \N$ such that $|x_{n,i}-x_{n,j}| \ge \eta_0$ for all $1 \le i \ne j \le k$.
In particular, every blow-up point is isolated, and by Proposition \ref{prop:isos}, it is isolated simple.
\end{prop}
\begin{proof}
Suppose not. Without loss of generality, we may assume that
\[\zeta_n := |x_{n,1}-x_{n,2}| = \min_{1 \le i \ne j \le k} |x_{n,i}-x_{n,j}| \to 0 \quad \text{as } n \to \infty.\]
Let
\[(w_n,z_n)(y) = \(\zeta_n^{\alpha_n} u_n(\zeta_n y + x_{n,1}),\,
\zeta_n^{\beta_n} v_n(\zeta_ny + x_{n,1})\) \quad \text{for } y \in (\Omega-x_{n,1})/\zeta_n.\]
It satisfies
\[\begin{cases}
-\Delta w_n = z_n^p &\text{in } (\Omega-x_{n,1})/\zeta_n,\\
-\Delta z_n = w_n^{q_n} &\text{in } (\Omega-x_{n,1})/\zeta_n,\\
w_n, z_n > 0 &\text{in } (\Omega-x_{n,1})/\zeta_n.
\end{cases}\]
Set $y_{n,2} = \zeta_n^{-1} (x_{n,2}-x_{n,1})$. By employing \eqref{eq:42} and arguing as in the proof of \cite[Proposition 4.2]{Li}, we see that
\[w_n(0),\, w_n(y_{n,2}),\, z_n(0),\, z_n(y_{n,2}) \to \infty \quad \text{as } n \to \infty.\]
Let $\by \in \R^N$ be such that $y_{n,2} \to \by \in \S^{N-1}$ as $n \to \infty$, passing to a subsequence.
Also, by \eqref{eq:41}, $0$ and $\by$ are isolated blow-up points of $\{(w_n,z_n)\}_{n \in \N}$. Thanks to Proposition \ref{prop:isos}, they are isolated simple.

Let $\wts$ be the set of blow-up points of $\{(w_n,z_n)\}_{n \in \N}$. Clearly, $0,\, \by \in \wts$ and $\inf\{|\by_1-\by_2|: \by_1,\, \by_2 \in \wts,\, \by_1 \ne \by_2\} \ge 1$.

At this moment, we divide the cases according to the value of $p$.

\medskip \noindent \textbf{Case 1: $p \in (\frac{N}{N-2}, \frac{N+2}{N-2})$.}
Denote $\tla_n = w_n^{1/\alpha_n}(0)$.
As in the proof of Proposition \ref{prop:uvdecay}, one can find constants $a_1,\, a_2,\, b_1,\, b_2 > 0$
and nonnegative harmonic functions $h_1^*$ and $h_2^*$ in $\R^N \setminus (\wts \setminus \{0,\by\})$ such that
\[\tla_n^{N-2-\alpha_n}w_n \to a_1|\cdot|^{2-N} + a_2|\cdot-\by|^{2-N}+h_1^* \quad \text{in } C^2_{\text{loc}}\(\R^N \setminus \wts\)\]
and
\begin{equation}\label{eq:dist1}
\tla_n^{N-2-\beta_n}z_n \to b_1|\cdot|^{2-N} + b_2|\cdot-\by|^{2-N}+h_2^* \quad \text{in } C^2_{\text{loc}}\(\R^N \setminus \wts\)
\end{equation}
as $n \to \infty$, up to a subsequence. In particular, the functions $a_2|\cdot-\by|^{2-N}+h_1^*$ and $b_2|\cdot-\by|^{2-N}+h_2^*$ are positive near $0$, which contradicts \eqref{eq:100}.

\medskip \noindent \textbf{Case 2: $p = \frac{N}{N-2}$.} The proof is similar to Case 1, so we omit it.

\medskip \noindent \textbf{Case 3: $p \in (\frac{2}{N-2}, \frac{N-1}{N-2})$.} The proof of Proposition \ref{prop:uvdecay} shows that \eqref{eq:dist1} holds and
\[\tla_n^{N-2-\alpha_n}w_n \to a_1|\cdot|^{2-(N-2)p} + a_2|\cdot-\by|^{2-(N-2)p}+h_1^* \quad \text{in } C^2_{\text{loc}}\(\R^N \setminus \wts\)\]
up to a subsequence, where $(a,b) = (a_1,b_1)$ or $(a_2,b_2)$ satisfies \eqref{eq:909}, and $h_1^* \in C^{1,\sigma}(\R^N \setminus (\wts \setminus \{0,\by\}))$ for some $\sigma \in (0,1)$.
By the maximum principle and the inequality
\[(1+x+y)^p \ge 1+x^p+y^p \quad \text{for } x, y \ge 0 \text{ and } p \ge 1,\]
$h_1^*$ is nonnegative. In particular, $a_2|\cdot-\by|^{2-(N-2)p}+h_1^* > 0$ near $0$, which contradicts \eqref{eq:101}.

\medskip \noindent \textbf{Case 4: $p \in [\frac{N-1}{N-2}, \frac{N}{N-2})$.}
Arguing as in Case 3, we deduce \eqref{eq:dist1} and
\[\tla_n^{N-2-\alpha_n}w_n \to a_1|\cdot|^{2-(N-2)p} + a_2|\cdot-\by|^{2-(N-2)p}+h_3^* \quad \text{in } C^2_{\text{loc}}\(\R^N \setminus \wts\)\]
up to a subsequence, where
\[0 \le h_3^* := a_1'|\cdot|^{N-(N-2)p}h_2^* + a_2'|\cdot-\by|^{N-(N-2)p}h_2^* + h_1^* \in C^{0,\sigma}\(\R^N \setminus \(\wts \setminus \{0,\by\}\)\)\]
for some $\sigma \in (0,1)$. In particular, $a_2|\cdot-\by|^{2-(N-2)p}+h_3^* > 0$ near $0$, which contradicts \eqref{eq:101}.
\end{proof}

\subsection{Proof of Theorem \ref{thm:main1}}\label{sec:Theorem01}
If $\{(u_{\e},v_{\e})\}_{\e \in (0,\e_0)}$ is bounded in $(L^{\infty}(\Omega))^2$, then elliptic regularity readily implies \eqref{eq:uvconv}.
In the following, we prove Theorem \ref{thm:main1} (1)--(3) after assuming that $(u_{\e},v_{\e})$ blows up at $k \in \N$ points $\xi_1, \ldots, \xi_k \in \Omega$ as $\e \to 0$.

\begin{proof}[Proof of Theorem \ref{thm:main1} (2)]
As already remarked, the convexity assumption on $\Omega$ leads us that each blow-up point $\xi_i$ is away from $\pa \Omega$.
Moreover, Proposition \ref{prop:uvdecay} implies that
\begin{equation}\label{eq:uvconv01}
(u_{\e}, v_{\e}) \to (0,0) \quad \text{in } C^1_{\textnormal{loc}}\(B\(\xi_i,\frac{\rho_2}{3}\) \setminus \{\xi_i\}\) \quad \text{as } \e \to 0
\end{equation}
for each index $i = 1, \ldots, k$ and a small number $\rho_2 > 0$. Then Harnack's inequality in \cite[Theorem 1.1]{CZ}, \eqref{eq:uvconv01}, 
and elliptic regularity 
yield
\begin{equation}\label{eq:uvconv02}
(u_{\e}, v_{\e}) \to (0,0) \quad \text{in } C^1_{\textnormal{loc}}\(\overline{\Omega} \setminus \cup_{i=1}^k B\(\xi_i,\frac{\rho_2}{6}\)\) \quad \text{as } \e \to 0.
\end{equation}
Combining \eqref{eq:uvconv01} and \eqref{eq:uvconv02}, we obtain \eqref{eq:uvconv0}.
\end{proof}

\begin{proof}[Proof of Theorem \ref{thm:main1} (3)]
The assertion \eqref{eq:uvconv1} follows from the proof of Lemma \ref{lemma:7}.

Let us prove \eqref{eq:uvconv11}. By \eqref{eq:901}, there is a constant $C > 0$ independent of $\e \in (0,\e_0)$ such that
\begin{equation}\label{eq:uvconv12}
v_{\e} \le C \lambda_{i\e}^{\beta_{\e}-(N-2)} \quad \text{in } B\(x_{i\e},\frac{\rho_2}{3}\) \setminus B\(x_{i\e},\frac{\rho_2}{6}\)
\end{equation}
for any $i = 1, \ldots, k$ and $\e \in (0,\e_0)$ small. On the other hand, by \eqref{eq:uvconv1} and $\alpha_{\e}q_{\e} - \beta_{\e} = 2$,
\begin{equation}\label{eq:uvconv13}
\begin{aligned}
v_{\e}(x) = \int_{\Omega} G(x,y) u_{\e}^{q_{\e}}(y)dy &\ge C \int_{B(x_{j\e},\frac{\rho_2}{3})} u_{\e}^{q_{\e}} \\
&\ge C \lambda_{j\e}^{\alpha_{\e}q_{\e}-N} \(\int_{\R^N} U_{1,0}^{q_0} + o(1)\) \ge C \lambda_{j\e}^{\beta_{\e}-(N-2)}
\end{aligned}
\end{equation}
for any $x \in B(x_{i\e},\frac{\rho_2}{3})$ and $1 \le j \ne i \le k$. From \eqref{eq:uvconv12} and \eqref{eq:uvconv13}, we discover \eqref{eq:uvconv11}.
\end{proof}

\begin{proof}[Proof of Theorem \ref{thm:main1} (1)]
By virtue of \eqref{eq:main0}, we may assume that
\[(u_{\e},v_{\e}) \to (u_0,v_0) \quad \begin{cases}
\text{weakly in } W^{2,{p+1 \over p}}(\Omega) \times W^{2,{q_0+1 \over q_0}}(\Omega)\\
\text{pointwise in } \Omega
\end{cases}
\text{as } \e \to 0,\]
up to a subsequence, and $(u_0,v_0)$ is a solution of \eqref{eq:LE} in which the positivity condition is replaced with $u,\, v \ge 0$ in $\Omega$.
On the other hand, thanks to \eqref{eq:uvconv01} or \eqref{eq:uvconv02}, $(u_0,v_0)$ vanishes at an interior point of $\Omega$. By the strong maximum principle, $(u_0,v_0) = (0,0)$ in $\Omega$.

Furthermore, by \eqref{eq:LE}, Theorem \ref{thm:main1} (3), Remark \ref{rmk:10}, and Fatou's lemma,
\begin{align*}
\liminf_{\e \to 0} \(\|u_{\e}\|_{W^{2,{p+1 \over p}}(\Omega)}^{p+1 \over p} + \|v_{\e}\|_{W^{2,{q_0+1 \over q_0}}(\Omega)}^{q_0+1 \over q_0}\)
&\ge \liminf_{\e \to 0} \left[\int_{B(x_{1\e},\lambda_{1\e}^{-1})} v_{\e}^{p+1} + \int_{B(x_{1\e},\lambda_{1\e}^{-1})} u_{\e}^{\frac{q_{\e}(q_0+1)}{q_0}} \right] \\
&\ge \int_{B(0,1)} V_{1,0}^{p+1} + \int_{B(0,1)} U_{1,0}^{q_0+1} > 0,
\end{align*}
which implies that there is no subsequence of $\{(u_{\e},v_{\e})\}_{\e \in (0,\e_0)}$ converging to $(0,0)$ strongly in $W^{2,(p+1)/p}(\Omega) \times W^{2,(q_0+1)/q_0}(\Omega)$.
\end{proof}

\section{Pointwise estimates for $(u_{\e},v_{\e})$}\label{sec:point}
Throughout this section, we assume that $N \ge 3$, $p \in (\frac{2}{N-2}, \frac{N+2}{N-2})$, and $p \ge 1$.

We fix $k \in \mathbb{N}$. Recalling $(\lambda_{i\e}, x_{i\e}) \in (0,\infty) \times \Omega$ in Theorem \ref{thm:main1} and the functions in \eqref{eq:bubble} and \eqref{eq:bubble21}--\eqref{eq:bubble22}, we write
\begin{equation}\label{eq:muUV}
\mu_{i\e} = \lambda_{i\e}^{-1}, \quad (U_{i\e}, V_{i\e}) = (U_{\mu_{i\e},x_{i\e}},V_{\mu_{i\e},x_{i\e}}),
\quad \text{and} \quad (\Psi_{i\e}^l,\Phi_{i\e}^l) = (\Psi_{\mu_{i\e},x_{i\e}}^l,\Phi_{\mu_{i\e},x_{i\e}}^l)
\end{equation}
for $i=1, \ldots, k$ and $l = 0, \ldots, N$.
Let also $(PU_{i\e}, PV_{i\e})$ be the unique solution of the system
\begin{equation}\label{eq:PUPV}
\begin{cases}
-\Delta PU_{i\e} = V_{i\e}^p &\text{in } \Omega,\\
-\Delta PV_{i\e} = U_{i\e}^{q_0} &\text{in } \Omega,\\
PU_{i\e} = PV_{i\e} = 0 &\text{on } \pa\Omega,
\end{cases}
\end{equation}
and $(P\Psi_{i\e}^l, P\Phi_{i\e}^l)$ be the unique solution of the system
\begin{equation}\label{eq:PsPh}
\begin{cases}
-\Delta P\Psi_{i\e}^l = pV_{i\e}^{p-1} \Phi_{i\e}^l & \text{ in } \Omega,\\
-\Delta P\Phi_{i\e}^l = q_0 U_{i\e}^{q_0-1} \Psi_{i\e}^l & \text{ in } \Omega,\\
P\Psi_{i\e}^l = P\Phi_{i\e}^l =0 & \text{ on } \pa\Omega.
\end{cases}
\end{equation}
Henceforth, we will often drop the subscript $\e$ for brevity, writing $\mu_i = \mu_{i\e}$, $U_i = U_{i\e}$, $\Psi_i^l = \Psi_{i\e}^l$, $PU_i = PU_{i\e}$, and so on.

By a standard comparison argument, we have
\begin{equation}\label{eq:PU_i}
PU_i(x) = \begin{cases}
U_i(x) - \frac{a_{N,p}}{\ga_N} \mu_i^{\frac{N}{p+1}} H(x,\xi_i) + o\(\mu_i^{\frac{N}{p+1}}\) &\text{if } p \in (\frac{N}{N-2}, \frac{N+2}{N-2}],\\
U_i(x) - \frac{a_{N,p}}{\ga_N} \mu_i^{\frac{N}{p+1}} \left[\whh(x,\xi_i) - \log\mu_i H(x,\xi_i)\right] + o\(\mu_i^{\frac{N}{p+1}}\) &\text{if } p = \frac{N}{N-2},\\
U_i(x) - \frac{a_{N,p}}{\ga_N} \mu_i^{\frac{Np}{q_0+1}} \whh(x,\xi_i) + o\(\mu_i^{\frac{Np}{q_0+1}}\) &\text{if } p \in (\frac{2}{N-2}, \frac{N}{N-2}),
\end{cases}
\end{equation}
and
\begin{equation}\label{eq:PV_i}
PV_i(x) = V_i(x) - \frac{b_{N,p}}{\ga_N} \mu_i^{\frac{N}{q_0+1}} H(x,\xi_i) + o\(\mu_i^{\frac{N}{q_0+1}}\)
\end{equation}
in $C^1(\overline{\Omega})$. Here, $H$ and $\whh$ are the functions satisfying \eqref{eq:H} and \eqref{eq:whh1}--\eqref{eq:whh2}, respectively. Similarly,
\begin{equation}\label{eq:PPs_i0}
P\Psi_i^0(x) = \begin{cases}
\Psi_i^0(x) + \frac{N}{p+1} \frac{a_{N,p}}{\ga_N} \mu_i^{\frac{N}{p+1}} H(x,\xi_i) + o\(\mu_i^{\frac{N}{p+1}}\) &\text{if } p \in (\frac{N}{N-2}, \frac{N+2}{N-2}],\\
\Psi_i^0(x) - \frac{N}{p+1} \frac{a_{N,p}}{\ga_N} \mu_i^{\frac{N}{p+1}} \log\mu_i H(x,\xi_i) + o\(\mu_i^{\frac{N}{p+1}} |\log\mu_i|\) &\text{if } p = \frac{N}{N-2},\\
\Psi_i^0(x) + \frac{Np}{q_0+1}\, \frac{a_{N,p}}{\ga_N} \mu_i^{\frac{Np}{q_0+1}} \whh(x,\xi_i) + o\(\mu_i^{\frac{Np}{q_0+1}}\) &\text{if } p \in (\frac{2}{N-2}, \frac{N}{N-2}),
\end{cases}
\end{equation}
\begin{equation}\label{eq:PPs_il}
P\Psi_i^l(x) = \begin{cases}
\Psi_i^l(x) + \frac{a_{N,p}}{\ga_N} \mu_i^{\frac{N}{p+1}} \pa_{\xi,l} H(x,\xi_i) + o\(\mu_i^{\frac{N}{p+1}}\) &\text{if } p \in (\frac{N}{N-2}, \frac{N+2}{N-2}],\\
\Psi_i^l(x) - \frac{a_{N,p}}{\ga_N} \mu_i^{\frac{N}{p+1}} \log\mu_i \pa_{\xi,l} H(x,\xi_i) + o\(\mu_i^{\frac{N}{p+1}} |\log\mu_i|\) &\text{if } p = \frac{N}{N-2},\\
\Psi_i^l(x) + \frac{a_{N,p}}{\ga_N} \mu_i^{\frac{Np}{q_0+1}} \pa_{\xi,l} \whh(x,\xi_i) + o\(\mu_i^{\frac{Np}{q_0+1}}\) &\text{if } p \in (\frac{2}{N-2}, \frac{N}{N-2}),
\end{cases}
\end{equation}
and
\begin{equation}\label{eq:PPh_il}
P\Phi_i^l(x) = \begin{cases}
\Phi_i^l(x) + \frac{N}{q_0+1} \frac{b_{N,p}}{\ga_N} \mu_i^{\frac{N}{q_0+1}} H(x,\xi_i) + o\(\mu_i^{\frac{N}{q_0+1}}\) &\text{if } l = 0,\\
\Phi_i^l(x) + \frac{b_{N,p}}{\ga_N} \mu_i^{\frac{N}{q_0+1}} \pa_{\xi,l} H(x,\xi_i) + o\(\mu_i^{\frac{N}{q_0+1}}\) &\text{if } l = 1, \ldots, N
\end{cases}
\end{equation}
in $C^1(\overline{\Omega})$. Here, $\pa_{\xi,l} H(x,\xi)$ and $\pa_{\xi,l} \whh(x,\xi)$ stand for the $l$-th components of $\nabla_{\xi} H(x,\xi)$ and $\nabla_{\xi} \whh(x,\xi)$, respectively.

Setting $\bsmu_{\e} = (\mu_{1\e}, \ldots, \mu_{k\e}) \in (0,\infty)^k$ and $\mbx_{\e} = (x_{1\e}, \ldots, x_{k\e}) \in \Omega^k$, we define $PU_{\bsmu_{\e},\mbx_{\e}}$ as the solution of the equation
\begin{equation}\label{eq:PU}
\begin{cases}
\displaystyle -\Delta PU_{\bsmu_{\e},\mbx_{\e}} = \(\sum_{i=1}^k PV_{i\e}\)^p &\text{in } \Omega,\\
PU_{\bsmu_{\e},\mbx_{\e}} = 0 &\text{on } \pa\Omega.
\end{cases}
\end{equation}
We will often drop the subscript $\e$, writing $\bsmu = \bsmu_{\e}$, $PU_{\mux} = PU_{\bsmu_{\e},\mbx_{\e}}$, and so on.

\begin{lemma}\label{lemma:PUdiff}
There exists a constant $C > 0$ independent of $\e \in (0,\e_0)$ such that
\begin{equation}\label{eq:PUdiff}
\left\|PU_{\bsmu_{\e},\mbx_{\e}} - \sum_{i=1}^k PU_{i\e}\right\|_{L^{\infty}(\Omega)} \le \begin{cases}
C\mu_{1\e}^{\frac{N}{p+1}} &\text{if } p \in (\frac{N}{N-2}, \frac{N+2}{N-2}],\\
C\mu_{1\e}^{\frac{N}{p+1}} |\log \mu_{1\e}| &\text{if } p = \frac{N}{N-2},\\
C\mu_{1\e}^{\frac{Np}{q_0+1}} &\text{if } p \in (\frac{2}{N-2}, \frac{N}{N-2})
\end{cases}
\end{equation}
for all $\e \in (0,\e_0)$.
\end{lemma}
\noindent For the proof of the above lemma, we need the following integral estimates.
\begin{lemma}\label{lemma:decay}
Let $x \in B(0,\lambda)$ for sufficiently large $\lambda > 0$. We have
\[\int_{B(0,\lambda)} \frac{1}{|x-y|^{N-2}}\frac{1}{(1+|y|)^a} dy \le
\begin{cases}
C(1+|x|)^{2-N} &\text{if } a > N,\\
C(1+|x|)^{2-a} &\text{if } 2 < a < N,\\
C\lambda^{2-a} &\text{if } 0 < a < 2.
\end{cases}\]
\end{lemma}
\begin{proof}
The proof is similar to that of \cite[Lemma B.2]{WY}. We omit the details.
\end{proof}
\begin{proof}[Proof of Lemma \ref{lemma:PUdiff}]
Fix a small number $\rho > 0$. By virtue of \eqref{eq:uvconv11}, \eqref{eq:bubble}, and \eqref{eq:PV_i},
\[|V_i(y)| + |PV_i(y)| = O\(\mu_1^{\frac{N}{q_0+1}}\) \quad \text{for any } y \in \Omega \setminus \cup_{i=1}^k B(x_i,\rho)\]
as $\e \to 0$. Hence
\begin{equation}\label{eq:PUdiff1}
PU_{\mux}(x)-\sum_{i=1}^k PU_i(x) = \sum_{i=1}^k \int_{B(x_i,\rho)} G(x,y) \left[\(\sum_{i=1}^k PV_i\)^p -\sum_{i=1}^k V_i^p\right](y)dy + O\(\mu_1^{\frac{Np}{q_0+1}}\)
\end{equation}
for $x \in \Omega$. Besides, by applying Lemma \ref{lemma:decay} and the fact that
\begin{equation}\label{eq:PUdiff2}
\frac{Np}{q_0+1} > \frac{N}{p+1} \quad \text{if and only if} \quad p > \frac{N}{N-2},
\end{equation}
one can show that the integral in the right-hand side of \eqref{eq:PUdiff1} is bounded by
\begin{equation}\label{eq:PUdiff3}
C \int_{B(x_i,\rho)} \frac{1}{|x-y|^{N-2}} \(\mu_i^{\frac{N}{q_0+1}}V_i^{p-1}(y) + \mu_i^{\frac{Np}{q_0+1}}\)dy \le \begin{cases}
C \mu_1^{\frac{N}{p+1}} &\text{if } p \in (\frac{N}{N-2}, \frac{N+2}{N-2}],\\
C \mu_1^{\frac{N}{p+1}} |\log \mu_1| &\text{if } p = \frac{N}{N-2},\\
C \mu_1^{\frac{Np}{q_0+1}} &\text{if } p \in (\frac{2}{N-2}, \frac{N}{N-2});
\end{cases}
\end{equation}
refer to the proof of \cite[Appendix A.4]{KP}. Therefore \eqref{eq:PUdiff} holds.
\end{proof}

In the rest of this section, we study the pointwise behavior of a family $\{(u_{\e},v_{\e})\}_{\e \in (0,\e_0)}$ of solutions of \eqref{eq:LE} with $q = q_{\e}$.
For simplicity, we assume that the results in Sections \ref{sec:dist} and \ref{sec:Theorem01} hold for the whole family $\{(u_{\e},v_{\e})\}_{\e \in (0,\e_0)}$, not only for its subsequence.

First, we examine $(u_{\e},v_{\e})$ near each blow-up point $\xi_i$ for $i = 1, \ldots, k$.
\begin{lemma}
Given $i = 1, \ldots, k$, let
\[(w_{i\e}(y), z_{i\e}(y)) = \(\lambda_{i\e}^{-\alpha_{\e}} u_{\e}(\lambda_{i\e}^{-1}y+x_{i\e}),
\lambda_{i\e}^{-\beta_{\e}} v_{\e}(\lambda_{i\e}^{-1}y+x_{i\e})\) \quad \text{for } y \in B(0,\lambda_{i\e}\rho_2)\]
where $\rho_2 > 0$ is a small number determined in the proof of Proposition \ref{prop:isos}.
Then there exists a constant $C > 0$ independent of $\e \in (0,\e_0)$ such that
\begin{multline}\label{eq:error11}
\|w_{i\e}-U_{1,0}\|_{L^{\infty}(B(0,\lambda_{i\e}\rho_2))} + \|z_{i\e}-V_{1,0}\|_{L^{\infty}(B(0,\lambda_{i\e}\rho_2))} \\
\le \begin{cases}
C \mu_{1\e}^{N-2} &\text{if } p \in (\frac{N}{N-2}, \frac{N+2}{N-2}),\\
C \mu_{1\e}^{N-2} |\log\mu_{1\e}| &\text{if } p = \frac{N}{N-2},\\
C \mu_{1\e}^{(N-2)p-2} &\text{if } p \in (\frac{2}{N-2}, \frac{N}{N-2})
\end{cases}
\end{multline}
for all $\e \in (0,\e_0)$.
\end{lemma}
\begin{proof}
Here, we deal with the case $p < \frac{N}{N-2}$ only. The remaining cases can be treated similarly.

\medskip
Assume that there exists a decreasing sequence $\{\e_n\}_{n \in \N} \subset (0,\e_0)$ such that $\e_n \to 0$ as $n \to \infty$ and
\[\|w_{in} - U_{1,0}\|_{L^{\infty}(B(0,\lambda_{in}\rho_2))} \ge \|z_{in} - V_{1,0}\|_{L^{\infty}(B(0,\lambda_{in}\rho_2))} \quad \text{for all } n \in \N\]
where $(w_{in},z_{in}) := (w_{i\e_n},z_{i\e_n})$ and $\lambda_{in} := \lambda_{i\e_n}$; if no such a sequence exists, we change the role of $w_{in} - U_{1,0}$ and $z_{in} - V_{1,0}$ in the following argument.
Let $\Lambda_n \ge 0$ and $y_n \in \overline{B(0,\lambda_{in}\rho_2)}$ satisfy
\[\Lambda_n = |w_{in}-U_{1,0}|(y_n) = \|w_{in} - U_{1,0}\|_{L^{\infty}(B(0,\lambda_{in}\rho_2))}.\]

Suppose that \eqref{eq:error11} does not hold. By \eqref{eq:uvconv11} and \eqref{eq:10}, we have
\begin{equation}\label{eq:error12}
\Lambda_n^{-1} \lambda_{in}^{2-(N-2)p} \to 0 \quad \text{and} \quad \Lambda_n^{-1} \e_n \to 0 \quad \text{as } n \to \infty.
\end{equation}
Also, we may assume that $y_n \in B(0,\frac{\lambda_{in}\rho_2}{2})$ for all $n \in \N$. If we set
\[(W_{in}, Z_{in}) = \(\Lambda_n^{-1} (w_{in} - U_{1,0}),\, \Lambda_n^{-1} (z_{in} - V_{1,0})\) \quad \text{in } B(0,\lambda_{in}\rho_2),\]
it satisfies
\[\begin{cases}
\displaystyle -\Delta W_{in} = \frac{z_{in}^p - V_{1,0}^p}{z_{in}-V_{1,0}}Z_{in} &\text{in } B(0,\lambda_{in}\rho_2),\\
\displaystyle -\Delta Z_{in} = \frac{w_{in}^{q_n} - U_{1,0}^{q_n}}{w_{in}-U_{1,0}}W_{in}
+ \Lambda_n^{-1}(U_{1,0}^{q_n}-U_{1,0}^{q_0}) &\text{in } B(0,\lambda_{in}\rho_2),\\
|W_{in}(y_n)| = 1 \text{ and } \nabla W_{in}(y_n) = 0
\end{cases}\]
where $q_n := q_{\e_n}$. By \eqref{eq:72}, and Propositions \ref{prop:dist} and \ref{prop:uvdecay}, it holds that
\begin{equation}\label{eq:error13}
w_{in}(y) \le C U_{1,0}(y) \quad \text{and} \quad z_{in}(y) \le C V_{1,0}(y) \quad \text{for } y \in B(0,\lambda_{in}\rho_2).
\end{equation}
Employing the representation formula, \eqref{eq:error13}, and Lemma \ref{lemma:decay}, we get
\begin{align*}
|Z_{in}(x)| &\le C \int_{B(0,\lambda_{in}\rho_2)} \frac{1}{|x-y|^{N-2}} \frac{1}{(1+|y|)^{((N-2)p-2)(q_n-1)}} dy \\
&\ + C \Lambda_n^{-1} \e_n \int_{B(0,\lambda_{in}\rho_2)} \frac{1}{|x-y|^{N-2}} \frac{\log(1+|y|)}{(1+|y|)^{((N-2)p-2)q_n}} dy + O\(\Lambda_n^{-1} \lambda_{in}^{2-N}\) \\
&\le \frac{C}{(1+|x|)^{\min\{N-(N-2)p+2,N-2\}}} 
+ \frac{C \Lambda_n^{-1} \e_n}{(1+|x|)^{N-2}} + O\(\Lambda_n^{-1} \lambda_{in}^{2-N}\)
\end{align*}
for all $x \in B(0,\frac{\lambda_{in}\rho_2}{2})$, as $n \to \infty$. From this, we compute
\begin{align*}
|W_{in}(x)| &\le C \int_{B(0,\lambda_{in}\rho_2)} \frac{1}{|x-y|^{N-2}} V_{1,0}^{p-1}(y) |Z_{in}(y)| dy + O\(\Lambda_n^{-1} \lambda_{in}^{2-(N-2)p}\) \\
&\le \frac{C}{(1+|x|)^{\min\{2,(N-2)p-2\}}} + \frac{C \Lambda_n^{-1} \e_n}{(1+|x|)^{(N-2)p-2}} + O\(\Lambda_n^{-1} \lambda_{in}^{2-(N-2)p}\) 
\end{align*}
for all $x \in B(0,\frac{\lambda_{in}\rho_2}{2})$. Therefore, together with \eqref{eq:error12},
we conclude that $|y_n|$ is bounded and there exists $y_0 \in \R^N$ such that $y_n \to y_0$ as $n \to \infty$, along a subsequence.
Furthermore, $(W_{in},Z_{in})$ converges to a pair $(W_{i0},Z_{i0})$, up to a subsequence, which satisfies
\[\begin{cases}
-\Delta W_{i0} = p V_{1,0}^{p-1} Z_{i0} \hspace{18pt} \text{in } \R^N,\\
-\Delta Z_{i0} = q_0 U_{1,0}^{q_0-1} W_{i0} \quad \text{in } \R^N,\\
|W_{i0}(x)| + |Z_{i0}(x)| \to 0 \quad \text{as } |x| \to \infty.
\end{cases}\]
Lemma \ref{lemma:FKP} and the conditions $W_{i0}(0) = \nabla W_{i0}(0) = 0$ show that $(W_{i0},Z_{i0}) = (0,0)$ in $\R^N$.
This is absurd because we also have that $|W_{i0}(y_0)| = 1$, so \eqref{eq:error11} is true.
\end{proof}

\noindent Combining \eqref{eq:uvconv11}, \eqref{eq:10}, \eqref{eq:PU_i}--\eqref{eq:PV_i}, \eqref{eq:PUdiff}, and \eqref{eq:error11}, we obtain
\begin{cor}
Set
\begin{equation}\label{eq:psph}
(\psi_{\e}, \phi_{\e}) = \(u_{\e}-PU_{\bsmu_{\e},\mbx_{\e}},\, v_{\e}-\sum_{i=1}^k PV_{i\e}\) \quad \text{in } \Omega.
\end{equation}
Then
\begin{equation}\label{eq:psph1}
\|\psi_{\e}\|_{L^{\infty}(B(x_{i\e},\rho_2))} \le \begin{cases}
C\mu_{1\e}^{\frac{N}{p+1}} &\text{if } p \in (\frac{N}{N-2}, \frac{N+2}{N-2}),\\
C\mu_{1\e}^{\frac{N}{p+1}} |\log \mu_{1\e}| &\text{if } p = \frac{N}{N-2},\\
C\mu_{1\e}^{\frac{Np}{q_0+1}} &\text{if } p \in (\frac{2}{N-2}, \frac{N}{N-2})
\end{cases}
\end{equation}
and
\begin{equation}\label{eq:psph2}
\|\phi_{\e}\|_{L^{\infty}(B(x_{i\e},\rho_2))} \le \begin{cases}
C\mu_{1\e}^{\frac{N}{q_0+1}} &\text{if } p \in (\frac{N}{N-2}, \frac{N+2}{N-2}),\\
C\mu_{1\e}^{\frac{N}{q_0+1}} |\log \mu_{1\e}| &\text{if } p = \frac{N}{N-2},\\
C\mu_{1\e}^{(N-2)p-2-\frac{N}{p+1}} &\text{if } p \in (\frac{2}{N-2}, \frac{N}{N-2})
\end{cases}
\end{equation}
for all $\e \in (0,\e_0)$.
\end{cor}

Next, we derive an upper bound of $(u_{\e},v_{\e})$ away from the points $\{\xi_1, \ldots, \xi_k\}$.
\begin{lemma}
Let $\mca_{\e,\rho_2} = \Omega \setminus \cup_{i=1}^k B(x_{i\e},\rho_2)$ where $\rho_2 > 0$ is determined in the proof of Proposition \ref{prop:isos}.
Then there exists a constant $C > 0$ independent of $\e \in (0,\e_0)$ such that
\begin{equation}\label{eq:outeru}
\|u_{\e}\|_{L^\infty (\mca_{\e,\rho_2})} \le \begin{cases}
C\mu_{1\e}^{N \over p+1} &\text{if } p \in (\frac{N}{N-2}, \frac{N+2}{N-2}),\\
C\mu_{1\e}^{N \over p+1} |\log \mu_{1\e}| &\text{if } p = \frac{N}{N-2},\\
C\mu_{1\e}^{Np \over q_0+1} &\text{if } p \in (\frac{2}{N-2}, \frac{N}{N-2})
\end{cases}
\end{equation}
and
\begin{equation}\label{eq:outerv}
\|v_{\e}\|_{L^\infty (\mca_{\e,\rho_2})} \le C \mu_{1\e}^{N \over q_0+1}
\end{equation}
for all $\e \in (0,\e_0)$.
\end{lemma}
\begin{proof}
Here, we deal with the case $p < \frac{N}{N-2}$ only. The remaining case can be treated similarly.

\medskip
By the representation formula and \eqref{eq:error13},
\begin{equation}\label{eq:outer1}
\begin{aligned}
u_{\e}(x) = \int_{\Omega} G(x,y) v_{\e}^p(y)dy &\le C \sum_{i=1}^k \int_{B(x_i,\frac{\rho_2}{2})} v_{\e}^p + C \|v_{\e}\|_{L^{\infty}(\mca_{\e,\rho_2/2})}^p \\
&\le C \mu_1^{Np \over q_0+1} + C \|v_{\e}\|_{L^{\infty}(\mca_{\e,\rho_2/2})}^p
\end{aligned}
\end{equation}
and
\begin{equation}\label{eq:outer2}
v_{\e}(x) \le C \sum_{i=1}^k \int_{B(x_i,\frac{\rho_2}{2})} u_{\e}^{q_{\e}} + C \|u_{\e}\|_{L^{\infty}(\mca_{\e,\rho_2/2})}^{q_{\e}}
\le C \mu_1^{N \over q_0+1} + C \|u_{\e}\|_{L^{\infty}(\mca_{\e,\rho_2/2})}^{q_{\e}}
\end{equation}
for any $x \in \mca_{\e,\rho_2}$. Applying \eqref{eq:outer1}, \eqref{eq:outer2}, \eqref{eq:uvconv0}, \eqref{eq:PUdiff2}, and the inequality $pq_{\e} > 1$, we obtain \eqref{eq:outeru}--\eqref{eq:outerv}.
\end{proof}

If $p \in (\frac{2}{N-2}, \frac{N}{N-2}]$, we have to improve the estimate of $\phi_{\e}$ which was given in \eqref{eq:psph2}.
\begin{lemma}
Let $a_0 = ((N-2)p-2)(q_0-1) > 4$. 
Let $\eta$ be $0$ for $a_0 < N$ and any small positive number for $a_0 = N$.

\medskip \noindent (1) Suppose that $p \in (\frac{2}{N-2}, \frac{N}{N-2})$. Then there exists a constant $C > 0$ independent of $\e \in (0,\e_0)$ such that
\begin{multline}\label{eq:error21}
|\phi_{\e}(x)| \le C \left[\mu_{1\e}^{\frac{Npq_0}{q_0+1}} + \mu_{1\e}^{\frac{N(p+2)}{q_0+1}} |\log\mu_{1\e}| \sum_{i=1}^k (\mu_{1\e}+|x-x_{i\e}|)^{2-N} \right. \\
\left. + \mu_{1\e}^{\frac{Npq_0}{q_0+1}-\eta} \sum_{i=1}^k (\mu_{1\e}+|x-x_{i\e}|)^{2-a_0+\eta} \mathbf{1}_{a_0 \le N}\right]
\end{multline}
for all $x \in \Omega$ and $\e \in (0,\e_0)$.

\medskip \noindent (2) Suppose that $p = \frac{N}{N-2}$. Then there exists a constant $C > 0$ independent of $\e \in (0,\e_0)$ such that
\begin{multline}\label{eq:error22}
|\phi_{\e}(x)| \le C\left[\mu_{1\e}^{\frac{Nq_0}{p+1}}|\log\mu_{1\e}|^{q_0} + C\mu_{1\e}^{\frac{N(p+2)}{q_0+1}} |\log\mu_{1\e}|^{q_0+1} \sum_{i=1}^k (\mu_{1\e}+|x-x_{i\e}|)^{2-N} \right. \\
\left. + C\mu_{1\e}^{\frac{Nq_0}{p+1}-\eta} |\log\mu_{1\e}|^{q_0} \sum_{i=1}^k (\mu_{1\e}+|x-x_{i\e}|)^{2-a_0+\eta} \mathbf{1}_{a_0 \le N}\right]
\end{multline}
for all $x \in \Omega$ and $\e \in (0,\e_0)$.
\end{lemma}
\begin{proof}
(1) By \eqref{eq:PUdiff}, \eqref{eq:PU_i}, \eqref{eq:psph1}, \eqref{eq:outeru}, \eqref{eq:10}, and the relation
\[\frac{N(p+2)}{q_0+1}-(N-a_0) 
= \frac{Npq_0}{q_0+1},\]
we have
\begin{equation}\label{eq:error23}
\begin{aligned}
|\phi_{\e}(x)| 
&= \sum_{i=1}^k \int_{B(x_i,\rho_2)} G(x,y)\left|\left\{U_i(y) + O\(\mu_1^{\frac{Np}{q_0+1}}\)\right\}^{q_{\e}} - U_i^{q_0}(y)\right| dy + O\(\mu_1^{\frac{Npq_0}{q_0+1}}\) \\
&\le C \sum_{i=1}^k \int_{B(x_i,\rho_2)} \frac{1}{|x-y|^{N-2}} \left[\e U_i^{q_{\e}}(y) |\log U_i(y)| + \mu_1^{\frac{Np}{q_0+1}} U_i^{q_0-1}(y)\right] dy + O\(\mu_1^{\frac{Npq_0}{q_0+1}}\) \\
&\le C\left[\mu_1^{\frac{N(p+2)}{q_0+1}} |\log\mu_1| \sum_{i=1}^k (\mu_1+|x-x_i|)^{2-N} \right. \\
&\hspace{125pt} \left. + \mu_1^{\frac{Npq_0}{q_0+1}-\eta} \sum_{i=1}^k (\mu_1+|x-x_i|)^{2-a_0+\eta} \mathbf{1}_{a_0 \le N}\right] + O\(\mu_1^{\frac{Npq_0}{q_0+1}}\)
\end{aligned}
\end{equation}
for all $x \in \Omega$ and $\e \in (0,\e_0)$. Here, the last inequality can be estimated as in \eqref{eq:PUdiff3}.

\medskip \noindent (2) Using
\[U_i(x) \le C\mu_i^{\frac{N}{p+1}} |\log\mu_i| (\mu_i+|x-x_i|)^{2-N}\]
for all $x \in \Omega$ and $i = 1, \ldots, k$, we can compute as in \eqref{eq:error23}.
\end{proof}

\section{Determining the blow-up locations}\label{sec:loc}
Letting $I_{p,q}$ be the energy functional \eqref{eq:energy} associated to \eqref{eq:LE}, we write $I_{\e} = I_{p,q_{\e}}$.
In this section, we will derive necessary conditions for the parameters $(\mux) = (\bsmu_{\e},\mbx_{\e}) \in (0,\infty)^k \times \Omega^k$ from the identities
\begin{equation}\label{eq:energy'}
I_{\e}'\(PU_{\mux} + \psi_{\e}, \sum_{i=1}^k PV_i + \phi_{\e}\) \(P\Psi_j^l, P\Phi_j^l\) = 0 \quad \text{for } j=1, \ldots, k \text{ and } l=0, \ldots, N
\end{equation}
where $(\psi_{\e},\phi_{\e})$ is the pair in \eqref{eq:psph}. Note that the left-hand side of \eqref{eq:energy'} equals
\begin{align}
&\ J_{1\e}(\mux) + J_{2\e}(\mux) + J_{3\e}(\mux) \nonumber \\
&:= \left[\int_{\Omega} \nabla PU_{\mux} \cdot \nabla P\Phi_j^l + \sum_{i=1}^k \int_{\Omega} \nabla PV_i \cdot \nabla P\Psi_j^l
- \int_{\Omega} \(\sum_{i=1}^k PV_i\)^p P\Phi_j^l - \int_{\Omega} (PU_{\mux})^{q_0} P\Psi_j^l\right] \nonumber \\
&\ + \left[\int_{\Omega} (PU_{\mux})^{q_0} P\Psi_j^l - \int_{\Omega} (PU_{\mux})^{q_{\e}} P\Psi_j^l\right] \label{eq:energy'2} \\
&\ + \left[I_{\e}'\(PU_{\mux} + \psi_{\e}, \sum_{i=1}^k PV_i + \phi_{\e}\) \(P\Psi_j^l, P\Phi_j^l\)
- I_{\e}'\(PU_{\mux}, \sum_{i=1}^k PV_i\) \(P\Psi_j^l,P\Phi_j^l\)\right]. \nonumber
\end{align}

A direct computation with \eqref{eq:PUPV}, \eqref{eq:PsPh}, \eqref{eq:PPs_i0}--\eqref{eq:PPh_il}, and \eqref{eq:PU} reveals that
\begin{equation}\label{eq:J1e}
J_{1\e}(\mux) = \int_{\Omega} \left[\sum_{i=1}^k U_i^{q_0} - (PU_{\mux})^{q_0}\right] P\Psi_j^l
\end{equation}
and
\begin{equation}\label{eq:J3e}
\begin{aligned}
J_{3\e}(\mux)
&= -\int_{\Omega} \left[(PU_{\mux} + \psi_{\e})^{q_{\e}} - (PU_{\mux})^{q_{\e}} - q_0 U_j^{q_0-1}\psi_{\e}\right] \Psi_j^l \\
&\ -\int_{\Omega} \left[\(\sum_{i=1}^k PV_i + \phi_{\e}\)^p - \(\sum_{i=1}^k PV_i\)^p - p V_j^{p-1} \phi_{\e} \right] \Phi_j^l\\
&\ + O\(\mu_1^{\frac{N}{q_0+1}}\) \int_{\Omega} \left|\(\sum_{i=1}^k PV_i + \phi_{\e}\)^p - \(\sum_{i=1}^k PV_i\)^p \right| \\
&\ + \begin{cases} O\(\mu_1^{\frac{N}{p+1}}\) \int_{\Omega} \left|(PU_{\mux} + \psi_{\e})^{q_{\e}} - (PU_{\mux})^{q_{\e}}\right| &\text{if } p \in (\frac{N}{N-2}, \frac{N+2}{N-2}), \\
O\(\mu_1^{\frac{N}{p+1}}|\log\mu_1|\) \int_{\Omega} \left|(PU_{\mux} + \psi_{\e})^{q_{\e}} - (PU_{\mux})^{q_{\e}}\right| &\text{if } p = \frac{N}{N-2}, \\
O\(\mu_1^{\frac{Np}{q_0+1}}\) \int_{\Omega} \left|(PU_{\mux} + \psi_{\e})^{q_{\e}} - (PU_{\mux})^{q_{\e}}\right| &\text{if } p \in (\frac{2}{N-2}, \frac{N}{N-2})
\end{cases}
\end{aligned}
\end{equation}
as $\e \to 0$. From now on, we study each of $J_{1\e}$, $J_{2\e}$, and $J_{3\e}$. There are four mutually exclusive cases according to the value of $p$.

\subsection{The case that $p \in (\frac{N}{N-2}, \frac{N+2}{N-2})$}\label{subsec:loc1}
We start with the following lemma.
\begin{lemma}\label{lemma:loc1}
If $N \ge 3$ and $p \in (\frac{N}{N-2}, \frac{N+2}{N-2})$, then there exists a constant $C > 0$ independent of $\e \in (0,\e_0)$ such that
\begin{align*}
\mu_i^{\frac{2N}{p+1}} \int_{B(x_i,\rho_2)} \(U_i^{q_{\e}-1} + U_i^{q_{\e}-2} \big|\Psi_i^0\big| \mathbf{1}_{q_{\e} \ge 2}\) &\le C\(\mu_1^{2(N-2)} + \mu_1^{\frac{N(q_0+1)}{p+1}}\), \\
\mu_i^{\frac{2N}{p+1}} \int_{B(x_i,\rho_2)} U_i^{q_{\e}-2}\big|\Psi_i^l\big| \mathbf{1}_{q_{\e} \ge 2} &\le C\(\mu_1^{2N-5} + \mu_1^{\frac{N(q_0+1)}{p+1}}\), \\
\mu_i^{\frac{2N}{q_0+1}} \int_{B(x_i,\rho_2)} \(V_i^{p-1} + V_i^{p-2} \big|\Phi_i^0\big| \mathbf{1}_{p \ge 2}\) &\le C\(\mu_1^{2(N-2)} + \mu_1^{(N-2)p-2}\), \\
\mu_i^{\frac{2N}{q_0+1}} \int_{B(x_i,\rho_2)} V_i^{p-2}\big|\Phi_i^l\big| \mathbf{1}_{p \ge 2} &\le C\(\mu_1^{2N-5} + \mu_1^{(N-2)p-2}\)
\end{align*}
for $i = 1, \ldots, k$ and $l = 1, \ldots, N$. In particular, all the integrals above are of order $o(\mu_1^{N-2})$ for $N \ge 4$.
\end{lemma}
\begin{proof}
Using \eqref{eq:bubble}, \eqref{eq:bubble21}--\eqref{eq:bubble22}, and \eqref{eq:V10est}--\eqref{eq:U10est1}, one can compute directly to derive the above results.
\end{proof}
In the rest of this subsection, we assume that $N \ge 4$.

\medskip \noindent \textsc{Estimate of $J_{3\e}$.}
According to \eqref{eq:J3e}, \eqref{eq:PUdiff}, and Lemma \ref{lemma:decay} and \ref{lemma:loc1}, we have
\begin{multline}\label{eq:r11}
|J_{3\e}(\mux)| \le \sum_{m=1}^k \int_{B(x_m,\rho_2)} \left|(PU_{\mux} + \psi_{\e})^{q_{\e}} - (PU_{\mux})^{q_{\e}} - q_0 U_j^{q_0-1}\psi_{\e}\right| \big|\Psi_j^l\big| \\
+ \sum_{m=1}^k \int_{B(x_m,\rho_2)} \left|\(\sum_{i=1}^k PV_i + \phi_{\e}\)^p - \(\sum_{i=1}^k PV_i\)^p - p V_j^{p-1} \phi_{\e} \right| \big|\Phi_j^l\big| + o\(\mu_1^{N-2}\).
\end{multline}
By \eqref{eq:psph1}, it holds that
\begin{multline}\label{eq:r12}
\left|(PU_{\mux} + \psi_{\e})^{q_{\e}} - (PU_{\mux})^{q_{\e}} - q_0 U_j^{q_0-1}\psi_{\e}\right| \big|\Psi_j^l\big| \\
\le C\(\mu_1^{\frac{2N}{p+1}} U_m^{q_{\e}-1} + \mu_1^{\frac{N(q_0+1)}{p+1}}\) \quad \text{in } B(x_m,\rho_2)
\end{multline}
for $m \ne j$, and
\begin{equation}\label{eq:r13}
\begin{aligned}
&\ \left|(PU_{\mux} + \psi_{\e})^{q_{\e}} - (PU_{\mux})^{q_{\e}} - q_0 U_j^{q_0-1}\psi_{\e}\right| \\
&\le C\left[(PU_{\mux})^{q_{\e}-2}\psi_{\e}^2 \mathbf{1}_{q_{\e} \ge 2} + |\psi_{\e}|^{q_{\e}} + \left|q_{\e} (PU_{\mux})^{q_{\e}-1} - q_0U_j^{q_0-1}\right| |\psi_{\e}|\right] \\
&\le C\left[\mu_1^{\frac{Nq_0}{p+1}} + \mu_1^{\frac{2N}{p+1}} U_j^{q_{\e}-2} \mathbf{1}_{q_{\e} \ge 2} + \e \mu_1^{N \over p+1} U_j^{q_0-1} \(1 + U_j^{O(\e)}|\log U_j|\)\right]
\end{aligned}
\quad \text{in } B(x_j,\rho_2).
\end{equation}
Thus we infer from Lemma \ref{lemma:loc1}, \eqref{eq:r12}--\eqref{eq:r13}, and \eqref{eq:10} that
\begin{equation}\label{eq:r130}
\sum_{m=1}^k \int_{B(x_m,\rho_2)} \left|(PU_{\mux} + \psi_{\e})^{q_{\e}}
- (PU_{\mux})^{q_{\e}} - q_0 U_j^{q_0-1}\psi_{\e}\right| \big|\Psi_j^l\big| = o\(\mu_1^{N-2}\).
\end{equation}
Similarly, we find
\begin{equation}\label{eq:r133}
\sum_{m=1}^k \int_{B(x_m,\rho_2)} \left|\(\sum_{i=1}^k PV_i + \phi_{\e}\)^p - \(\sum_{i=1}^k PV_i\)^p - p V_j^{p-1} \phi_{\e} \right| \big|\Phi_j^l\big| = o\(\mu_1^{N-2}\).
\end{equation}
As a result,
\begin{equation}\label{eq:r131}
J_{3\e}(\mux) = o\(\mu_1^{N-2}\).
\end{equation}

\medskip \noindent \textsc{Estimate of $J_{1\e}$.}
We observe from \eqref{eq:J1e}, \eqref{eq:PU_i}, \eqref{eq:PPs_i0}--\eqref{eq:PPs_il}, \eqref{eq:PUdiff}, and Lemma \ref{lemma:decay} that
\begin{equation}\label{eq:r132}
J_{1\e}(\mux) = - \sum_{m=1}^k \int_{B(x_m,\rho_2)} \left[(PU_{\mux})^{q_0} - \sum_{i=1}^k U_i^{q_0}\right] P\Psi_j^l + o\(\mu_1^{N-2}\).
\end{equation}
Define
\begin{equation}\label{eq:hps}
\hps_{\mux} = PU_{\mux} - \sum_{i=1}^k PU_i \quad \text{in } \Omega.
\end{equation}
Then, \eqref{eq:PUdiff} yields that $\hps_{\mux} = O(\mu_1^{\frac{N}{p+1}})$, and \eqref{eq:PU_i} leads to
\begin{equation}\label{eq:r134}
PU_{\mux} = U_m + \sum_{i=1,\, i \ne m}^k U_i - \frac{a_{N,p}}{\ga_N} \sum_{i=1}^k \mu_i^{\frac{N}{p+1}}
H(\cdot,\xi_i) + \hps_{\mux} + \mcr_1 \quad \text{in } B(x_m,\rho_2)
\end{equation}
for each $m = 1, \ldots, k$. Here, $\mcr_1$ is a remainder term satisfying
\begin{equation}\label{eq:r135}
|\mcr_1| + |\nabla \mcr_1| = o\(\mu_1^{\frac{N}{p+1}}\).
\end{equation}
From \eqref{eq:r134}, we obtain
\begin{multline}\label{eq:r14}
(PU_{\mux})^{q_0} - \sum_{i=1}^k U_i^{q_0} = q_0 U_m^{q_0-1} \left[\sum_{i=1,\, i \ne m}^k U_i - \frac{a_{N,p}}{\ga_N} \sum_{i=1}^k \mu_i^{\frac{N}{p+1}} H(\cdot,\xi_i) + \hps_{\mux} + \mcr_1\right] \\
+ O\(\mu_1^{\frac{2N}{p+1}}U_m^{q_0-2} \mathbf{1}_{q_0 \ge 2}\) + O\(\mu_1^{\frac{Nq_0}{p+1}}\) \quad \text{in } B(x_m,\rho_2).
\end{multline}
Also, using \eqref{eq:U10est1} and \eqref{eq:r135}, we see
\[U_i(x) = U_i(x_m) + \nabla U_i(x_m) \cdot (x-x_m) + O\(\mu_1^{\frac{N}{p+1}} |x-x_m|^2\) \quad \text{for } x \in B(x_m,\rho_2),\]
\[\nabla U_i(x) = -\mu_i^{\frac{N}{p+1}} (N-2)a_{N,p} \frac{x-x_i}{|x-x_i|^N} + o\(\mu_1^{\frac{N}{p+1}}\) \quad \text{for } x \in B(x_m,\rho_2) \text{ and } i \ne m,\]
and
\[\mcr_1(x) = \mcr_1(x_m) + o\(\mu_1^{\frac{N}{p+1}}|x-x_m|\) \quad \text{for } x \in B(x_m,\rho_2).\]
For the moment, we assume that $l = 1, \ldots, N$. Applying the above estimates, \eqref{eq:r14}, Lemma \ref{lemma:loc1}, and the oddness of $\Psi_{1,0}^l$ in the $l$-th variable, we deduce
\begin{align}
&\ \sum_{m=1}^k \int_{B(x_m,\rho_2)} \left[(PU_{\mux})^{q_0} - \sum_{i=1}^k U_i^{q_0}\right] P\Psi_j^l \nonumber \\
&= q_0 \sum_{m=1}^k \int_{B(x_m,\rho_2)} U_m^{q_0-1}(x) \left[\sum_{i=1,\, i \ne m}^k U_i(x)
- \frac{a_{N,p}}{\ga_N} \sum_{i=1}^k \mu_i^{\frac{N}{p+1}} H(x,\xi_i) + \hps_{\mux}(x) + \mcr_1(x)\right] \nonumber \\
&\hspace{122pt} \times \left[\Psi_j^l(x) + \frac{a_{N,p}}{\ga_N} \mu_j^{\frac{N}{p+1}} \pa_{\xi,l} H(x,\xi_j)
+ o\(\mu_j^{\frac{N}{p+1}}\)\right] dx + o\(\mu_1^{N-2}\) \nonumber \\
&= q_0 \int_{B(x_j,\rho_2)} U_j^{q_0-1}(x) \left[\sum_{i=1,\, i \ne j}^k U_i(x) - \frac{a_{N,p}}{\ga_N} \sum_{i=1}^k \mu_i^{\frac{N}{p+1}} H(x,\xi_i)
+ \hps_{\mux}(x) + \mcr_1(x)\right] \Psi_j^l(x)dx \nonumber \\
&\ + o\(\mu_1^{N-2}\) \label{eq:r15} \\
&= - q_0(N-2) a_{N,p} \sum_{i=1,\, i \ne j}^k \mu_i^{\frac{N}{p+1}} \int_{B(\xi_j,\rho_2)} U_j^{q_0-1}(x) \frac{(\xi_j-\xi_i) \cdot (x-\xi_j)}{|\xi_j-\xi_i|^N} \Psi_j^l(x) dx \nonumber \\
&\ - q_0 \frac{a_{N,p}}{\ga_N} \sum_{i=1}^k \mu_i^{\frac{N}{p+1}} \int_{B(\xi_j,\rho_2)} U_j^{q_0-1}(x) \nabla_x H(\xi_j,\xi_i) \cdot (x-\xi_j) \Psi_j^l(x) dx \nonumber \\
&\ + q_0 \int_{B(x_j,\rho_2)} U_j^{q_0-1}\hps_{\mux} \Psi_j^l + o\(\mu_1^{\frac{N}{p+1}} \int_{B(x_j,\rho_2)} U_j^{q_0-1}(x) |x-x_j| \big|\Psi_j^l(x)\big| dx\) \nonumber \\
&\ + O\(\mu_1^{\frac{N}{p+1}} \int_{B(x_j,\rho_2)} U_j^{q_0-1}(x) |x-x_j|^2 \big|\Psi_j^l(x)\big| dx\) + o\(\mu_1^{N-2}\) \nonumber
\end{align}
where $\nabla_x H(\xi_j,\xi_i) := \nabla_x H(x,\xi_i)|_{x=\xi_j}$.
By reminding \eqref{eq:A12} and using integration by parts, it is plain to check the identities
\begin{equation}\label{eq:A1}
A_1 = - q_0 \int_{\R^N} U_{1,0}^{q_0-1} \Psi_{1,0}^1 x_1 = \cdots = - q_0 \int_{\R^N} U_{1,0}^{q_0-1} \Psi_{1,0}^N x_N = - \frac{q_0(q_0+1)}{N} \int_{\R^N} U_{1,0}^{q_0-1} \Psi_{1,0}^0
\end{equation}
and
\begin{equation}\label{eq:A2}
A_2 = -p \int_{\R^N} V_{1,0}^{p-1} \Phi_{1,0}^1 x_1 = \cdots
= -p \int_{\R^N} V_{1,0}^{p-1} \Phi_{1,0}^N x_N = -\frac{p(p+1)}{N} \int_{\R^N} V_{1,0}^{p-1} \Phi_{1,0}^0.
\end{equation}
Using \eqref{eq:A1}, we evaluate each term of the rightmost side of \eqref{eq:r15}. Then we get
\begin{multline}\label{eq:r16}
\sum_{m=1}^k \int_{B(x_m,\rho_2)} \left[(PU_{\mux})^{q_0} - \sum_{i=1}^k U_i^{q_0}\right] P\Psi_j^l
= \mu_j^{N-2} \frac{a_{N,p}}{\ga_N} A_1 \pa_l H(\xi_j,\xi_j) \\
- \frac{a_{N,p}}{\ga_N} A_1 \sum_{i=1,\, i \ne j}^k \mu_i^{\frac{N}{p+1}} \mu_j^{\frac{N}{q_0+1}} \pa_l G(\xi_j,\xi_i)
+ q_0 \int_{B(x_j,\rho_2)} U_j^{q_0-1}\hps_{\mux} \Psi_j^l + o\(\mu_1^{N-2}\)
\end{multline}
where $\pa_l H(\xi_j,\xi_j) := \pa_l H(x,\xi_j)|_{x=\xi_j}$ and $\pa_l G(\xi_j,\xi_i) := \pa_l G(x,\xi_i)|_{x=\xi_j}$, respectively.
On the other hand, by the representation formula for $\hps_{\mux}$, Fubini's theorem and the bound
\[\left|U_j^{q_0-1}\Psi_j^l\right| = O\(\mu_1^{\frac{Nq_0}{p+1}}\) = o\(\mu_1^{\frac{N}{q_0+1}}\) \quad \text{in } C^1\(\overline{\Omega} \setminus B(x_j,\rho_2)\),\]
we have
\begin{equation}\label{eq:r17}
\begin{aligned}
&\ q_0 \int_{B(x_j,\rho_2)} U_j^{q_0-1}\hps_{\mux} \Psi_j^l \\
&= q_0 \int_{\Omega} \left[\(\sum_{i=1}^k PV_i\)^p - \sum_{i=1}^k V_i^p \right](y) \int_{B(x_j,\rho_2)} G(y,x) U_j^{q_0-1}(x) \Psi_j^l(x) dx dy\\
&= \int_{\Omega} \left[\(\sum_{i=1}^k PV_i\)^p -\sum_{i=1}^k V_i^p \right] \(P\Phi_j^l + \mcr_2\)
\end{aligned}
\end{equation}
where $\mcr_2$ is a remainder term such that $|\mcr_2| + |\nabla \mcr_2| = o(\mu_1^{\frac{N}{q_0+1}})$.
Applying \eqref{eq:A2} and arguing as in \eqref{eq:r15} and \eqref{eq:r16}, we see that the rightmost side of \eqref{eq:r17} is equal to
\begin{equation}\label{eq:r18}
\begin{aligned}
&\ p \int_{B(x_j,\rho_2)} V_j^{p-1}(x) \left[\sum_{i=1,\, i \ne j}^k V_i(x) - \frac{b_{N,p}}{\ga_N}
\sum_{i=1}^k \mu_i^{\frac{N}{q_0+1}} H(x,\xi_i) + \mcr_2'(x)\right] \Phi_j^l(x) dx + o\(\mu_1^{N-2}\)\\
&= \mu_j^{N-2} \frac{b_{N,p}}{\ga_N} A_2 \pa_l H(\xi_j,\xi_j)
- \frac{b_{N,p}}{\ga_N} A_2 \sum_{i=1,\, i \ne j}^k \mu_i^{\frac{N}{q_0+1}} \mu_j^{\frac{N}{p+1}} \pa_l G(\xi_j,\xi_i) + o\(\mu_1^{N-2}\)
\end{aligned}
\end{equation}
where $\mcr_2'$ is a term satisfying $|\mcr_2| + |\nabla \mcr_2| = o(\mu_1^{\frac{N}{q_0+1}})$.
Consequently, by combining \eqref{eq:r16} and \eqref{eq:r17}--\eqref{eq:r18}, and employing \eqref{eq:ab}, we obtain
\begin{multline}\label{eq:r19}
\sum_{m=1}^k \int_{B(x_m,\rho_2)} \left[(PU_{\mux})^{q_0} - \sum_{i=1}^k U_i^{q_0}\right] P\Psi_j^l
= 2\mu_j^{N-2} \frac{a_{N,p}}{\ga_N} A_1 \pa_l H(\xi_j,\xi_j) \\
- \frac{a_{N,p}}{\ga_N} A_1 \sum_{i=1,\, i \ne j}^k \(\mu_i^{\frac{N}{p+1}} \mu_j^{\frac{N}{q_0+1}}
+ \mu_i^{\frac{N}{q_0+1}} \mu_j^{\frac{N}{p+1}}\) \pa_l G(\xi_j,\xi_i) + o\(\mu_1^{N-2}\).
\end{multline}

Analogously, if $l = 0$, then using the evenness of $\Psi_m^0$, \eqref{eq:A1}--\eqref{eq:A2}, and \eqref{eq:cr-hy}, we derive
\begin{multline}\label{eq:r191}
\sum_{m=1}^k \int_{B(x_m,\rho_2)} \left[(PU_{\mux})^{q_0} - \sum_{i=1}^k U_i^{q_0}\right] P\Psi_j^0
= \mu_j^{N-2} \frac{a_{N,p}}{\ga_N} (N-2)A_1 H(\xi_j, \xi_j) \\
- \frac{a_{N,p}}{\ga_N} NA_1 \sum_{i=1,\, i \ne j}^k \(\frac{\mu_j^{\frac{N}{q_0+1}} \mu_i^{\frac{N}{p+1}}}{q_0+1}
+ \frac{\mu_j^{\frac{N}{p+1}} \mu_i^{\frac{N}{q_0+1}}}{p+1}\) G(\xi_j, \xi_i) + o\(\mu_1^{N-2}\).
\end{multline}

\medskip \noindent \textsc{Estimate of $J_{2\e}$.}
Assume that $l = 1, \ldots, N$. Using Talyor's theorem, the relation that $|q_0-q_{\e}| = O(\e)$, \eqref{eq:10}, $N \ge 4$, and the oddness of $\Psi_{1,0}^l$ in the $l$-th variable, we get
\begin{align}
&\int_{\Omega} \left[(PU_{\mux})^{q_0} - (PU_{\mux})^{q_{\e}}\right] P\Psi_j^l \nonumber \\
&= (q_0-q_{\e}) \int_{\Omega} (PU_{\mux})^{q_0} \log (PU_{\mux}) P\Psi_j^l + O\(\e^2 \int_{\Omega} (PU_{\mux})^{q_0+O(\e)} P\Psi_j^l\) \nonumber \\
&= (q_0-q_{\e}) \sum_{m=1}^k \int_{B(x_m,\rho_2)} (PU_{\mux})^{q_0} \log (PU_{\mux}) P\Psi_j^l + o\(\e \mu_1^{N-2}\) + O\(\e^2\) \label{eq:r192} \\
&= (q_0-q_{\e}) q_0 \int_{B(x_j,\rho_2)} U_j^{q_0}\log U_j \left[\Psi_j^l(x)
+ \frac{a_{N,p}}{\ga_N} \mu_j^{\frac{N}{p+1}} \pa_{\xi,l} H(x,\xi_j) + o\(\mu_j^{\frac{N}{p+1}}\)\right] dx + o\(\mu_1^{N-2}\) \nonumber \\
&= o\(\mu_1^{N-2}\). \nonumber
\end{align}
We note that
\[\int_{\R^N} U_{\mu,0}^{q_0+1} = \int_{\R^N} U_{1,0}^{q_0+1} \quad \text{and} \quad
\int_{\R^N} U_{\mu,0}^{q_0+1} \log U_{\mu,0} = -\frac{N}{q_0+1} \log \mu \int_{\R^N} U_{1,0}^{q_0+1} + \int_{\R^N} U_{1,0}^{q_0+1} \log U_{1,0}\]
for all $\mu > 0$. Differentiating the both equalities with respect to $\mu$ and putting $\mu = 1$ on the results, we obtain
\begin{equation}\label{eq:r194}
\int_{\R^N} U_{1,0}^p \Psi_{1,0}^0 = 0 \quad \text{and} \quad
A_3 := \int_{\R^N} U_{1,0}^{q_0} \log U_{1,0} \Psi_{1,0}^0 = \frac{N}{(q_0+1)^2} \int_{\R^N} U_{1,0}^{q_0+1} dx > 0.
\end{equation}
Using \eqref{eq:r194} and computing as in \eqref{eq:r192}, we have
\begin{align}
&\ \int_{\Omega} \left[(PU_{\mux})^{q_0} - (PU_{\mux})^{q_{\e}}\right] P\Psi_j^0 \nonumber \\
&= (q_0-q_{\e}) q_0 \int_{B(x_j,\rho_2)} U_j^{q_0}\log U_j \left[\Psi_j^0(x)
+ \frac{N}{p+1} \frac{a_{N,p}}{\ga_N} \mu_j^{\frac{N}{p+1}} H(x,\xi_j) + o\(\mu_j^{\frac{N}{p+1}}\)\right] dx + o\(\mu_1^{N-2}\) \nonumber \\
&= (q_0-q_{\e}) (A_3+o(1)) + o\(\mu_1^{N-2}\) \label{eq:r193}.
\end{align}

\medskip \noindent \textsc{Conclusion.}
Recalling \eqref{eq:uvconv11}, we write $d_i = d_{i\e} = \mu_{i\e} \mu_{1\e}^{-1} \in (0,\infty)$ for $i = 1, \ldots, k$ and $\e \in (0,\e_0)$.
From \eqref{eq:energy'}--\eqref{eq:energy'2}, \eqref{eq:r131}, \eqref{eq:r132}, \eqref{eq:r19}--\eqref{eq:r191}, \eqref{eq:r192}, and \eqref{eq:r193}, we see
\begin{equation}\label{eq:r195}
2d_j^{N-2} (\nabla_{\xi} H)(\xi_j,\xi_j) - \sum_{i=1,\, i \ne j}^k \(d_i^{\frac{N}{p+1}} d_j^{\frac{N}{q_0+1}}
+ d_i^{\frac{N}{q_0+1}} d_j^{\frac{N}{p+1}}\) (\nabla_{\xi} G)(\xi_i,\xi_j) = o(1)
\end{equation}
where $(\nabla_{\xi} H)(\xi_j,\xi_j) := \nabla_{\xi} H(\xi_j,\xi)|_{\xi=\xi_j}$ and $(\nabla_{\xi} G)(\xi_i,\xi_j) := \nabla_{\xi} G(\xi_i,\xi)|_{\xi=\xi_j}$, respectively, and
\begin{multline}\label{eq:r196}
(N-2) d_j^{N-2} H(\xi_j, \xi_j) - N \sum_{i=1,\, i \ne j}^k \(\frac{d_i^{\frac{N}{p+1}} d_j^{\frac{N}{q_0+1}}}{q_0+1}
+ \frac{d_i^{\frac{N}{q_0+1}} d_j^{\frac{N}{p+1}}}{p+1}\) G(\xi_i,\xi_j) \\
= \frac{(A_3+o(1))\ga_N}{a_{N,p}A_1} \frac{q_0-q_{\e}}{\mu_1^{N-2}} + o(1)
\end{multline}
for $j = 1, \ldots, k$.

\subsection{The case that $p = \frac{N}{N-2}$}
It turns out that analysis of this threshold situation resembles to that of the previous case, albeit we also have to exploit a refined estimate of $\phi_{\e}$ given in \eqref{eq:error22}.
In this subsection, we will briefly explain how the computations need to be changed.

\medskip \noindent \textsc{Estimate of $J_{3\e}$.}
A slight alteration of the argument in the previous subsection yield \eqref{eq:r11} and \eqref{eq:r130}.
However, an analogous argument no more yields \eqref{eq:r133}.
Instead, one can deduce it by applying \eqref{eq:error22} as a substitute and conducting tedious but straightforward computations. In conclusion,
\begin{equation}\label{eq:r221}
J_{3\e}(\mux) = o\(\mu_1^{N-2}\)
\end{equation}
for $N \ge 4$.

\medskip \noindent \textsc{Estimate of $J_{1\e}$.}
We argue as in the previous subsection, suitably caring for the modified terms in the expansion of $PU_i$.
Since the integral $A_2$ in \eqref{eq:A12} is divergent when $p = \frac{N}{N-2}$, we instead work with the quantity
\[A_2[R] = \int_{B(0,R)} V_{1,0}^p = -p\int_{B(0,R)} V_{1,0}^{p-1}\Phi_{1,0}^l x_l = -\frac{p(p+1)}{N} \int_{B(0,R)} V_{1,0}^{p-1}\Phi_{1,0}^0\]
where $R > 0$ and $l = 1, \ldots, N$. By \eqref{eq:U10est2},
\[A_2[R] = b_{N,p}^p \left|\S^{N-1}\right| \log R\, (1+o_R(1)) \quad \text{as } R \to \infty\]
where $o_R(1) \to 0$ as $R \to \infty$. In particular, we have a formula analogous to \eqref{eq:r18}:
\begin{equation}\label{eq:r222}
\begin{aligned}
&\ p \int_{B(x_j,\rho_2)} V_j^{p-1}(x) \left[\sum_{i=1,\, i \ne j}^k V_i(x) - \frac{b_{N,p}}{\ga_N} \sum_{i=1}^k \mu_i^{\frac{N}{q_0+1}} H(x,\xi_i) + \mcr_2'(x)\right] \Phi_j^l(x) dx + O\(\mu_1^{N-2}\) \\
&= \mu_j^{N-2} \frac{b_{N,p}}{\ga_N} A_2\left[\rho_2 \mu_j^{-1}\right] \pa_l H(\xi_j,\xi_j)
- \frac{b_{N,p}}{\ga_N} A_2\left[\rho_2 \mu_j^{-1}\right] \sum_{i=1,\, i \ne j}^k \mu_i^{\frac{N}{q_0+1}} \mu_j^{\frac{N}{p+1}} \pa_l G(\xi_j,\xi_i) \\
&\ + o\(\mu_1^{N-2}|\log\mu_1|\) \\
&= \mu_j^{N-2} |\log\mu_j| \frac{b_{N,p}^{p+1}}{\ga_N} \left|\S^{N-1}\right| \pa_l H(\xi_j,\xi_j)
- \frac{b_{N,p}^{p+1}}{\ga_N} \left|\S^{N-1}\right| \sum_{i=1,\, i \ne j}^k \mu_i^{\frac{N}{q_0+1}} \mu_j^{\frac{N}{p+1}} |\log\mu_j| \pa_l G(\xi_j,\xi_i) \\
&\ + o\(\mu_1^{N-2}|\log\mu_1|\).
\end{aligned}
\end{equation}
From \eqref{eq:r222} and \eqref{eq:ab1}, we deduce \eqref{eq:r132}, \begin{multline}\label{eq:r28}
\sum_{m=1}^k \int_{B(x_m,\rho_2)} \left[(PU_{\mux})^{q_0} - \sum_{i=1}^k U_i^{q_0}\right] P\Psi_j^l
= 2\mu_j^{N-2} |\log\mu_1| \frac{a_{N,p}}{\ga_N} A_1 \pa_l H(\xi_j,\xi_j) \\
- \frac{a_{N,p}}{\ga_N} A_1 \sum_{i=1,\, i \ne j}^k \(\mu_i^{\frac{N}{p+1}} \mu_j^{\frac{N}{q_0+1}}
+ \mu_i^{\frac{N}{q_0+1}} \mu_j^{\frac{N}{p+1}}\)|\log\mu_1| \pa_l G(\xi_j,\xi_i) + o\(\mu_1^{N-2}|\log\mu_1|\)
\end{multline}
for $l = 1, \ldots, N$, and
\begin{multline}\label{eq:r281}
\sum_{m=1}^k \int_{B(x_m,\rho_2)} \left[(PU_{\mux})^{q_0} - \sum_{i=1}^k U_i^{q_0}\right] P\Psi_j^0
= \mu_j^{N-2} |\log\mu_1| \frac{a_{N,p}}{\ga_N} (N-2)A_1 H(\xi_j, \xi_j) \\
- \frac{a_{N,p}}{\ga_N} NA_1 \sum_{i=1,\, i \ne j}^k \(\frac{\mu_j^{\frac{N}{q_0+1}} \mu_i^{\frac{N}{p+1}}}{q_0+1}
+ \frac{\mu_j^{\frac{N}{p+1}} \mu_i^{\frac{N}{q_0+1}}}{p+1}\)|\log\mu_1| G(\xi_j, \xi_i) + o\(\mu_1^{N-2}|\log\mu_1|\).
\end{multline}

\medskip \noindent \textsc{Estimate of $J_{2\e}$.}
Arguing as in \eqref{eq:r192}, we get
\begin{equation}\label{eq:r282}
\int_{\Omega} \left[(PU_{\mux})^{q_0} - (PU_{\mux})^{q_{\e}}\right] P\Psi_j^l = o\(\mu_1^{N-2}\)
\end{equation}
for $l = 1, \ldots, N$. Moreover, as in \eqref{eq:r193}, we have
\begin{align}
&\ \int_{\Omega} \left[(PU_{\mux})^{q_0} - (PU_{\mux})^{q_{\e}}\right] P\Psi_j^0 \nonumber \\
&= (q_0-q_{\e}) q_0 \int_{B(x_j,\rho_2)} U_j^{q_0}\log U_j \left[\Psi_j^0(x)
+ \frac{N}{p+1} \frac{a_{N,p}}{\ga_N} \mu_j^{\frac{N}{p+1}} |\log\mu_j| H(x,\xi_j) + o\(\mu_j^{\frac{N}{p+1}}|\log\mu_j|\)\right] dx \nonumber \\
&\ + o\(\mu_1^{N-2}\) \label{eq:r283} \\
&= (q_0-q_{\e})(A_3+o(1)) + o\(\mu_1^{N-2}|\log\mu_1|\). \nonumber
\end{align}

\medskip \noindent \textsc{Conclusion.} As before, we write $d_i = d_{i\e} = \mu_{i\e} \mu_{1\e}^{-1} \in (0,\infty)$ for $i = 1, \ldots, k$ and $\e \in (0,\e_0)$.
From \eqref{eq:energy'}--\eqref{eq:energy'2}, \eqref{eq:r221}, \eqref{eq:r132}, \eqref{eq:r28}--\eqref{eq:r281}, and \eqref{eq:r282}--\eqref{eq:r283}, we see
\begin{equation}\label{eq:r29}
2d_j^{N-2} (\nabla_{\xi} H)(\xi_j,\xi_j) - \sum_{i=1,\, i \ne j}^k \(d_i^{\frac{N}{p+1}} d_j^{\frac{N}{q_0+1}}
+ d_i^{\frac{N}{q_0+1}} d_j^{\frac{N}{p+1}}\) (\nabla_{\xi} G)(\xi_i,\xi_j) = o(1)
\end{equation}
and
\begin{multline*}
(N-2) d_j^{N-2} H(\xi_j, \xi_j) - N \sum_{i=1,\, i \ne j}^k \(\frac{d_i^{\frac{N}{p+1}} d_j^{\frac{N}{q_0+1}}}{q_0+1}
+ \frac{d_i^{\frac{N}{q_0+1}} d_j^{\frac{N}{p+1}}}{p+1}\) G(\xi_i,\xi_j) \\
= \frac{(A_3+o(1)) \ga_N}{a_{N,p}A_1} \frac{q_0-q_{\e}}{\mu_1^{N-2}|\log\mu_1|} + o(1)
\end{multline*}
for $j = 1, \ldots, k$.

\subsection{The case that $p \in [\frac{N-1}{N-2}, \frac{N}{N-2}) \cap (\frac{3}{N-2}, \frac{N}{N-2})$}\label{subsec:loc3}
We start with the following lemma that will be used in this and the next subsections.
Unlike the previous subsections, we need an auxiliary parameter $\kappa > 0$ to make the `remainder' terms truly small.
\begin{lemma}\label{lemma:loc2}
Let $N \ge 3$, $p \in (\frac{2}{N-2}, \frac{N}{N-2})$, $a_0 = ((N-2)p-2)(q_0-1) > 4$, and $\kappa \in (0,1)$ be a sufficiently small number. Then there exists a constant $C > 0$ independent of $\e \in (0,\e_0)$ such that
\begin{align*}
\mu_i^{\frac{2pN}{q_0+1}} \int_{B(x_i,\mu_1^{\kappa})} \(U_i^{q_{\e}-1} + U_i^{q_{\e}-2} \big|\Psi_i^0\big| \mathbf{1}_{q_{\e} \ge 2}\)
&\le C\(\mu_1^{2((N-2)p-2)} + \mu_1^{pN} \mu_1^{\kappa(N-a_0)}|\log\mu_1|\), \\
\mu_i^{\frac{2pN}{q_0+1}} \int_{B(x_i,\mu_1^{\kappa})} U_i^{q_{\e}-2} \big|\Psi_i^l\big| \mathbf{1}_{q_{\e} \ge 2}
&\le C\(\mu_1^{2((N-2)p-2)-1} + \mu_1^{pN} \mu_1^{\kappa(N-1-a_0)}|\log\mu_1|\), \\
\mu_i^{\frac{N}{q_0+1}} \int_{B(x_i,\mu_1^{\kappa})} V_i^p &\le C\(\mu_1^{(N-2)p-2}\mu_1^{\kappa(N-(N-2)p)}\) = o\(\mu_1^{(N-2)p-2}\),\\
\int_{B(x_i,\mu_1^{\kappa})} V_i^{p-2} \phi_{\e}^2 \(\big|\Phi_i^0\big| + \big|\Phi_i^l\big|\)
&= O\(\mu_1^{2((N-2)p-2)-1}|\log\mu_1|\) + o\(\mu_1^{(N-2)p-2}\), \\
\mu_1^{\frac{N}{q_0+1}} \int_{B(x_i,\mu_1^{\kappa})} V_i^{p-2} |\phi_{\e}| \(\big|\Phi_i^0\big| + \big|\Phi_i^l\big|\) &= O\(\mu_1^{(N-2)p-2+(N-3)}\) + o\(\mu_1^{(N-2)p-2}\)
\end{align*}
for $i = 1, \ldots, k$ and $l = 1, \ldots, N$. In particular, all the integrals above are of order $o(\mu_1^{(N-2)p-2})$ provided $N \ge 4$ and $p > \frac{3}{N-2}$.
\end{lemma}
\begin{proof}
We can obtain the first, second, and third estimates as in Lemma \ref{lemma:loc1}.
To derive the fourth and fifth estimates, we utilize \eqref{eq:error21}. For instance, applying $pq_0 > 1$ and $\frac{N(p+1)}{q_0+1} = (N-2)p-2$ also, we compute
\begin{align*}
&\ \int_{B(x_i,\mu_1^{\kappa})} V_i^{p-2} \phi_{\e}^2 \big|\Phi_i^0\big| \\
&\le C \left[\mu_1^{\frac{2N}{p+1}} \(\mu_1^{\frac{2Npq_0}{q_0+1}} + \mu_1^{\frac{2N(p+2)}{q_0+1}} |\log \mu_1| + \mu_1^{\frac{2Npq_0}{q_0+1}-2\eta}\)
\int_{B(0,\mu_1^{\kappa-1})} \frac{dy}{(1+|y|)^{(N-2)(p-1)}} \right. \\
&\hspace{25pt} + \mu_1^{\frac{2N}{p+1}+\frac{2N(p+2)}{q_0+1}-2(N-2)} \int_{B(0,\mu_1^{\kappa-1})} \frac{dy}{(1+|y|)^{(N-2)(p+1)}} \\
&\hspace{25pt} \left. + \mathbf{1}_{a_0 \le N} \mu_1^{\frac{2N}{p+1} + \frac{2Npq_0}{q_0+1} + 2(2-a_0)}
\int_{B(0,\mu_1^{\kappa})} \frac{dy}{(1+|y|)^{(N-2)(p-1)+2(a_0-2)}} \right]\\
&= o\(\mu_1^{(N-2)p-2}\) + O\(\mu_1^{2((N-2)p-2)}|\log\mu_1|\) + O\(\mu_1^{2((N-2)p-2)}\).
\end{align*}
The other estimates can be deduced similarly.
\end{proof}

\medskip \noindent \textsc{Estimate of $J_{3\e}$.}
In this estimate, we assume that $N \ge 4$ and $p \in (\max\{1,\frac{3}{N-2}\}, \frac{N}{N-2})$ so that it can be also applied in the next subsection.

\medskip
According to \eqref{eq:J3e}, \eqref{eq:PUdiff}, \eqref{eq:error21}, and Lemmas \ref{lemma:decay} and \ref{lemma:loc2}, we have
\begin{align*}
|J_{3\e}(\mux)| &\le \sum_{m=1}^k \int_{B(x_m,\mu_1^{\kappa})} \left|(PU_{\mux} + \psi_{\e})^{q_{\e}} - (PU_{\mux})^{q_{\e}} - q_0 U_j^{q_0-1}\psi_{\e}\right| \big|\Psi_j^l\big| \\
&\ + \sum_{m=1}^k \int_{B(x_m,\mu_1^{\kappa})} \left|\(\sum_{i=1}^k PV_i + \phi_{\e}\)^p - \(\sum_{i=1}^k PV_i\)^p - p V_j^{p-1} \phi_{\e} \right| \big|\Phi_j^l\big| + o\(\mu_1^{(N-2)p-2}\)
\end{align*}
for $\kappa \in (0,1)$ small enough. By \eqref{eq:psph1}, it holds that
\begin{equation}\label{eq:r30}
\left|(PU_{\mux} + \psi_{\e})^{q_{\e}} - (PU_{\mux})^{q_{\e}} - q_0 U_j^{q_0-1}\psi_{\e}\right| \big|\Psi_j^l\big|
\le C\(\mu_1^{\frac{2Np}{q_0+1}} U_m^{q_{\e}-1} + \mu_1^{Np}\) \quad \text{in } B(x_m, \mu_1^{\kappa})
\end{equation}
for $m \ne j$, and
\begin{equation}\label{eq:r301}
\begin{aligned}
&\ \left|(PU_{\mux} + \psi_{\e})^{q_{\e}} - (PU_{\mux})^{q_{\e}} - q_0 U_j^{q_0-1}\psi_{\e}\right| \\
&\le C\left[\mu_1^{\frac{Npq_0}{q_0+1}} + \mu_1^{\frac{2Np}{q_0+1}} U_j^{q_{\e}-2} \mathbf{1}_{q_{\e} \ge 2} + \e \mu_1^{Np \over q_0+1} U_j^{q_0-1} \(1 + U_j^{O(\e)}|\log U_j|\)\right] \quad \text{in } B(x_j,\mu_1^{\kappa}).
\end{aligned}
\end{equation}
Thus we infer from Lemma \ref{lemma:loc2}, \eqref{eq:r30}--\eqref{eq:r301}, and \eqref{eq:10} that
\[\sum_{m=1}^k \int_{B(x_m,\mu_1^{\kappa})} \left|(PU_{\mux} + \psi_{\e})^{q_{\e}} - (PU_{\mux})^{q_{\e}} - q_0 U_j^{q_0-1}\psi_{\e}\right| \big|\Psi_j^l\big| = o\(\mu_1^{(N-2)p-2}\).\]
Moreover, since $p < \frac{N}{N-2} \le 2$, we know
\begin{equation}\label{eq:r302}
\left|\(\sum_{i=1}^k PV_i + \phi_{\e}\)^p - \(\sum_{i=1}^k PV_i\)^p - p V_j^{p-1} \phi_{\e} \right| \le C V_m^p \quad \text{in } B(x_m,\mu_1^{\kappa})
\end{equation}
for $m \ne j$, and
\begin{equation}\label{eq:r303}
\begin{aligned}
&\ \left|\(\sum_{i=1}^k PV_i + \phi_{\e}\)^p - \(\sum_{i=1}^k PV_i\)^p - p V_j^{p-1} \phi_{\e} \right| \\
&\le C \left[\(\sum_{i=1}^k PV_i\)^{p-2}\phi_{\e}^2 + \left|\(\sum_{i=1}^k PV_i\)^{p-1} - V_j^{p-1}\right| |\phi_{\e}|\right] \\
&\le C V_j^{p-2} \(\phi_{\e}^2 + \mu_1^{\frac{N}{q_0+1}} |\phi_{\e}|\)
\end{aligned}
\quad \text{in } B(x_j,\mu_1^{\kappa}).
\end{equation}
By Lemmas \ref{lemma:loc1} and \ref{lemma:loc2}, and \eqref{eq:r302}--\eqref{eq:r303},
\[\sum_{m=1}^k \int_{B(x_m,\mu_1^{\kappa})} \left|\(\sum_{i=1}^k PV_i + \phi_{\e}\)^p - \(\sum_{i=1}^k PV_i\)^p - p V_j^{p-1} \phi_{\e} \right| |\Phi_j^l| = o\(\mu_1^{(N-2)p-2}\).\]
Therefore
\begin{equation}\label{eq:r31}
J_{3\e}(\mux) = o\(\mu_1^{(N-2)p-2}\).
\end{equation}

\medskip \noindent \textsc{Estimate of $J_{1\e}$.}
In the remaining part of this subsection, we assume that $N \ge 4$ and $p \in [\frac{N-1}{N-2}, \frac{N}{N-2}) \cap (\frac{3}{N-2}, \frac{N}{N-2})$.

\medskip
Fix any $l = 0, \ldots, N$. Arguing as in Subsection \ref{subsec:loc1}, we obtain equalities analogous to \eqref{eq:r132} and \eqref{eq:r16}:
\begin{equation}\label{eq:r320}
J_{1\e}(\mux) = - \sum_{m=1}^k \int_{B(x_m,\mu_1^{\kappa})} \left[(PU_{\mux})^{q_0} - \sum_{i=1}^k U_i^{q_0}\right] P\Psi_j^l + o\(\mu_1^{(N-2)p-2}\).
\end{equation}
and
\begin{multline}\label{eq:r32}
\sum_{m=1}^k \int_{B(x_m,\mu_1^{\kappa})} \left[(PU_{\mux})^{q_0} - \sum_{i=1}^k U_i^{q_0}\right] P\Psi_j^l
= \frac{a_{N,p}}{\ga_N} A_1 \mu_j^{(N-2)p-2} \pa_l \whh(\xi_j,\xi_j) \\
- \frac{a_{N,p}}{\ga_N} A_1 \sum_{i=1,\, i \ne j}^k \mu_i^{\frac{Np}{q_0+1}} \mu_j^{\frac{N}{q_0+1}} \pa_l \whg(\xi_j,\xi_i)
+ q_0 \int_{B(x_j,\mu_1^{\kappa})} U_j^{q_0-1} \hps_{\mux} \Psi_j^l + o\(\mu_1^{(N-2)p-2}\)
\end{multline}
where $\hps_{\mux}$, $\whh$, and $\whg$ are the functions defined by \eqref{eq:hps}, \eqref{eq:whh1}--\eqref{eq:whh2}, and \eqref{eq:whg}, respectively.
Unfortunately, analogous computations to \eqref{eq:r17}--\eqref{eq:r18} only give a crude estimate on the integral in the right-hand side of \eqref{eq:r32}.
We will overcome this technical issue by performing a more delicate analysis. Define the regular part $\tps_{\mux}$ of $\hps_{\mux}$ by
\[\tps_{\mux}(x) = \int_{\mca} G(x,y) \left[\(\sum_{i=1}^k PV_i\)^p - \sum_{i=1}^k V_i^p + \frac{p b_{N,p}}{\ga_N} \mu_1^{\frac{N}{q_0+1}}
\sum_{i=1}^k A_{\mbd,\bsxi,i} V_i^{p-1} \right](y) dy \quad \text{for } x \in \Omega.\]
Here, $(\mux) = (\bsmu_{\e},\mbx_{\e}) = (\mu_{1\e}, \ldots, \mu_{k\e}, x_{1\e}, \ldots, x_{k\e}) \in (0,\infty)^k \times \Omega^k$,
\[\mbx_{\e} \to \bsxi = (\xi_1, \ldots, \xi_k) \in \Omega^k \quad \text{and} \quad \mbd = \mbd_{\e} := \mu_{1\e}^{-1} \bsmu_{\e} \to \bsde \in (0,\infty)^k \quad \text{as } \e \to 0,\]
$A_{\mbd,\bsxi,i}$ is the quantity in \eqref{eq:Adxi} (where $\bsde$ is replaced with $\mbd$), and
\[\mca = \mca_{\e} := \Omega \setminus \cup_{i=1}^k B(x_{i\e},\vrh_{i\e}) \quad \text{where} \quad \vrh_i = \vrh_{i\e} := \mu_{1\e}^{\kappa} + \rho_2|x_{i\e}-\xi_i|^{\frac{1}{N+1}},\]
$\kappa \in (0,\frac{1}{(N-2)p+1})$ is a number small enough, and $\rho_2 > 0$ is one in the proof of Proposition \ref{prop:isos}.
\begin{lemma}
Assume that $N \ge 4$ and $p \in [\frac{N-1}{N-2}, \frac{N}{N-2})$.
Recall the functions $\wth_{\mbd,\bsxi}$ in \eqref{eq:wthdx} (where $\bsde$ is replaced with $\mbd$) and $\ovh$ in \eqref{eq:ovh}. Then
\begin{equation}\label{eq:r33}
\begin{aligned}
\tps_{\mux}(x) &= \(\frac{b_{N,p}}{\ga_N}\)^p \left[-\mu_1^{\frac{Np}{q_0+1}} \wth_{\mbd,\bsxi}(x) + \frac{\tga_{N,p,1}}{\ga_N} \sum_{i=1}^k \mu_i^{\frac{Np}{q_0+1}}\whh(x,\xi_i) \right. \\
&\hspace{60pt} \left. - \tga_{N,p,2} \mu_1^{\frac{N}{q_0+1}} \sum_{i=1}^k \mu_i^{\frac{N(p-1)}{q_0+1}} A_{\mbd,\bsxi,i} \ovh(x,\xi_i) \right] + \mcr_3(x) \quad \text{for } x \in \Omega
\end{aligned}
\end{equation}
where $\mcr_3$ is a remainder term satisfying $|\mcr_3(x)| + |\nabla \mcr_3(x)| = o(\mu_1^{\frac{Np}{q_0+1}})$
and the values of the constants $b_{N,p}$, $\ga_N$, $\tga_{N,p,1}$, and $\tga_{N,p,2}$ are found in \eqref{eq:HV}, \eqref{eq:H}, and \eqref{eq:tga}.
\end{lemma}
\begin{proof}
By virtue of \eqref{eq:V10est} and \eqref{eq:PV_i}, we have
\[\left|V_i(x) - \mu_i^{\frac{N}{q_0+1}} \frac{b_{N,p}}{|x-\xi_i|^{N-2}}\right| \le C \mu_1^{\frac{N}{q_0+1}} \frac{\mu_i + |x_i-\xi_i|}{|x-\xi_i|^{N-1}}\]
and
\[\left|PV_i(x) - \mu_i^{\frac{N}{q_0+1}} \frac{b_{N,p}}{\ga_N} G(x,\xi_i)\right| \le C \mu_1^{\frac{N}{q_0+1}} \frac{\mu_i + |x_i-\xi_i|}{|x-\xi_i|^{N-1}} 
\]
for $x \in \Omega \setminus B(x_i, \vrh_i)$. Thus
\begin{equation}\label{eq:r331}
\begin{aligned}
&\ \left[\left(\sum_{i=1}^k PV_i \right)^p - \sum_{i=1}^k V_i^p + \frac{pb_{N,p}}{\ga_N} \mu_1^{\frac{N}{q_0+1}} \sum_{i=1}^k A_{\mbd,\bsxi,i} V_i^{p-1}\right](y) \\
&= \(\frac{b_{N,p}}{\ga_N}\)^p \mu_1^{\frac{Np}{q_0+1}} \left[\(\sum_{i=1}^k d_i^{\frac{N}{q_0+1}} G(y,\xi_i)\)^p - \sum_{i=1}^k \(\frac{d_i^{\frac{N}{q_0+1}} \ga_N}{|y-\xi_i|^{N-2}}\)^p\right] \\
&\ + p \(\frac{b_{N,p}}{\ga_N}\)^p \mu_1^{\frac{N}{q_0+1}} \sum_{i=1}^k
A_{\mbd,\bsxi,i} \frac{\mu_i^{\frac{N(p-1)}{q_0+1}} \ga_N^{p-1}}{|y-\xi_i|^{(N-2)(p-1)}} + \sum_{i=1}^k \mu_1^{\frac{Np}{q_0+1}} 
O\(\frac{\mu_i + |x_i-\xi_i|}{|y-\xi_i|^{(N-2)p+1}}\) 
\end{aligned}
\end{equation}
for $y \in \mca$ where $d_i = \mu_i \mu_1^{-1}$. Besides,
\begin{multline}\label{eq:r332}
\left|\(\sum_{j=1}^k d_j^{\frac{N}{q_0+1}} G(y,\xi_j)\)^p - \(\frac{d_i^{\frac{N}{q_0+1}} \ga_N}{|y-\xi_i|^{N-2}}\)^p + p A_{\mbd,\bsxi,i} \(\frac{d_i^{\frac{N}{q_0+1}} \ga_N}{|y-\xi_i|^{N-2}}\)^{p-1}\right| \\
\le \frac{C}{|y-\xi_i|^{(N-2)(p-1)-1}} \quad \text{for } y \in B(x_i,\rho_2).
\end{multline}
From \eqref{eq:r331}, \eqref{eq:r332}, and the dominated convergence theorem, it follows that
\begin{align*}
\tps_{\mux}(x) &= \(\frac{b_{N,p}}{\ga_N}\)^p \mu_1^{\frac{Np}{q_0+1}} \int_{\mca} G(x,y) \left[\(\sum_{i=1}^k d_i^{\frac{N}{q_0+1}} G(y,\xi_i)\)^p
- \sum_{i=1}^k \(\frac{d_i^{\frac{N}{q_0+1}}\ga_N}{|y-\xi_i|^{N-2}}\)^p \right. \\
&\hspace{140pt} \left. + p \sum_{i=1}^k A_{\mbd,\bsxi,i} \(\frac{d_i^{\frac{N}{q_0+1}} \ga_N}{|y-\xi_i|^{N-2}}\)^{p-1} \right]dy + \mcr_3'(x) \\
&= \(\frac{b_{N,p}}{\ga_N}\)^p \mu_1^{\frac{Np}{q_0+1}} \left[\wtg_{\mbd,\bsxi}(x) - \tga_{N,p,1} \sum_{i=1}^k \frac{d_i^{\frac{Np}{q_0+1}}}{|x-\xi_i|^{(N-2)p-2}}
+ \tga_{N,p,2} \sum_{i=1}^k A_{\mbd,\bsxi,i} \frac{d_i^{\frac{N(p-1)}{q_0+1}}}{|x-\xi_i|^{(N-2)p-N}} \right.\\
&\hspace{90pt} \left. + \frac{\tga_{N,p,1}}{\ga_N} \sum_{i=1}^k d_i^{\frac{Np}{q_0+1}} \whh(x,\xi_i)
- \tga_{N,p,2} \sum_{i=1}^k A_{\mbd,\bsxi,i} d_i^{\frac{N(p-1)}{q_0+1}} \ovh(x,\xi_i)\right] +\mcr_3(x)
\end{align*}
for $x \in \Omega$ provided $\kappa \in (0,\frac{1}{(N-2)p+1})$ small. Here, $\mcr_3'$ is a function such that $|\mcr_3'(x)|+|\nabla \mcr_3'(x)| = o(\mu_1^{\frac{Np}{q_0+1}})$.
This and \eqref{eq:wthdx} give us the desired equality \eqref{eq:r33}.
\end{proof}

By the definition of $\tps_{\mux}$, we have
\begin{align*}
\hps_{\mux}(x) &= \tps_{\mux}(x) - \frac{pb_{N,p}}{\ga_N} \mu_1^{\frac{N}{q_0+1}} \sum_{i=1}^k A_{\mbd,\bsxi,i} \int_{\mca} G(x,y) V_i^{p-1}(y) dy \\
&\ + \sum_{i=1}^k \int_{B(x_i,\vrh_i)} G(x,y) \left[\(\sum_{m=1}^k PV_m\)^p -\sum_{m=1}^k V_m^p\right](y)dy \quad \text{for } x \in \Omega.
\end{align*}
Hence
\begin{equation}\label{eq:r34}
\begin{aligned}
&\ q_0 \int_{B(x_j, \mu_1^{\kappa})} U_j^{q_0-1} \hps_{\mux} \Psi_j^l - q_0 \int_{B(x_j, \mu_1^{\kappa})} U_j^{q_0-1} \tps_{\mux} \Psi_j^l \\
&= \sum_{i=1}^k \int_{B(x_j, \mu_1^{\kappa})} \int_{B(x_i,\vrh_i)} \left[\(\sum_{m=1}^k PV_m\)^p -\sum_{m=1}^k V_m^p\right](y) G(y,x) \(q_0 U_j^{q_0-1} \Psi_j^l\)(x) dy dx\\
&\ - \frac{b_{N,p}}{\ga_N} \mu_1^{\frac{N}{q_0+1}} \sum_{i=1}^k A_{\mbd,\bsxi,i} \int_{B(x_j, \mu_1^{\kappa})} \int_{\mca} G(x,y) p V_i^{p-1}(y) \(q_0 U_j^{q_0-1} \Psi_j^l\)(x) dy dx\\
&=: \mci_1 - \frac{b_{N,p}}{\ga_N} \mu_1^{\frac{N}{q_0+1}} \sum_{i=1}^k A_{\mbd,\bsxi,i} \mci_{2i}.
\end{aligned}
\end{equation}
Applying Fubini's theorem, we compute
\begin{align}
\mci_1 &= \sum_{i=1}^k \int_{B(x_i,\vrh_i)} \left[\(\sum_{m=1}^k PV_m\)^p - \sum_{m=1}^k V_m^p\right](y) \left[P\Phi_j^l(y) + O\(\mu_1^{\frac{Npq_0}{q_0+1}}\mu_1^{-\kappa[(N-2)p-2]q_0-\kappa}\)\right] dy \nonumber \\
&= o\(\mu_1^{(N-2)p-2}\). \label{eq:r341}
\end{align}
Let us calculate $\mci_{2i}$. Because of \eqref{eq:V10est},
\[\left|V_i^{p-1}(y) - \mu_i^{\frac{N(p-1)}{q_0+1}} \frac{b_{N,p}^{p-1}}{|y-x_i|^{(N-2)(p-1)}}\right| \le C\mu_1^{\frac{N(p-1)}{q_0+1}+1} \frac{1}{|y-x_i|^{(N-2)(p-1)+1}} \quad \text{for } y \in \mca.\]
Accordingly,
\begin{equation}\label{eq:r35}
\begin{aligned}
\int_{\mca} G(x,y) pV_i^{p-1}(y) dy &= \mu_i^{\frac{N(p-1)}{q_0+1}} \int_{\mca} G(x,y) \frac{p b_{N,p}^{p-1}}{|y-x_i|^{(N-2)(p-1)}} dy + \mcr_4(x) \\
&= \(\frac{b_{N,p}}{\ga_N}\)^{p-1} \tga_{N,p,2} \mu_i^{\frac{N(p-1)}{q_0+1}} \(|x-x_i|^{N-(N-2)p} - \ovh(x,x_i)\) \\
&\ - \mu_i^{\frac{N(p-1)}{q_0+1}} \sum_{m=1}^k \int_{B(x_m,\vrh_m)} G(x,y) \frac{p b_{N,p}^{p-1}}{|y-x_i|^{(N-2)(p-1)}} dy + \mcr_4(x)
\end{aligned}
\end{equation}
where $\mcr_4$ is a function satisfying $|\mcr_4(x)|+|\nabla \mcr_4(x)| = o(\mu_1^{\frac{N(p-1)}{q_0+1}})$.
On the other hand, integration by parts and the oddness of $\Phi_{1,0}^l$ in the $l$-th variable yield
\begin{multline}\label{eq:r351}
\mu_i^{\frac{N(p-1)}{q_0+1}} \int_{B(x_j, \mu_1^{\kappa})} \(|x-x_i|^{N-(N-2)p} - \ovh(x,x_i)\) \(q_0 U_j^{q_0-1} \Psi_j^l\)(x) dx \\
= A_1 \mu_i^{\frac{N(p-1)}{q_0+1}} \mu_j^{\frac{N}{q_0+1}} \left[\pa_l \ovh(x_j,x_i) - (1-\delta_{ij}) \frac{(N-(N-2)p) (x_j-x_i)_l}{|x_j-x_i|^{(N-2)(p-1)}}\right] + o\(\mu_1^{\frac{Np}{q_0+1}}\)
\end{multline}
and
\begin{multline}\label{eq:r352}
\mu_i^{\frac{N(p-1)}{q_0+1}} \sum_{m=1}^k \int_{B(x_j,\mu_1^{\kappa})} \int_{B(x_m,\vrh_m)} \frac{G(x,y)}{|y-x_i|^{(N-2)(p-1)}} \(q_0 U_j^{q_0-1} \Psi_j^l\)(x) dydx \\
= \mu_i^{\frac{N(p-1)}{q_0+1}} \int_{B(x_j,\vrh_j)} \frac{1}{|y-x_i|^{(N-2)(p-1)}} \left[\Phi_j^l(y) + O\(\mu_1^{\frac{N}{q_0+1}}\)\right] dy
+ o\(\mu_1^{\frac{Np}{q_0+1}}\) = o\(\mu_1^{\frac{Np}{q_0+1}}\)
\end{multline}
where $A_1 > 0$ is the constant in \eqref{eq:A12} and $\delta_{ij}$ stands for the Kronecker delta.
By integrating \eqref{eq:r35} over $B(x_j, \mu_1^{\kappa})$, applying \eqref{eq:r351}, \eqref{eq:r352}, and the oddness of $\Psi_{1,0}^l$ in the $l$-th variable, we obtain
\begin{equation}\label{eq:r36}
\begin{aligned}
\mci_{2i} 
&= A_1 \(\frac{b_{N,p}}{\ga_N}\)^{p-1} \tga_{N,p,2} \mu_i^{\frac{N(p-1)}{q_0+1}} \mu_j^{\frac{N}{q_0+1}} \left[\pa_l \ovh(x_j,x_i) - (1-\delta_{ij}) \frac{(N-(N-2)p) (x_j-x_i)_l}{|x_j-x_i|^{(N-2)(p-1)}}\right] \\
&\ + o\(\mu_1^{\frac{Np}{q_0+1}}\)
\end{aligned}
\end{equation}
for $l = 1, \ldots, N$. Combining \eqref{eq:r32}, \eqref{eq:r34}, \eqref{eq:r341}, \eqref{eq:r36}, \eqref{eq:whh2}, the identities
\[b_{N,p} = \ga_N A_1 \quad \text{and} \quad \(\frac{b_{N,p}}{\ga_N}\)^p \tga_{N,p,1} = a_{N,p},\]
(see \eqref{eq:tga} and \eqref{eq:HVab} for the second identity), and the estimate
\[q_0 \int_{B(x_j, \mu_1^{\kappa})} U_j^{q_0-1} \tps_{\mux} \Psi_j^l = - A_1 \mu_j^{\frac{N}{q_0+1}} \pa_l \tps_{\mux}(\xi_j) + o\(\mu_1^{(N-2)p-2}\),\]
we establish
\begin{align}
&\ \sum_{m=1}^k \int_{B(x_m,\mu_1^{\kappa})} \left[(PU_{\mux})^{q_0} - \sum_{i=1}^k U_i^{q_0}\right] P\Psi_j^l \nonumber \\
&= \frac{a_{N,p}}{\ga_N} A_1 \mu_j^{(N-2)p-2} \pa_l \whh(\xi_j, \xi_j) - \frac{a_{N,p}}{\ga_N} A_1
\sum_{i=1,\, i \ne j}^k \mu_i^{\frac{Np}{q_0+1}} \mu_j^{\frac{N}{q_0+1}} \pa_l \whg(\xi_j, \xi_i) - A_1 \mu_j^{\frac{N}{q_0+1}} \pa_l \tps_{\mux}(\xi_j) \nonumber \\
&\ - A_1^{p+1} \tga_{N,p,2} \mu_1^{\frac{N}{q_0+1}} \sum_{i=1}^k \mu_i^{\frac{N(p-1)}{q_0+1}} \mu_j^{\frac{N}{q_0+1}} A_{\mbd,\bsxi,i} \pa_l \ovh(x_j,x_i) \label{eq:r381} \\
&\ + A_1^{p+1} \tga_{N,p,2} \mu_1^{\frac{N}{q_0+1}} \sum_{i=1,\, i \ne j}^k \mu_i^{\frac{N(p-1)}{q_0+1}} \mu_j^{\frac{N}{q_0+1}}
A_{\mbd,\bsxi,i} \frac{(N-(N-2)p) (x_j-x_i)_l}{|x_j-x_i|^{(N-2)(p-1)}} + o\(\mu_1^{(N-2)p-2}\) \nonumber \\
&= A_1^{p+1} \mu_j^{\frac{N}{q_0+1}} \left[\mu_1^{\frac{Np}{q_0+1}} \pa_l \wth_{\mbd,\bsxi}(\xi_j)
+ ((N-2)p-2)\tga_{N,p,1} \sum_{i=1,\, i \ne j}^k \frac{\mu_i^{\frac{Np}{q_0+1}} (\xi_j-\xi_i)_l}{|\xi_j-\xi_i|^{(N-2)p}}\right. \nonumber \\
&\hspace{75pt} + \left. (N-(N-2)p) \tga_{N,p,2} \mu_1^{\frac{N}{q_0+1}} \sum_{i=1,\, i \ne j}^k A_{\mbd,\bsxi,i} \frac{\mu_i^{\frac{N(p-1)}{q_0+1}} (\xi_j-\xi_i)_l}{|\xi_j-\xi_i|^{(N-2)(p-1)}}\right] + o\(\mu_1^{(N-2)p-2}\) \nonumber
\end{align}
for $l = 1, \ldots, N$. Similarly, we discover
\begin{align}
&\ \sum_{m=1}^k \int_{B(x_m,\mu_1^{\kappa})} \left[(PU_{\mux})^{q_0} - \sum_{i=1}^k U_i^{q_0}\right] P\Psi_j^0 \nonumber \\
&= \frac{NA_1^{p+1}}{q_0+1} \mu_j^{\frac{N}{q_0+1}} \left[\mu_1^{\frac{Np}{q_0+1}} \wth_{\mbd,\bsxi}(\xi_j) - \sum_{i=1,\, i \ne j}^k \frac{\tga_{N,p,1} \mu_i^{\frac{Np}{q_0+1}}}{|\xi_j-\xi_i|^{(N-2)p-2}}
+ \mu_1^{\frac{N}{q_0+1}} \sum_{i=1,\, i \ne j}^k \frac{\tga_{N,p,2} A_{\mbd,\bsxi,i} \mu_i^{\frac{N(p-1)}{q_0+1}}}{|\xi_j-\xi_i|^{(N-2)p-N}}\right] \nonumber \\
&\ + o\(\mu_1^{(N-2)p-2}\). \label{eq:r382}
\end{align}

\medskip \noindent \textsc{Estimate of $J_{2\e}$.} Arguing as in \eqref{eq:r192} and using $(N-2)p-2 > 1$, we get
\begin{equation}\label{eq:r371}
\int_{\Omega} \left[(PU_{\mux})^{q_0} - (PU_{\mux})^{q_{\e}}\right] P\Psi_j^l = o\(\mu_1^{(N-2)p-2}\)
\end{equation} 
for $l = 1, \ldots, N$. Moreover, as in \eqref{eq:r193}, we have
\begin{equation}\label{eq:r372}
\int_{\Omega} \left[(PU_{\mux})^{q_0} - (PU_{\mux})^{q_{\e}}\right] P\Psi_j^0 = (q_0-q_{\e})(A_3+o(1)) + o\(\mu_1^{(N-2)p-2}\).
\end{equation}

\medskip \noindent \textsc{Conclusion.} From \eqref{eq:energy'}--\eqref{eq:energy'2}, \eqref{eq:r31}, \eqref{eq:r320}, \eqref{eq:r381}--\eqref{eq:r382}, and \eqref{eq:r371}--\eqref{eq:r372}, we see
\begin{multline}\label{eq:r391}
(\nabla \wth_{\mbd,\bsxi})(\xi_j) + ((N-2)p-2) \tga_{N,p,1} \sum_{i=1,\, i \ne j}^k d_i^{\frac{Np}{q_0+1}} \frac{\xi_j-\xi_i}{|\xi_j-\xi_i|^{(N-2)p}} \\
+ (N-(N-2)p) \tga_{N,p,2} \sum_{i=1,\, i \ne j}^k A_{\mbd,\bsxi,i} d_i^{\frac{N(p-1)}{q_0+1}} \frac{\xi_j-\xi_i}{|\xi_j-\xi_i|^{(N-2)(p-1)}} = o(1)
\end{multline}
and
\begin{multline}\label{eq:r392}
\frac{NA_1^{p+1}}{q_0+1} d_j^{\frac{N}{q_0+1}} \left[\wth_{\mbd,\bsxi}(\xi_j) - \tga_{N,p,1} \sum_{i=1,\, i \ne j}^k \frac{d_i^{\frac{Np}{q_0+1}}}{|\xi_j-\xi_i|^{(N-2)p-2}}
+ \tga_{N,p,2} \sum_{i=1,\, i \ne j}^k \frac{A_{\mbd,\bsxi,i} d_i^{\frac{N(p-1)}{q_0+1}}}{|\xi_j-\xi_i|^{(N-2)p-N}}\right] \\
= (A_3+o(1)) \frac{q_0-q_{\e}}{\mu_1^{(N-2)p-2}} + o(1)
\end{multline}
for $j = 1, \ldots, k$.

\subsection{The case that $p \in (\max\{1,\frac{3}{N-2}\}, \frac{N-1}{N-2})$}
Compared to the previous case $p \in [\frac{N-1}{N-2}, \frac{N}{N-2})$, the analysis here is more straightforward.
Arguing as in the previous subsection, we observe that \eqref{eq:r31}, \eqref{eq:r371}, and \eqref{eq:r372} are still true.
Furthermore, we have the equalities that are analogous to \eqref{eq:r381} and \eqref{eq:r382}:
\begin{multline*}
\int_{\Omega}\left[(PU_{\mux})^{q_0} - \sum_{i=1}^k U_i^{q_0} \right] P\Psi_j^l \\
= A_1^{p+1} \mu_j^{\frac{N}{q_0+1}} \left[\mu_1^{\frac{Np}{q_0+1}} \pa_l \wth_{\mbd,\bsxi}(\xi_j)
+ ((N-2)p-2) \tga_{N,p,1} \sum_{i=1,\, i \ne j}^k \frac{\mu_i^{\frac{Np}{q_0+1}}(\xi_j-\xi_i)_l }{|\xi_j-\xi_i|^{(N-2)p}}\right] + o\(\mu_1^{(N-2)p-2}\)
\end{multline*}
for $l=1, \ldots, N$, and
\begin{multline*}
\int_{\Omega} \left[(PU_{\mux})^{q_0} - \sum_{i=1}^k U_i^{q_0} \right] P\Psi_j^0 \\
= \frac{NA_1^{p+1}}{q_0+1} \mu_j^{\frac{N}{q_0+1}} \left[\mu_1^{\frac{Np}{q_0+1}} \wth_{\mbd,\bsxi}(\xi_j) - \sum_{i=1,\, i \ne j}^k
\frac{\tga_{N,p,1} \mu_i^{\frac{Np}{q_0+1}}}{|\xi_j-\xi_i|^{(N-2)p-2}}\right] + o\(\mu_1^{(N-2)p-2}\).
\end{multline*}
As a consequence,
\begin{equation}\label{eq:r393}
(\nabla \wth_{\mbd,\bsxi})(\xi_j) + ((N-2)p-2) \tga_{N,p,1} \sum_{i=1,\, i \ne j}^k d_i^{\frac{Np}{q_0+1}} \frac{\xi_j-\xi_i}{|\xi_j-\xi_i|^{(N-2)p}} = o(1)
\end{equation}
and
\begin{equation}\label{eq:r394}
\frac{NA_1^{p+1}}{q_0+1} d_j^{\frac{N}{q_0+1}} \left[\wth_{\mbd,\bsxi}(\xi_j)
- \tga_{N,p,1} \sum_{i=1,\, i \ne j}^k \frac{d_i^{\frac{Np}{q_0+1}}}{|\xi_j-\xi_i|^{(N-2)p-2}}\right]
= (A_3+o(1)) \frac{q_0-q_{\e}}{\mu_1^{(N-2)p-2}} + o(1)
\end{equation}
for $j = 1, \ldots, k$.

\subsection{Proof of Theorems \ref{thm:main2} and \ref{thm:main3}}
We are ready to establish Theorems \ref{thm:main2} and \ref{thm:main3}.

\begin{proof}[Proof of Theorem \ref{thm:main2}]
We assume that $N \ge 4$ and $p \in (\max\{1,\frac{3}{N-2}\}, \frac{N}{N-2})$.

\medskip \noindent (1) Fix any $j = 1, \ldots, k$ and $l = 1, \ldots, N$. By \eqref{eq:wthdx} and \eqref{eq:wthdxi},
\begin{equation}\label{eq:r40}
\begin{aligned}
&\ \wth_{\mbd_{\e},\bsxi}(x) \\
&= \begin{cases}
\displaystyle \wth_{\mbd_{\e},\bsxi,j}(x) + \sum_{i=1,\, i \ne j}^k \frac{\tga_{N,p,1} d_{i\e}^{\frac{Np}{q_0+1}}}{|x-\xi_i|^{(N-2)p-2}} &\text{if } p \in (\frac{2}{N-2}, \frac{N-1}{N-2}),\\
\displaystyle \wth_{\mbd_{\e},\bsxi,j}(x) + \sum_{i=1,\, i \ne j}^k \left[\frac{\tga_{N,p,1} d_{i\e}^{\frac{Np}{q_0+1}}}{|x-\xi_i|^{(N-2)p-2}}
- \frac{\tga_{N,p,2} A_{\mbd_{\e},\bsxi,i} d_{i\e}^{\frac{N(p-1)}{q_0+1}}}{|x-\xi_i|^{(N-2)p-N}}\right] &\text{if } p \in [\frac{N-1}{N-2}, \frac{N}{N-2})
\end{cases}
\end{aligned}
\end{equation}
for $x \in \Omega$. Therefore \eqref{eq:r391} and \eqref{eq:r393} are equivalent to $\pa_l \wth_{\mbd_{\e},\bsxi,j}(\xi_j) = o(1)$.
Taking $\e \to 0$ leads to $\pa_l \wth_{\dexi,j}(\xi_j) = 0$.

\medskip \noindent (2) By \eqref{eq:r392}, \eqref{eq:r394}, and \eqref{eq:r40},
\[\frac{NA_1^{p+1}}{q_0+1} d_{j\e}^{\frac{N}{q_0+1}} \wth_{\mbd_{\e},\bsxi,j}(\xi_j) = (A_3+o(1)) \frac{q_0-q_{\e}}{\mu_{1\e}^{(N-2)p-2}} + o(1).\]
Note that $\mu_{1\e} = u_{\e}^{-1/\alpha_{\e}}(x_{1\e})$, $\lim_{\e \to 0} \alpha_{\e} = \frac{N}{q_0+1}$, $\lim_{\e \to 0} \mu_{1\e}^{\e} = 1$, and
\begin{equation}\label{eq:q0qe}
q_0-q_{\e} = \frac{\e (q_0+1)^2}{N+\e (q_0+1)}.
\end{equation}
From these observations, \eqref{eq:S}, \eqref{eq:A12}, and \eqref{eq:r194}, we see
\[\lim_{\e \to 0} \e u_{\e}^{p+1}(x_{1\e}) = \frac{N^2A_1^{p+1}}{A_3(q_0+1)^3} \delta_j^{\frac{N}{q_0+1}} \wth_{\dexi,j}(\xi_j)
= \frac{N}{q_0+1} S^{\frac{p(q_0+1)}{1-pq_0}} \|U_{1,0}\|_{L^{q_0}(\R^N)}^{(p+1)q_0} \delta_j^{\frac{N}{q_0+1}} \wth_{\dexi,j}(\xi_j)\]
for $j = 1, \ldots, k$. From this, the assertions in the statement easily follow.

\medskip \noindent (3) Fixing a small number $\eta \in (0,\rho_2)$, we assume that $x \in \Omega \setminus \cup_{i=1}^k B(\xi_i,\eta)$. Then \eqref{eq:902} and \eqref{eq:outeru} imply
\[\lambda_{1\e}^{\alpha_{\e}} u_{\e}^{q_{\e}}(x) \le C \sum_{i=1}^k \frac{\lambda_{1\e}^{\alpha_{\e} (q_{\e}+1)} \lambda_{i\e}^{-q_{\e} ((N-2)p-2)}}{|x-x_{i\e}|^{q_{\e} ((N-2)p-2)}}
\le C \eta^{-q_0((N-2)p-2)} \lambda_{1\e}^{-2(p+1)+O(\e)} \to 0 \quad \text{as } \e \to 0.\]
Using this, \eqref{eq:uvconv1}, \eqref{eq:error13}, and the dominated convergence theorem, we conclude
\begin{align*}
\|u_{\e}\|_{L^\infty(\Omega)} v_{\e}(x) 
&= \sum_{i=1}^k G(x,\xi_i) \int_{B(x_{i\e},\frac{\eta}{2})} \lambda_{1\e}^{\alpha_{\e}} u_{\e}^{q_{\e}}(y) dy
+ \sum_{i=1}^k \int_{B(x_{i\e},\frac{\eta}{2})} [G(x,y)-G(x,\xi_i)] \lambda_{1\e}^{\alpha_{\e}} u_{\e}^{q_{\e}}(y) dy \\
&\ + \int_{\Omega \setminus \cup_{i=1}^k B(x_{i\e},\frac{\eta}{2})} G(x,y) \lambda_{1\e}^{\alpha_{\e}} u_{\e}^{q_{\e}}(y) dy \\
&\to \|U_{1,0}\|_{L^{q_0}(\R^N)}^{q_0} \sum_{i=1}^k \delta_i^{\frac{N}{q_0+1}}G(x,\xi_i) + 0 + 0
\end{align*}
in $C^1(\Omega \setminus \cup_{i=1}^k B(\xi_i,\eta))$ as $\e \to 0$.

Besides, because
\[-\Delta \(\lambda_{1\e}^{\alpha_{\e} p} u_{\e}\) = (\lambda_{1\e}^{\alpha_{\e}}v_{\e})^p = (\|u_{\e}\|_{L^\infty(\Omega)} v_{\e})^p
\to \|U_{1,0}\|_{L^{q_0}(\R^N)}^{pq_0} \(\sum_{i=1}^k \delta_i^{\frac{N}{q_0+1}} G(\cdot,\xi_i)\)^p,\]
we know
\[\|u_{\e}\|_{L^\infty(\Omega)}^p u_{\e}(x) = \lambda_{1\e}^{\alpha_{\e} p} u_{\e}(x) \to \|U_{1,0}\|_{L^{q_0}(\R^N)}^{pq_0} \wtg_{\dexi}(x)
\quad \text{in } C^1_{\textnormal{loc}}(\Omega \setminus \{\xi_1, \ldots, \xi_k\})\]
as $\e \to 0$. This completes the proof.
\end{proof}

\begin{proof}[Proof of Theorem \ref{thm:main3}]
In light of \cite[Theorems 1.1-1.3]{G}, it suffices to show that $\{(u_{\e}, v_{\e})\}_{\e \in (0,\e_0)}$ has only one blow-up point. Recall from \eqref{eq:r195} and \eqref{eq:r29} that
\begin{equation}\label{eq:r41}
\delta_j^{N-2} \nabla \tau(\xi_j) - \sum_{i=1,\, i \ne j}^k \(\delta_i^{\frac{N}{p+1}} \delta_j^{\frac{N}{q_0+1}}
+ \delta_i^{\frac{N}{q_0+1}} \delta_j^{\frac{N}{p+1}}\) (\nabla_{\xi} G)(\xi_i,\xi_j) = 0
\end{equation}
for $j = 1, \ldots, k$, provided $p \in [\frac{N}{N-2}, \frac{N+2}{N-2})$.
According to \cite[Corollary 3.2]{CT}, the Robin function $\tau$ is strictly convex on bounded convex domains in $\R^N$ for $N \ge 3$.
Therefore, by slightly modifying the proof of \cite[Theorem 3.2]{GT},
we conclude that there does not exist $\bsxi = (\xi_1, \ldots, \xi_k)$ satisfying \eqref{eq:r41} for $k \ge 2$.
\end{proof}

\section{Existence results}\label{sec:existence}
Let $X_{p,q_0} = W^{2,(p+1)/p}(\Omega) \times W^{2,(q_0+1)/q_0}(\Omega)$.
Given $\trh \in (0,1)$ small, we define the configuration set
\begin{multline}\label{eq:Lambda}
\Lambda_{\trh} = \left\{(\tdx) = (\td_1, \ldots, \td_k, x_1, \ldots, x_k) \in (\trh, \trh^{-1})^k \times \Omega^k: \right. \\
\left. \text{dist}(x_i, \pa\Omega) \ge \trh \text{ for } 1 \le i \le k,\,
\text{dist}(x_i, x_j) \ge \trh \text{ for } 1 \le i \ne j \le k \right\}.
\end{multline}
We write $\bsmu_{\e} = (\mu_{1\e}, \ldots, \mu_{k\e})$ where $\mu_{i\e}$ is the parameter defined by \eqref{eq:muie}.
Using this $\bsmu_{\e}$ and $\mbx$ in \eqref{eq:Lambda}, we also set $(U_{i\e}, V_{i\e})$ and $(\Psi_{i\e}^l,\Phi_{i\e}^l)$ by \eqref{eq:muUV} (replacing $x_i$ for $x_{i\e}$), and $PU_{\bsmu_{\e},\mbx}$ by \eqref{eq:PU}.
As before, we will often drop the subscript $\e$.

\medskip
The argument in \cite{KP}, based on the classical Lyapunov-Schmidt reduction, gives the following
\begin{prop}
Suppose that the assumptions in Theorem \ref{thm:mainby} hold.
Then, for each $\e \in (0,\e_0)$ and $(\tdx) \in \Lambda_{\trh}$, there exist a unique pair $(\psi_{\mux}, \phi_{\mux}) \in X_{p,q_0}$
and numbers $\{c_{il}\}_{i=1, \ldots, k,\, l=0, \ldots, N} \subset \R^{(N+1)k}$ such that
\begin{equation}\label{eq:LEaux}
\begin{cases}
\displaystyle -\Delta \(PU_{\mux} + \psi_{\mux}\) = \left|\sum_{i=1}^k PV_i + \phi_{\mux}\right|^{p-1}\(\sum_{i=1}^k PV_i + \phi_{\mux}\)
+ p \sum_{i=1}^k \sum_{l=0}^N c_{il} V_i^{p-1} \Phi_i^l &\text{in } \Omega, \\
\displaystyle -\Delta \(\sum_{i=1}^k PV_i + \phi_{\mux}\) = \left|PU_{\mux} + \psi_{\mux}\right|^{q_{\e}-1}\(PU_{\mux} + \psi_{\mux}\)
+ q_0 \sum_{i=1}^k \sum_{l=0}^N c_{il} U_i^{q_0-1} \Psi_i^l &\text{in } \Omega,\\
\displaystyle \psi_{\mux} = \phi_{\mux} = 0 &\text{on } \pa\Omega
\end{cases}
\end{equation}
satisfying
\begin{equation}\label{eq:LEaux_1}
\int_{\Omega} \(pV_i^{p-1}\Phi_i^l \phi_{\mux} + q_0U_i^{q_0-1}\Psi_i^l \psi_{\mux}\) = 0 \quad \text{for } i=1, \ldots, k \text{ and } l=0, \ldots, N
\end{equation}
and
\begin{equation}\label{eq:LEaux_2}
\begin{aligned}
&\ \left\|\Delta \psi_{\mux}\right\|_{W^{2,{p+1 \over p}}(\Omega)} + \left\|\Delta \phi_{\mux}\right\|_{W^{2,{q_0+1 \over q_0}}(\Omega)} \\
&\le C \left\|\sum_{i=1}^k U_i^{q_0} - (PU_{\mux})^{q_{\e}}\right\|_{L^{q_0+1 \over q_0}(\Omega)} \\
&\le C \begin{cases}
\mu_1^{N-2} + \mu_1^{Nq_0 \over p+1} |\log \mu_1|^{q_0 \over q_0+1} + \e |\log \mu_1| &\text{if } p \in (\frac{N}{N-2}, \frac{N+2}{N-2}),\\
\mu_1^{N-2} |\log \mu_1| + \mu_1^{Nq_0 \over p+1} |\log \mu_1|^{q_0(q_0+2) \over q_0+1} + \e |\log \mu_1| &\text{if } p = \frac{N}{N-2},\\
\mu_1^{(N-2)p-2} + \mu_1^{Npq_0 \over q_0+1} |\log \mu_1|^{q_0 \over q_0+1} + \e |\log \mu_1| &\text{if } p \in (\frac{2}{N-2}, \frac{N}{N-2})
\end{cases}
\end{aligned}
\end{equation}
for some $C > 0$ depending only on $N$, $p$, $\Omega$, $\e_0$, and $\trh$.
\end{prop}
\begin{proof}
One can argue as in \cite[Proposition 4.6, Lemma 3.3]{KP}.
\end{proof}
\noindent Employing the above proposition, we will prove Theorem \ref{thm:mainby}. To find the desired solution of \eqref{eq:LE}, it is sufficient to find $(\tdx) \in \Lambda_{\trh}$ such that
\begin{equation}\label{eq:cjl}
c_{jl} = 0 \quad \text{for } j=1, \ldots, k \text{ and } l=0, \ldots, N.
\end{equation}
Owing to \eqref{eq:LEaux}, equation \eqref{eq:cjl} is equivalent to
\begin{equation}\label{eq:energy'3}
I_{\e}'\(PU_{\mux} + \psi_{\mux}, \sum_{i=1}^k PV_i + \phi_{\mux}\) \(P\Psi_j^l, P\Phi_j^l\) = 0 \quad \text{for } j=1, \ldots, k \text{ and } l=0, \ldots, N;
\end{equation}
cf. \eqref{eq:energy'}.\footnote{To guarantee the positivity of the solution of \eqref{eq:LE}, we work with the energy functional \eqref{eq:energy}
whose nonlinear terms $|v|^{p+1}$ and $|u|^{q+1}$ are substituted with $v_+^{p+1}$ and $u_+^{q+1}$, respectively.
Thanks to the maximum principle for cooperative elliptic systems in \cite{dFM}, all components of the solution are positive.}
In the following, we will express \eqref{eq:energy'3} in terms of the function $\Upsilon_{p,k}$ in \eqref{eq:Upsilon}.
As in the previous section, we will split the cases according to the value of $p$.

\subsection{The case that $p \in (\frac{N}{N-2}, \frac{N+2}{N-2})$}
A direct computation using \eqref{eq:LEaux}--\eqref{eq:LEaux_2}, H\"older's inequality, and the Sobolev inequality shows
\begin{multline}\label{eq:ext11}
I_{\e}'\(PU_{\mux} + \psi_{\mux}, \sum_{i=1}^k PV_i + \phi_{\mux}\) \(P\Psi_j^l, P\Phi_j^l\)
- I_{\e}'\(PU_{\mux}, \sum_{i=1}^k PV_i\) \(P\Psi_j^l,P\Phi_j^l\) \\
= \int_{\Omega} \left[(PU_{\mux})^{q_{\e}} - \sum_{i=1}^k U_i^{q_0}\right] \psi_{\mux} - 2 \int_{\Omega} \nabla \psi_{\mux} \cdot \nabla \phi_{\mux} = o\(\mu_1^{N-2}\)
\end{multline}
uniformly in $(\tdx) \in \Lambda_{\trh}$ as $\e \to 0$; cf. \eqref{eq:r131}.
By Estimates of $J_{1\e}$ and $J_{2\e}$ in Subsection \ref{subsec:loc1}, \eqref{eq:ext11}, \eqref{eq:q0qe}, \eqref{eq:S}, \eqref{eq:A12}, \eqref{eq:r194}, and \eqref{eq:agaA2}, equality \eqref{eq:energy'3} reads
\begin{equation}\label{eq:ext12}
\td_j^{N-2} \nabla \tau(x_j) - \sum_{i=1,\, i \ne j}^k \(\td_i^{\frac{N}{p+1}} \td_j^{\frac{N}{q_0+1}}
+ \td_i^{\frac{N}{q_0+1}} \td_j^{\frac{N}{p+1}}\) (\nabla_{\xi} G)(x_i,x_j) = o(1)
\end{equation}
and
\begin{equation}\label{eq:ext13}
(N-2) \td_j^{N-2} \tau(x_j) - N \sum_{i=1,\, i \ne j}^k \(\frac{\td_i^{\frac{N}{p+1}} \td_j^{\frac{N}{q_0+1}}}{q_0+1}
+ \frac{\td_i^{\frac{N}{q_0+1}} \td_j^{\frac{N}{p+1}}}{p+1}\) G(x_i,x_j) = \frac{1}{A_1A_2} S^{\frac{p(q_0+1)}{pq_0-1}} + o(1)
\end{equation}
for $j = 1, \ldots, k$; cf. \eqref{eq:r195} and \eqref{eq:r196}.
We set $C_0 = (A_1A_2)^{-1} S^{p(q_0+1)/(pq_0-1)} > 0$ in the definition of $\Upsilon_{p,k}$ in \eqref{eq:Upsilon}.
Then it is easy to see that \eqref{eq:ext12} and \eqref{eq:ext13} are equivalent to
\begin{equation}\label{eq:ext14}
\nabla_{x_j} \Upsilon_{p,k}(\tdx) + o(1) = 0 \quad \text{and} \quad \frac{\pa \Upsilon_{p,k}}{\pa \td_j}(\tdx) + o(1) = 0 \quad \text{for } j = 1, \ldots, k,
\end{equation}
respectively.

If $(\dexi) \in (0,\infty)^k \times \Omega^k$ is an isolated critical point of $\Upsilon_{p,k}$ for which \eqref{eq:deg} holds,
then properties of the Brouwer degree yield a solution $(\wtmbd_{\e},\mbx_{\e}) \in \Lambda$ of \eqref{eq:ext14} for $\e > 0$ small
such that $(\wtmbd_{\e},\mbx_{\e}) \to (\dexi)$ as $\e \to 0$. The proof is completed.

\subsection{The case that $p = \frac{N}{N-2}$}
Let $\mu_{\e} = \left[\frac{-(N-2)\e}{W_{-1}(-(N-2)\e)}\right]^{1 \over N-2}$, which is a small positive number for $\e > 0$ small. By the definition of the function $W_{-1}$, it holds that
\begin{align*}
-(N-2)\e \mu_{\e}^{2-N} = W_{-1}(-(N-2)\e) &\Rightarrow -(N-2)\e = -(N-2)\e \mu_{\e}^{2-N} e^{-(N-2)\e \mu_{\e}^{2-N}} \\
&\Leftrightarrow \mu_{\e} = e^{-\e \mu_{\e}^{2-N}} \Leftrightarrow \mu_{\e}^{N-2} \log \mu_{\e} = -\e.
\end{align*}
Using these relations and arguing as in the previous subsection, we see that \eqref{eq:energy'3} is equivalent to \eqref{eq:ext12} and
\begin{multline}\label{eq:ext15}
(N-2) \td_j^{N-2} \tau(x_j) - N \sum_{i=1,\, i \ne j}^k \(\frac{\td_i^{\frac{N}{p+1}} \td_j^{\frac{N}{q_0+1}}}{q_0+1}
+ \frac{\td_i^{\frac{N}{q_0+1}} \td_j^{\frac{N}{p+1}}}{p+1}\) G(x_i,x_j) \\
= \frac{(A_3+o(1)) \ga_N}{a_{N,p}A_1} \frac{q_0-q_{\e}}{\mu_1^{N-2} \log\mu_1^{-1}} + o(1) = \frac{1}{\left|\S^{N-1}\right| b_{N,p}^p A_1} S^{\frac{p(q_0+1)}{pq_0-1}} + o(1).
\end{multline}
If we set $C_0 = (|\S^{N-1}| b_{N,p}^p A_1)^{-1} S^{p(q_0+1)/(pq_0-1)} > 0$, then \eqref{eq:ext12} and \eqref{eq:ext15} are reduced to \eqref{eq:ext14}.
Thus the desired conclusion follows from the preceding discussion.

\subsection{The case that $p \in [\frac{N-1}{N-2}, \frac{N}{N-2}) \cap (\frac{3}{N-2}, \frac{N}{N-2})$}
Reasoning as in \eqref{eq:ext11}, we obtain
\begin{multline}\label{eq:ext31}
I_{\e}'\(PU_{\mux} + \psi_{\mux}, \sum_{i=1}^k PV_i + \phi_{\mux}\) \(P\Psi_j^l, P\Phi_j^l\)
- I_{\e}'\(PU_{\mux}, \sum_{i=1}^k PV_i\) \(P\Psi_j^l,P\Phi_j^l\) \\
= o\(\mu_1^{(N-2)p-2}\)
\end{multline}
uniformly in $(\tdx) \in \Lambda_{\trh}$ as $\e \to 0$; cf. \eqref{eq:r31}.
By Estimates of $J_{1\e}$ and $J_{2\e}$ in Subsection \ref{subsec:loc3}, \eqref{eq:ext31}, \eqref{eq:q0qe}, \eqref{eq:S}, \eqref{eq:r194}, and the identity
\[\begin{cases}
\td_1^{\frac{Np}{q_0+1}} \wth_{\mbd,\mbx}(x) = \wth_{\tdx}(x),\\
\td_1^{\frac{N}{q_0+1}} A_{\mbd,\mbx,i} = A_{\tdx,i}
\end{cases}
\quad \text{for } \wtmbd = \td_1 \mbd = (\td_1d_1, \ldots, \td_1d_k),\ x \in \Omega,\ i = 1, \ldots, k,\]
equality \eqref{eq:energy'3} reads
\begin{multline}\label{eq:ext32}
\nabla \wth_{\tdx}(x_j) + ((N-2)p-2) \tga_{N,p,1} \sum_{i=1,\, i \ne j}^k \td_i^{\frac{Np}{q_0+1}} \frac{x_j-x_i}{|x_j-x_i|^{(N-2)p}} \\
+ (N-(N-2)p) \tga_{N,p,2} \sum_{i=1,\, i \ne j}^k A_{\tdx,i} \td_i^{\frac{N(p-1)}{q_0+1}} \frac{x_j-x_i}{|x_j-x_i|^{(N-2)(p-1)}} = o(1)
\end{multline}
and
\begin{multline}\label{eq:ext33}
\frac{N}{q_0+1} \td_j^{\frac{N}{q_0+1}} \left[\wth_{\tdx}(x_j) - \sum_{i=1,\, i \ne j}^k \frac{\tga_{N,p,1} \td_i^{\frac{Np}{q_0+1}}}{|x_j-x_i|^{(N-2)p-2}} \right. \\
\left. + \sum_{i=1,\, i \ne j}^k \frac{\tga_{N,p,2} A_{\tdx,i} \td_i^{\frac{N(p-1)}{q_0+1}}}{|x_j-x_i|^{(N-2)p-N}}\right] - \frac{S^{\frac{p(q_0+1)}{pq_0-1}}}{A_1^{p+1}} = o(1)
\end{multline}
for $j = 1, \ldots, k$; cf. \eqref{eq:r391} and \eqref{eq:r392}. To conclude the proof, it is enough to prove the following.
\begin{lemma}
Let $\textnormal{I}$ and $\textnormal{II}$ be the left-hand sides of \eqref{eq:ext32} and \eqref{eq:ext33}, respectively.
We set $C_0 = (p+1) A_1^{-(p+1)} S^{p(q_0+1)/(pq_0-1)} > 0$ in the definition of $\Upsilon_{p,k}$ in \eqref{eq:Upsilon}. For $j = 1, \ldots, k$, we have
\begin{equation}\label{eq:ext331}
\nabla_{x_j} \Upsilon_{p,k}(\tdx) = (p+1) \td_j^{N \over q_0+1} \textnormal{I}
\quad \text{and} \quad
\td_j \frac{\pa \Upsilon_{p,k}}{\pa \td_j}(\tdx) = (p+1) \textnormal{II}.
\end{equation}
\end{lemma}
\begin{proof}
(1) We set a function
\[P(x,\xi) = -(\nabla_{\xi} G)(x,\xi) \cdot \nu(\xi) \quad \text{for } x \in \Omega \text{ and } \xi \in \pa \Omega,\]
and fix $j = 1, \ldots, k$. By \eqref{eq:wthdxi} and \eqref{eq:wtgdx},
\begin{equation}\label{eq:ext341}
\begin{aligned}
\wth_{\tdx,j}(x) &= \int_{\Omega} G(x,z) \left[\(\frac{\ga_N \td_j^{N \over q_0+1}}{|z-x_j|^{N-2}}\)^p - \(\sum_{m=1}^k \td_m^{N \over q_0+1} G(z,x_m)\)^p\right] dz \\
&\ + \int_{\pa \Omega} P(x,z) \frac{\tga_{N,p,1} \td_j^{Np \over q_0+1}}{|z-x_j|^{(N-2)p-2}} dS_z - \tga_{N,p,2} A_{\tdx,j} \td_j^{N(p-1) \over q_0+1} |x-x_j|^{N-(N-2)p}
\end{aligned}
\end{equation}
for $x \in \Omega$, where the subscript $z$ in $dS_z$ refers to the variable of integration. Also, the representation formula tells us that
\begin{multline}\label{eq:ext342}
\tga_{N,p,2} |x-x_j|^{N-(N-2)p} = p \int_{\Omega} G(x,z) \(\frac{\ga_N}{|z-x_j|^{N-2}}\)^{p-1} dz \\
+ \tga_{N,p,2} \int_{\pa \Omega} P(x,z) |z-x_j|^{N-(N-2)p} dS_z.
\end{multline}
Plugging \eqref{eq:ext342} into \eqref{eq:ext341} and differentiating the result in the $x$-variable, we see
\begin{align}
\td_j^{N \over q_0+1} \(\nabla_x \wth_{\tdx,j}\)(x_j)
&= \td_j^{N \over q_0+1} \int_{\Omega} (\nabla_x G)(x_j,z) \left[\(\frac{\ga_N \td_j^{N \over q_0+1}}{|z-x_j|^{N-2}}\)^p - \(\sum_{m=1}^k \td_m^{N \over q_0+1} G(z,x_m)\)^p \right. \nonumber \\
&\hspace{120pt} \left. - p A_{\tdx,j} \(\frac{\ga_N \td_j^{N \over q_0+1}}{|z-x_j|^{N-2}}\)^{p-1} \right] dz \nonumber \\
&\ + \tga_{N,p,1} \int_{\pa \Omega} (\nabla_x P)(x_j,z) \frac{\td_j^{N(p+1) \over q_0+1}}{|z-x_j|^{(N-2)p-2}} dS_z \label{eq:ext34} \\
&\ - \tga_{N,p,2} A_{\tdx,j} \td_j^{Np \over q_0+1} \int_{\pa \Omega} (\nabla_x P)(x_j,z) |z-x_j|^{N-(N-2)p} dS_z. \nonumber
\end{align}
Note that the first term on the right-hand side of \eqref{eq:ext34} is well-defined for $(N-2)p < N$, because its integrand is bounded by
\[\frac{C}{|x_j-z|^{N-1}}\(\frac{1}{|z-x_j|^{(N-2)(p-2)}} + \frac{1}{|z-x_j|^{(N-2)(p-1)-1}}\) \quad \text{for } z \in \Omega \text{ near } x_j.\]
Moreover, rewriting \eqref{eq:ext341} as
\begin{align}
\wth_{\tdx,i}(x) &= \int_{\Omega} G(x,z) \left[\sum_{m=1}^k \(\frac{\ga_N \td_m^{N \over q_0+1}}{|z-x_m|^{N-2}}\)^p - \(\sum_{m=1}^k \td_m^{N \over q_0+1} G(z,x_m)\)^p\right] dz \nonumber \\
&\ + \tga_{N,p,1} \sum_{m=1}^k \int_{\pa \Omega} P(x,z) \frac{\td_m^{Np \over q_0+1}}{|z-x_m|^{(N-2)p-2}} dS_z
- \tga_{N,p,2} A_{\tdx,i} \td_i^{N(p-1) \over q_0+1} |x-x_i|^{N-(N-2)p} \nonumber \\
&\ - \tga_{N,p,1} \sum_{m = 1,\, m \ne i}^k \frac{\td_m^{Np \over q_0+1}}{|x-x_m|^{(N-2)p-2}} \label{eq:ext350}
\end{align}
for $x \in \Omega$ and $i = 1, \ldots, k$, and then differentiating it with respect to the parameter $x_j$, we observe
\begin{equation}\label{eq:ext351}
\begin{aligned}
\(\nabla_{x_j} \wth_{\tdx,i}\)(x)
&= p \int_{\Omega} G(x,z) \left[(N-2)\ga_N^p \td_j^{Np \over q_0+1} \frac{z-x_j}{|z-x_j|^{(N-2)p+2}} \right. \\
&\hspace{75pt} \left. - \(\sum_{m=1}^k \td_m^{N \over q_0+1} G(z,x_m)\)^{p-1} \td_j^{N \over q_0+1} (\nabla_{\xi} G)(z,x_j) \right] dz \\
&\ + ((N-2)p-2) \tga_{N,p,1} \td_j^{Np \over q_0+1} \int_{\pa \Omega} P(x,z) \frac{z-x_j}{|z-x_j|^{(N-2)p}} dS_z \\
&\ - \tga_{N,p,2} \(\nabla_{x_j} A_{\tdx,i}\) \td_i^{N(p-1) \over q_0+1} |x-x_i|^{N-(N-2)p} \\
&\ - ((N-2)p-2) \tga_{N,p,1} \td_j^{Np \over q_0+1} \frac{x-x_j}{|x-x_j|^{(N-2)p}}
\end{aligned}
\end{equation}
for $i \ne j$ and $x \in \Omega$ near $x_i$. Putting \eqref{eq:ext342} into \eqref{eq:ext350} and differentiating the result, we also obtain
\begin{align}
\(\nabla_{x_j} \wth_{\tdx,j}\)(x)
&= p \int_{\Omega} G(x,z) \left[(N-2)\ga_N^p \td_j^{Np \over q_0+1} \frac{z-x_j}{|z-x_j|^{(N-2)p+2}} \right. \nonumber \\
&\hspace{75pt} - \(\sum_{m=1}^k \td_m^{N \over q_0+1} G(z,x_m)\)^{p-1} \td_j^{N \over q_0+1} (\nabla_{\xi} G)(z,x_j) \nonumber \\
&\hspace{75pt} \left. - (N-2)(p-1) \ga_N^{p-1} A_{\tdx,j} \td_j^{N(p-1) \over q_0+1} \frac{z-x_j}{|z-x_j|^{(N-2)(p-1)+2}} \right] dz \nonumber \\
&\ + ((N-2)p-2) \tga_{N,p,1} \td_j^{Np \over q_0+1} \int_{\pa \Omega} P(x,z) \frac{z-x_j}{|z-x_j|^{(N-2)p}} dS_z \label{eq:ext352} \\
&\ + (N-(N-2)p) \tga_{N,p,2} A_{\tdx,j} \td_j^{N(p-1) \over q_0+1} \int_{\pa \Omega} P(x,z) \frac{z-x_j}{|z-x_j|^{(N-2)(p-1)}} dS_z \nonumber \\
&\ - \tga_{N,p,2} \(\nabla_{x_j} A_{\tdx,j}\) \td_j^{N(p-1) \over q_0+1} |x-x_j|^{N-(N-2)p} \nonumber
\end{align}
for $x \in \Omega$ near $x_j$. Combining \eqref{eq:ext351} and \eqref{eq:ext352}, we find
\begin{align}
&\ \sum_{i=1}^k \td_i^{N \over q_0+1} \(\nabla_{x_j} \wth_{\tdx,i}\)(x_i) \nonumber \\
&= p \td_j^{N \over q_0+1} \int_{\Omega} G(x_j,z) \left[(N-2)\ga_N^p \td_j^{Np \over q_0+1} \frac{z-x_j}{|z-x_j|^{(N-2)p+2}} \right. \nonumber \\
&\hspace{105pt} - \(\sum_{m=1}^k \td_m^{N \over q_0+1} G(z,x_m)\)^{p-1} \td_j^{N \over q_0+1} (\nabla_x G)(x_j,z) \nonumber \\
&\hspace{105pt} \left. - (N-2)(p-1) \ga_N^{p-1} A_{\tdx,j} \td_j^{N(p-1) \over q_0+1} \frac{z-x_j}{|z-x_j|^{(N-2)(p-1)+2}} \right] dz \nonumber \\
&\ + ((N-2)p-2) \tga_{N,p,1} \td_j^{N(p+1) \over q_0+1} \int_{\pa \Omega} P(x_j,z) \frac{z-x_j}{|z-x_j|^{(N-2)p}} dS_z \label{eq:ext35} \\
&\ + (N-(N-2)p) \tga_{N,p,2} A_{\tdx,j} \td_j^{Np \over q_0+1} \int_{\pa \Omega} P(x_j,z) \frac{z-x_j}{|z-x_j|^{(N-2)(p-1)}} dS_z \nonumber \\
&\ + p \sum_{i=1,\, i \ne j}^k \td_i^{N \over q_0+1} \int_{\Omega} G(x_i,z) \left[(N-2)\ga_N^p \td_j^{Np \over q_0+1} \frac{z-x_j}{|z-x_j|^{(N-2)p+2}} \right. \nonumber \\
&\hspace{140pt} \left. - \(\sum_{m=1}^k \td_m^{N \over q_0+1} G(z,x_m)\)^{p-1} \td_j^{N \over q_0+1} (\nabla_x G)(x_j,z) \right] dz \nonumber \\
&\ + ((N-2)p-2) \tga_{N,p,1} \sum_{i=1,\, i \ne j}^k \td_i^{N \over q_0+1} \td_j^{Np \over q_0+1} \left[\int_{\pa \Omega} P(x_i,z) \frac{z-x_j}{|z-x_j|^{(N-2)p}} dS_z - \frac{x_i-x_j}{|x_i-x_j|^{(N-2)p}}\right]. \nonumber
\end{align}

On the other hand, \eqref{eq:Upsilon} and the argument in the proof of Theorem \ref{thm:main2} (1) show
\begin{multline}\label{eq:ext36}
\nabla_{x_j} \Upsilon_{p,k}(\tdx) - (p+1) \td_j^{N \over q_0+1} \textnormal{I} \\
= \left[\td_j^{N \over q_0+1} \(\nabla_x \wth_{\tdx,j}\)(x_j) + \sum_{i=1}^k \td_i^{N \over q_0+1} \(\nabla_{x_j} \wth_{\tdx,i}\)(x_i)\right]
- (p+1) \td_j^{N \over q_0+1} \(\nabla_x \wth_{\tdx,j}\)(x_j).
\end{multline}
Hence, if we verify that the right-hand side of \eqref{eq:ext36} vanishes, we will obtain the first equality of \eqref{eq:ext331}.
To achieve the aim, we will employ Green's second identity: For a smooth domain $D$ and $f,\, g \in C^2(\overline{D})$,
\begin{equation}\label{eq:Greenid}
\int_D (g \Delta f - f\Delta g) (z) dz = \int_{\pa D} \(\frac{\pa f}{\pa \nu} g - \frac{\pa g}{\pa \nu} f\)(z) dS_z.
\end{equation}
For simplicity, let us assume that $x_j = 0$ and write $\Omega_r = \Omega \setminus B(0,r)$ for $r > 0$ small. If
\[\begin{cases}
\displaystyle f_1(z) := G(0,z),\\
\displaystyle g_1(z) := ((N-2)p-2) \tga_{N,p,1} \td_j^{N(p+1) \over q_0+1} \frac{z}{|z|^{(N-2)p}} + (N-(N-2)p) \tga_{N,p,2} A_{\tdx,j} \td_j^{Np \over q_0+1} \frac{z}{|z|^{(N-2)(p-1)}}
\end{cases}\]
for $z \in \Omega_r$, then \eqref{eq:Greenid} with $(f,g) = (f_1,g_1)$ and $D = \Omega_r$ is reduced to
\begin{multline}\label{eq:ext371}
p \int_{\Omega_r} G(0,z) \left[(N-2)\ga_N^p \td_j^{N(p+1) \over q_0+1} \frac{z}{|z|^{(N-2)p+2}} - (N-2)(p-1) \ga_N^{p-1} A_{\tdx,j} \td_j^{Np \over q_0+1} \frac{z}{|z|^{(N-2)(p-1)+2}} \right] dz \\
= -\int_{\pa \Omega} P(0,z) g_1(z) dS_z + o_r(1)
\end{multline}
where $o_r(1) \to 0$ as $r \to 0$. Indeed, it is a consequence of \eqref{eq:Green}, the symmetry of $\frac{z}{|z|}$ on $\pa B(0,r)$ with respect to the origin, and the relations
\[\Delta f_1(z) = 0 \quad \text{for } z \in \Omega_r \quad \text{and} \quad \frac{\pa f_1}{\pa \nu}(z) = -P(0,z) \quad \text{for } z \in \pa \Omega.\]
Also, if
\[\begin{cases}
\displaystyle f_2(z) := G(x_i,z) \quad (\text{for any } i \ne j),\\
\displaystyle g_2(z) := ((N-2)p-2) \tga_{N,p,1} \frac{z}{|z|^{(N-2)p}}
\end{cases}\]
for $z \in \Omega_r \setminus B(x_i,s)$, then plugging $(f,g) = (f_2,g_2)$ and $D = \Omega_r \setminus B(x_i,s)$ into \eqref{eq:Greenid}, taking $s \to 0$, and using the estimates
\begin{align*}
\int_{\pa B(x_i,s)} \(\frac{\pa f_2}{\pa \nu} g_2\)(z) dS_z
&= ((N-2)p-2) \tga_{N,p,1} \int_{\pa B(x_i,s)} \frac{(N-2)\ga_N}{|z-x_i|^{N-1}} \frac{z}{|z|^{(N-2)p}} dS_z + o_s(1) \\
&= ((N-2)p-2) \tga_{N,p,1} \frac{x_i}{|x_i|^{(N-2)p}} + o_s(1) \quad \text{(by \eqref{eq:H})}
\end{align*}
give
\begin{multline}\label{eq:ext372}
p (N-2)\ga_N^p \int_{\Omega_r} G(x_i,z) \frac{z}{|z|^{(N-2)p+2}} dz \\
= ((N-2)p-2) \tga_{N,p,1} \left[- \int_{\pa \Omega} P(x_i,z) \frac{z}{|z|^{(N-2)p}} dS_z + \frac{x_i}{|x_i|^{(N-2)p}}\right].
\end{multline}
Finally, if
\[\begin{cases}
\displaystyle f_3(z) := (\nabla_x G)(0,z),\\
\displaystyle g_3(z) := \tga_{N,p,1} \td_j^{N(p+1) \over q_0+1} \frac{1}{|z|^{(N-2)p-2}} - \tga_{N,p,2} A_{\tdx,j} \td_j^{Np \over q_0+1} |z|^{N-(N-2)p}
\end{cases}\]
for $z \in \Omega_r$, then \eqref{eq:Greenid} with $(f,g) = (f_3,g_3)$ and $D = \Omega_r$ reads
\begin{multline}\label{eq:ext373}
\int_{\Omega} (\nabla_x G)(0,z) \left[\ga_N^p \td_j^{N(p+1) \over q_0+1} \frac{1}{|z|^{(N-2)p}} - p \ga_N^{p-1} A_{\tdx,j} \td_j^{Np \over q_0+1} \frac{1}{|z|^{(N-2)(p-1)}}\right] dz \\
= -\int_{\pa \Omega} (\nabla_x P)(0,z) g_3(z) dS_z + o(1).
\end{multline}
In view of \eqref{eq:ext36}, \eqref{eq:ext34}, \eqref{eq:ext35}, \eqref{eq:ext371}--\eqref{eq:ext373}, and the dominated convergence theorem, we have
\begin{align*}
&\ \nabla_{x_j} \Upsilon_{p,k}(\tdx) - (p+1) \td_j^{N \over q_0+1} \textnormal{I} \\
&= \lim_{r \to 0} \left[- p \td_j^{N \over q_0+1} \int_{\Omega_r} G(x_j,z) \(\sum_{m=1}^k \td_m^{N \over q_0+1} G(z,x_m)\)^{p-1} \td_j^{N \over q_0+1} (\nabla_x G)(x_j,z) dz + o_r(1) \right. \\
&\hspace{35pt} - p \sum_{i=1,\, i \ne j}^k \td_i^{N \over q_0+1} \int_{\Omega_r} G(x_i,z) \(\sum_{m=1}^k \td_m^{N \over q_0+1} G(z,x_m)\)^{p-1} \td_j^{N \over q_0+1} (\nabla_x G)(x_j,z) dz + o_r(1) \\
&\hspace{35pt} \left. + p \td_j^{N \over q_0+1} \int_{\Omega_r} (\nabla_x G)(x_j,z) \(\sum_{m=1}^k \td_m^{N \over q_0+1} G(z,x_m)\)^p dz + o_r(1) \right] \\
&= 0.
\end{align*}
As a result, the first equality in \eqref{eq:ext331} is true.

\medskip \noindent (2) By \eqref{eq:wthdxi}, \eqref{eq:wtgdx}, and the assumption that $(N-2)p < N$, we have
\begin{align*}
&\ \sum_{i=1}^k \td_i^{\frac{N}{q_0+1}} \wth_{\tdx,i}(x_i) \\
&= \sum_{i=1}^k \int_{\Omega} \left[\(\frac{\ga_N \td_i^{N \over q_0+1}}{|z-x_i|^{N-2}}\)^{p+1}
- \td_i^{N \over q_0+1} G(x_i,z) \(\sum_{m=1}^k \td_m^{N \over q_0+1} G(z,x_m)\)^p\right] dz \\
&\ + \ga_N^{p+1} \sum_{i=1}^k \int_{\R^N \setminus \Omega} \frac{\td_i^{N(p+1) \over q_0+1}}{|z-x_i|^{(N-2)(p+1)}} dz
- \tga_{N,p,2} \sum_{i=1}^k A_{\tdx,i} \td_i^{\frac{Np}{q_0+1}} \underbrace{|x_i-x_i|^{N-(N-2)p}}_{=0}
\end{align*}
where the first term on the right-hand side is well-defined. It follows from \eqref{eq:Upsilon}, \eqref{eq:wthdx}, and \eqref{eq:ext33} that
\begin{align*}
&\ \frac{q_0+1}{N(p+1)} \left[\td_j \frac{\pa \Upsilon_{p,k}}{\pa \td_j}(\tdx) - (p+1) \textnormal{II}\right] \\
&= \int_{\Omega} \left[\(\frac{\ga_N \td_j^{N \over q_0+1}}{|z-x_j|^{N-2}}\)^{p+1}
- \td_j^{\frac{N}{q_0+1}} G(z,x_j) \(\sum_{m=1}^k \td_m^{N \over q_0+1} G(z,x_m)\)^p\right] dz\\
&\ + \ga_N^{p+1} \int_{\R^N \setminus \Omega} \frac{\td_j^{N(p+1) \over q_0+1}}{|z-x_j|^{(N-2)(p+1)}} dz - \frac{q_0+1}{N(p+1)} C_0 \\
&\ - \int_{\Omega} \left[\ga_N^{p+1} \sum_{i=1}^k \frac{\td_j^{\frac{N}{q_0+1}}}{|z-x_j|^{N-2}} \frac{\td_i^{Np \over q_0+1}}{|z-x_i|^{(N-2)p}}
- \td_j^{\frac{N}{q_0+1}} G(z,x_j) \(\sum_{m=1}^k \td_m^{N \over q_0+1} G(z,x_m)\)^p\right] dz \\
&\ - \ga_N^{p+1} \sum_{i=1}^k \int_{\R^N \setminus \Omega} \frac{\td_j^{N \over q_0+1}}{|z-x_j|^{N-2}}\frac{\td_i^{Np \over q_0+1}}{|z-x_i|^{(N-2)p}} dz
+ \tga_{N,p,2} \td_j^{\frac{N}{q_0+1}} \sum_{i=1}^k A_{\tdx,i} \td_i^{\frac{N(p-1)}{q_0+1}} |x_j-x_i|^{N-(N-2)p} \\
&\ + \td_j^{\frac{N}{q_0+1}} \left[\tga_{N,p,1} \sum_{i=1,\, i \ne j}^k \frac{\td_i^{\frac{Np}{q_0+1}}}{|x_j-x_i|^{(N-2)p-2}}
- \tga_{N,p,2} \sum_{i=1,\, i \ne j}^k \frac{A_{\tdx,i} \td_i^{\frac{N(p-1)}{q_0+1}}}{|x_j-x_i|^{(N-2)p-N}}\right] + \frac{q_0+1}{N} \frac{S^{\frac{p(q_0+1)}{pq_0-1}}}{A_1^{p+1}} = 0.
\end{align*}
Consequently, the second equality in \eqref{eq:ext331} is true.
\end{proof}

\subsection{The case that $p \in (\max\{1,\frac{3}{N-2}\}, \frac{N-1}{N-2})$}
This case is easier to handle than the previous one, and the value of $C_0 > 0$ remains the same. We skip the details.

\bigskip \noindent \small{\textbf{Acknowledgement.} S. Kim was supported by Basic Science Research Program through the National
Research Foundation of Korea (NRF) funded by the Ministry of Education (NRF2020R1C1C1A01010133, NRF2020R1A4A3079066).}

\end{document}